\documentclass{article}

 \usepackage[dblblindworkshop, final]{tagds_2025}

 \workshoptitle{Topology, Algebra, and Geometry in Data Science}

\usepackage[utf8]{inputenc} 
\usepackage[T1]{fontenc}    
\usepackage{hyperref}       
\usepackage{url}            
\usepackage{booktabs}       
\usepackage{nicefrac}       
\usepackage{microtype}      
\usepackage{xcolor}         

\usepackage{amssymb,amsmath,amsthm}       %
\usepackage{amsfonts}      %

\newtheorem{theorem}{Theorem}
\newtheorem{lemma}{Lemma}
\newtheorem{proposition}{Proposition}
\newtheorem{corollary}{Corollary}
\newtheorem{example}{Example}

\theoremstyle{definition}

\theoremstyle{remark}
\newtheorem{remark}{Remark}

\newcommand{\R}{\mathbb{R}}

\newcommand{\supp}{\operatorname{supp}}
\newcommand{\aff}{\operatorname{aff}}
\newcommand{\conv}{\operatorname{conv}}
\DeclareMathOperator*{\argmax}{arg\,max}

\newcommand{\chat}{\hat c}
\newcommand{\dmax}{d_{\max}}

\newcommand{\simplex}{\Delta_{n-1}}

\usepackage{mathtools}     %

\usepackage{graphicx}
\usepackage{algorithmic,algorithm}

\usepackage{url}

\title{Peeling metric spaces of strict negative type}

\author{%
  Steve Huntsman\\ 
  \texttt{steve.huntsman@cynnovative.com} \\
}

\begin{document}

\maketitle

\begin{abstract}
We describe a unified and computationally tractable framework for finding outliers in, and maximum-diversity subsets of, finite metric spaces of strict negative type. Examples of such spaces include finite subsets of Euclidean space and finite subsets of a sphere without antipodal points. The latter accounts for state-of-the-art text embeddings, and we apply our framework in this context to sketch a hallucination mitigation strategy and separately to a class of path diversity optimization problems with a real-world example.
\end{abstract}

\section{\label{sec:Introduction}Introduction}

Many problems in data science and machine learning can be distilled to identifying outliers \cite{boukerche2020outlier}, anomalies \cite{chandola2009anomaly}, and diverse subsets \cite{leinster2016maximizing, mahabadi2023core}. A vast and almost totally disconnected literature that we shall not attempt to capture with additional references is devoted to these problems. 

This paper details a unified natural interpretation of, and framework for, solving these problems in a broad class of situations. Specifically, so-called \emph{strict negative type} finite metric spaces (including, but not limited to, finite subsets of Euclidean space) admit a natural notion of outliers or boundary elements that we call a \emph{peel} and that simultaneously maximizes a natural measure of diversity \cite{leinster2021entropy}. The notion of a peel involves no ambiguity (e.g., free parameters) and a peel can be approximated by solving a finite (and in practice, short) sequence of linear equations, as detailed in Algorithm \ref{alg:ScaleZeroArgMaxDiversity} below. (\S \ref{sec:patch} details how to efficiently compute peels exactly with certification, anytime approximation guarantees, and marginal extra effort.) We then detail the applicability of peels to mitigating hallucinations in large language models. We then discuss product metrics, with an eye towards computing outlying/diverse sequences or paths, including a detailed real-world example. A supplement contains appendices with proofs, discussions of extensions, and auxiliary experimental results.

\section{\label{sec:Preliminaries}Weightings, magnitude, and diversity}

A square matrix $Z \ge 0$ is a \emph{similarity matrix} if $\text{diag}(Z) > 0$. We are concerned with the class of similarity matrices of the form $Z = \exp[-td]$ where $(f[M])_{jk} := f(M_{jk})$, i.e., the exponential is componentwise, $t \in (0,\infty)$, and $d$ is a square matrix whose entries are in $[0,\infty]$ and satisfy the triangle inequality. In this paper we will always assume that $d$ is the matrix of an actual metric (so in particular, symmetric along with $Z$) on a finite space. 

We say that $d$ is \emph{negative type} if $x^T d x \le 0$ for $1^T x = 0$ and $x^T x = 1$ (equivalently to this last, $x \ne 0$). If the inequality is strict, we say that $d$ is \emph{strict negative type}: this entails that $Z$ is positive semidefinite for \emph{all} $t>0$. Important examples of negative type metrics on finite spaces are finite subsets of Euclidean space with the $L^1$ or $L^2$ distances, finite subsets of spheres with the geodesic distance, finite subsets of hyperbolic space, and ultrametrics (i.e., metrics satisfying $d(x,z) \le \max \{d(x,y),d(y,z)\}$) on finite spaces. However, not all of these are strict negative type: e.g., spheres with antipodal points are not strict negative type \cite{hjorth1998finite}. 

A \emph{weighting} $w$ is a solution to $Zw = 1$, where $1$ indicates a vector of all ones. If $Z$ has a weighting $w$, then its \emph{magnitude} is $\text{Mag}(Z) := \sum_j w_j$. If $d$ is negative type, then $Z$ is positive definite, so it has a unique weighting. It turns out that weightings are excellent scale-dependent boundary or outlier detectors in Euclidean space \cite{willerton2009heuristic,bunch2020practical,huntsman2023diversity}: in fact, behavior evocative of boundary detection applies more generally \cite{huntsman2023magnitude}. A technical explanation of the Euclidean boundary-detecting behavior draws on the notion of Bessel capacities \cite{meckes2015magnitude}.

\begin{example}
\label{ex:3PointSpace}
Consider $\{x_j\}_{j=1}^3 \subset \mathbb{R}^2$ with $d_{jk} := d(x_j,x_k)$ given by $d_{12} = d_{13} = 1 = d_{21} = d_{31}$ and $d_{23} = \delta = d_{32}$ with $\delta \ll 1$. It turns out that
\begin{equation}
w_1 = \frac{e^{(\delta+2)t}-2e^{(\delta+1)t}+e^{2t}}{e^{(\delta+2)t}-2e^{\delta t}+e^{2t}}; \quad w_2 = w_3 = \frac{e^{(\delta+2)t}-e^{(\delta+1)t}}{e^{(\delta+2)t}-2e^{\delta t}+e^{2t}}. \nonumber
\end{equation}
For $t \ll 1$, $w \approx (1/4,1/4,1/2)^T$; for $t \gg 1$, $w \approx (1,1,1)^T$, and it turns out that for $t \approx 10$, $w \approx (1/2,1/2,1)^T$: see Figure \ref{fig:3pointSpace}. I.e., the two nearby points have ``effective sizes'' near 1/4, then 1/2, then 1; meanwhile, the far point has effective size near 1/2, then 1, where it remains; the ``effective number of points'' goes from $\approx 1/4+1/4+1/2 = 1$, to $\approx 1/2 + 1/2 + 1 = 2$, to $\approx 1+1+1 = 3$.


\begin{figure}[h]
  \centering
  \includegraphics[trim = 40mm 110mm 40mm 105mm, clip, width=.75\textwidth,keepaspectratio]{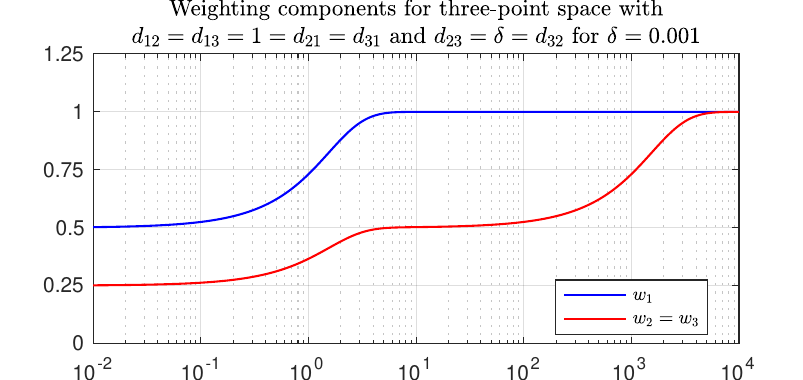}
  \caption{Weighting for an ``isoceles'' metric space. The magnitude function $w_1+w_2+w_3$ is a scale-dependent ``effective number of points.''}
  \label{fig:3pointSpace}
\end{figure}
\end{example}

Fairly recent mathematical developments have clarified the role that magnitude and weightings play in \emph{maximizing} a general and axiomatically supported notion of diversity \cite{leinster2016maximizing,leinster2021entropy}. 
Specifically, the \emph{diversity of order $q$} for a probability distribution $p$ and similarity matrix $Z$ is 
\begin{equation}
\label{eq:diversity}
\exp \left ( \frac{1}{1-q} \log \sum_{j: p_j > 0} p_j (Zp)_j^{q-1} \right )
\end{equation} 
for $1 < q < \infty$, and via limits for $q = 1,\infty$. This is a ``correct'' measure of diversity in much the same way that Shannon entropy is a ``correct'' measure of information. In fact, the logarithm of diversity is a geometrical generalization of the R\'enyi entropy of order $q$. The usual R\'enyi entropy is recovered for $Z = I$, and Shannon entropy subsequently for $q = 1$.

\begin{theorem}
\label{thm:maxDiversity}
If $Z$ is symmetric, positive definite, and has a unique positive weighting $w$, then for all $q$, $w$ is proportional to the diversity-maximizing distribution \cite{leinster2016maximizing}.
\end{theorem}

The situation described by Theorem \ref{thm:maxDiversity} reduces diversity maximization to a standard linear algebra problem while simultaneously removing any ambiguity regarding the parameter $q$. It is possible to efficiently compute a ``cutoff scale'' \cite{huntsman2023diversity} such that we can optimally enforce this desirable situation for similarity matrices of the form $Z = \exp[-td]$. However, in practice this scale is often quite large, and the resulting weighting will have many components with values close to unity, degrading the utility of this construction. It is frequently desirable to work in the limit $t \downarrow 0$: for example, in Figure \ref{fig:3pointSpace}, this limit successfully identifies one point as an outlier. We turn to this limit in the sequel.

\section{\label{sec:peelingTheorem}The peeling theorem}

For a probability distribution $p$ in $\Delta_{n-1} := \{ p \in [0,1]^n : 1^T p = 1 \}$, the diversity of order $1$ is
\begin{equation}
\label{eq:diversity1}
D_1^{Z}(p) := \prod_{j:p_j > 0} (Zp)_j^{-p_j}
\end{equation}
and the corresponding generalized entropy is
\begin{equation}
\label{eq:ssEntropy1}
\log D_1^{Z}(p) = - \sum_{j:p_j > 0} p_j \log (Zp)_j.
\end{equation}
These can be efficiently optimized for $Z = \exp[-td]$ in the limit $t \downarrow 0$ when $d$ is strict negative type.

The first-order approximation $Z = \exp[-td] \approx 11^T - td$ generically yields 
\begin{equation}
\label{eq:ssEntropy1Approx}
\log D_1^{Z}(p) \approx t p^T d p.
\end{equation}
The quantity $p^T d p$ is called the \emph{quadratic entropy} of $d$: it is convex if $d$ is strict negative type. (For details, see Theorem 4.3 of \cite{rao1984convexity} and Proposition 5.20 of \cite{devriendt2022graph} as well as \cite{leinster2012measuring,meckes2013positive,leinster2016maximizing,leinster2021entropy}.) Therefore if $d$ is strict negative type, \eqref{eq:ssEntropy1Approx} can be efficiently maximized over any sufficiently simple polytope via quadratic programming. However, in \S \ref{sec:peeling} we will give a more practical (i.e., much faster and more sparsity-accurate) algorithm for maximizing the quadratic entropy of strict negative type metrics.

\subsection{\label{sec:peeling}Maximizing quadratic entropy of strict negative type metrics}

We want to compute
\begin{equation}
\label{eq:argMaxQE}
p_*(d) := \arg \max_{p \in \Delta_{n-1}} p^T d p.
\end{equation}
For a strict negative type metric $d$, we call $p_*(d)$ (or, depending on context, its support) the \emph{peel} of $d$, for reasons that figures below make obvious.


The preceding version of this paper claimed that Algorithm \ref{alg:ScaleZeroArgMaxDiversity} below produced $p_*(d)$ exactly. This is not generally true: exhaustive details and remedies are provided in \S \ref{sec:patch}. In the interest of continuity with the preceding version, we make informal statements here that are also fully detailed in \S \ref{sec:patch}. The existence of a weakly polynomial algorithm for peeling (as a special case of convex quadratic programming) is ensured by (e.g.) \cite{monteiro1989interior}.

\begin{theorem}[peeling theorem]
\label{thm:peeling}
For $d$ strict negative type and metric, Algorithm \ref{alg:ScaleZeroArgMaxDiversity} returns an approximation to $p_*(d)$ with strong computable \emph{a priori} error bounds (and in practice, very good if not exact results) in time $O(n^{\omega+1})$, where $\omega \le 3$ is the exponent characterizing the complexity of matrix multiplication and inversion. Algorithm \ref{alg:ScaleZeroArgMaxDiversity} can be augmented with an efficient (in practice, very brief) ascent phase that terminates by producing the exact result in finitely many steps while also providing an anytime convergence guarantee.  \qed
\end{theorem}

\begin{algorithm}
  \caption{\textsc{ScaleZeroArgMaxDiversity}$(d)$}
  \label{alg:ScaleZeroArgMaxDiversity}
\begin{algorithmic}[1]
  \REQUIRE Strict negative type metric $d$ on $[n] \equiv \{1, \dots, n\}$
  \STATE $p \leftarrow \frac{d^{-1} 1}{1^T d^{-1} 1}$
  \WHILE{$\exists i : p_i < 0$}
  \STATE $\mathcal{J} \leftarrow \{ j : p_j > 0 \}$
  \COMMENT{Restriction of support}
  \STATE $p \leftarrow 0_{[n]}$
  \STATE $p_\mathcal{J} \leftarrow \frac{d_{\mathcal{J},\mathcal{J}}^{-1} 1_\mathcal{J}}{1_\mathcal{J}^T d_{\mathcal{J},\mathcal{J}}^{-1} 1_\mathcal{J}}$ 
  \ENDWHILE
  \RETURN $p = p_*(d)$
\end{algorithmic}
\end{algorithm}

\begin{corollary}
For $d$ strict negative type and metric and for all $q$, Algorithm \ref{alg:ScaleZeroArgMaxDiversity} efficiently approximates $\arg \max_{p \in \Delta_{n-1}} \lim_{t \downarrow 0} D_q^Z(p).$
\end{corollary}

\begin{figure}[h]
  \centering
  \includegraphics[trim = 65mm 110mm 65mm 100mm, clip, width=.4\columnwidth,keepaspectratio]{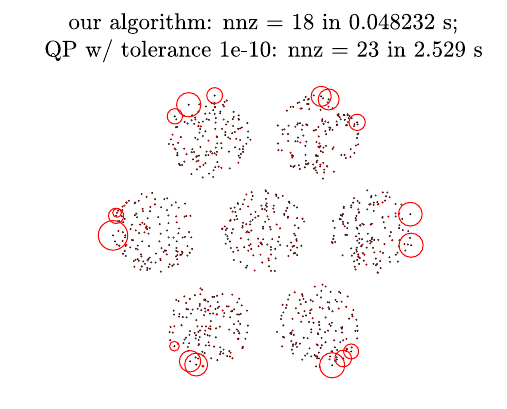}
  \includegraphics[trim = 65mm 110mm 65mm 100mm, clip, width=.4\columnwidth,keepaspectratio]{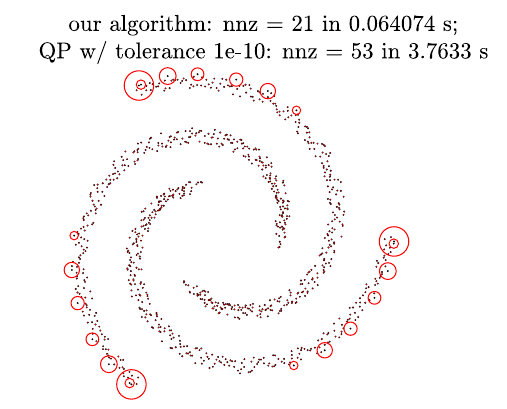}
\caption{Peels produced by Algorithm \ref{alg:ScaleZeroArgMaxDiversity} acting on the Euclidean distance matrix of the $\approx 1000$ black points, indicated by {\color{red}red} circles with radius proportional to the corresponding entries of $p$. The numbers of nonzero (nnz) entries of the output are indicated along with the runtimes of the algorithm; the same numbers are reported for a quadratic programming run with tolerance $10^{-10}$. 
  }
  \label{fig:maxQE}
\end{figure}

As a practical matter, Algorithm \ref{alg:ScaleZeroArgMaxDiversity} performs better than a quadratic programming solver: it is much faster (e.g., in MATLAB on $\approx 1000$ points, a few hundredths of a second versus several seconds for a quadratic programming solver with tolerance $10^{-10}$) and more accurate, in particular by handling sparsity exactly. Figure \ref{fig:maxQE} shows representative results.\footnote{We verified that these results are exact: this is typical for the approximation of Theorem \ref{thm:peeling} (see \S \ref{sec:patch}). The impact on runtime of an exact certified algorithm was negligible here, so we retain the original figure. We did \emph{not} attempt to redo the other examples below, as they are simply not amenable to this sort of verification.} It is also very simple to implement: excepting any preliminary checks on inputs, each line of the algorithm can be (somewhat wastefully) implemented in a standard-length line of MATLAB or Python. The exact algorithms in \S \ref{sec:patch} are also not complicated to implement, though some care should be paid to tolerances in predicate checks.

\subsection{\label{sec:peelingExamples}Iterated peeling of text embeddings}

We selected 150 named RGB color codes from the large-scale color survey \cite{xkcd,lo2024exploring} by restricting consideration to colors with a single word in their name, and then further restricting by human judgment to get a desired number while trying to avoid ambiguity. We then fed prompts of the form 
\begin{quote}
Describe the color of \underline{\hspace{1cm}} in relation to other colors.
\end{quote}
to \texttt{gpt-4o}, where the placeholder is for a color name. We embedded prompts and responses using \texttt{voyage-3.5} 
\footnote{See \url{https://blog.voyageai.com/2025/05/20/voyage-3-5/}. ModernBERT \cite{warner2024smarter} produced visually inferior embeddings (not shown, but see \cite{lo2024exploring}).} 
and repeatedly peeled the results using spherical distance of normalizations, as shown in Figures \ref{fig:red1}-\ref{fig:blueGreen2}. Appendix \S C in the supplement shows an example along the same lines with all 150 colors at once. 

All of our examples here and below were produced in seconds or less on a MacBook Pro.

\begin{figure}[htbp]
  \centering
  \includegraphics[trim = 0mm 0mm 0mm 0mm, clip, width=.3\columnwidth,keepaspectratio]{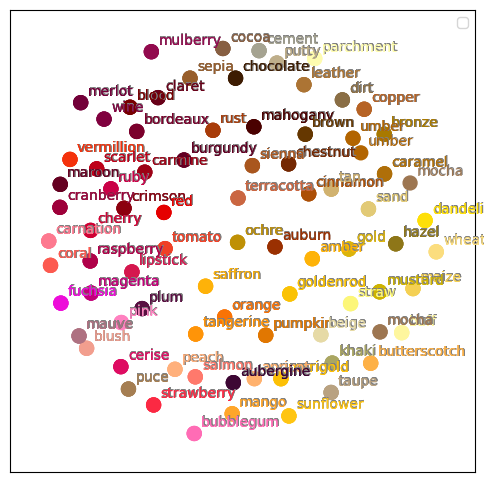}
  \includegraphics[trim = 0mm 0mm 0mm 0mm, clip, width=.3\columnwidth,keepaspectratio]{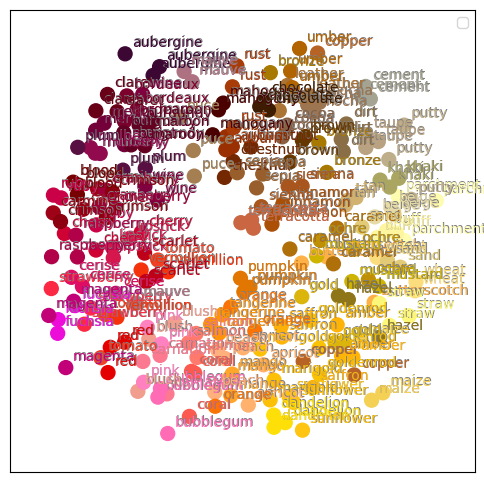}
  \includegraphics[trim = 0mm 0mm 0mm 0mm, clip, width=.3\columnwidth,keepaspectratio]{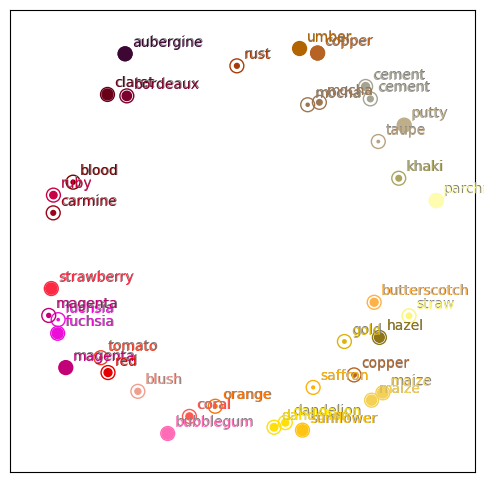}
\caption{Left: multidimensional scaling (MDS) of 3 prompt embeddings for each of the 80 predominantly red colors. Center: MDS of response embeddings. Since the same prompt yields different responses, $3 \cdot 80 = 240$ distinct points are shown. Right: The peel of response embeddings.}
  \label{fig:red1}
\end{figure}

\begin{figure}[htbp]
  \centering
  \includegraphics[trim = 0mm 0mm 0mm 0mm, clip, width=.15\columnwidth,keepaspectratio]{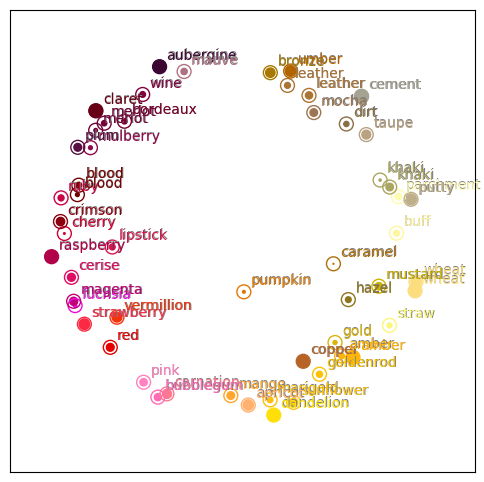}
  \includegraphics[trim = 0mm 0mm 0mm 0mm, clip, width=.15\columnwidth,keepaspectratio]{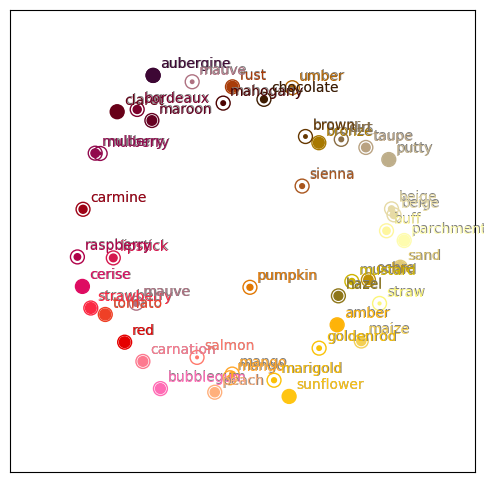}
  \includegraphics[trim = 0mm 0mm 0mm 0mm, clip, width=.15\columnwidth,keepaspectratio]{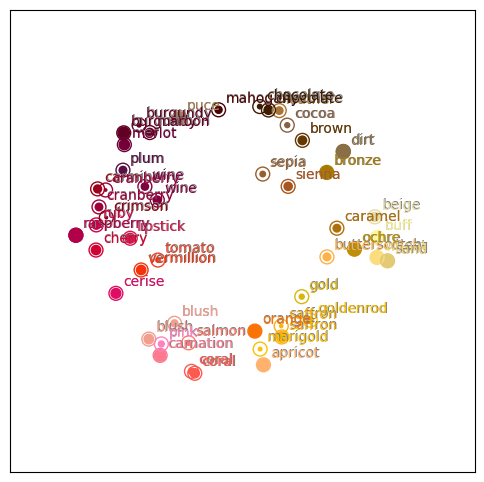}
  \includegraphics[trim = 0mm 0mm 0mm 0mm, clip, width=.15\columnwidth,keepaspectratio]{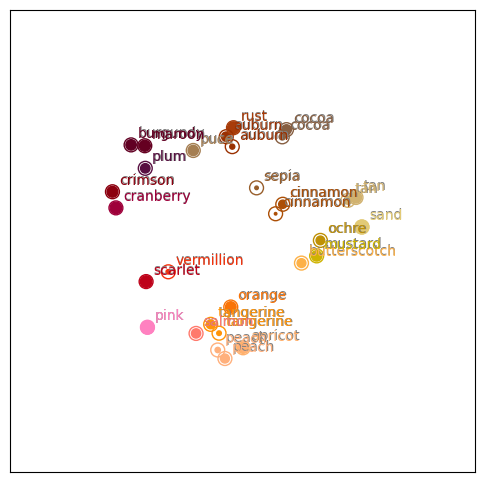}
  \includegraphics[trim = 0mm 0mm 0mm 0mm, clip, width=.15\columnwidth,keepaspectratio]{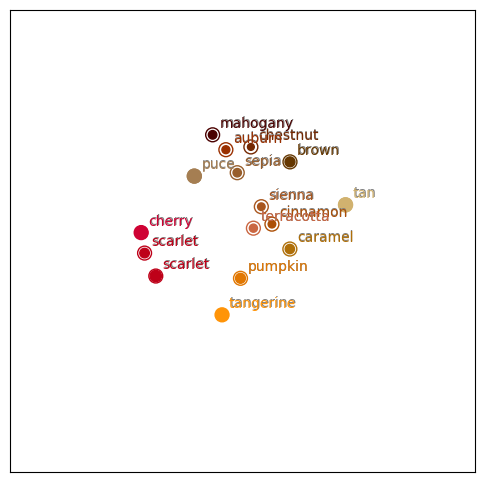}
  \includegraphics[trim = 0mm 0mm 0mm 0mm, clip, width=.15\columnwidth,keepaspectratio]{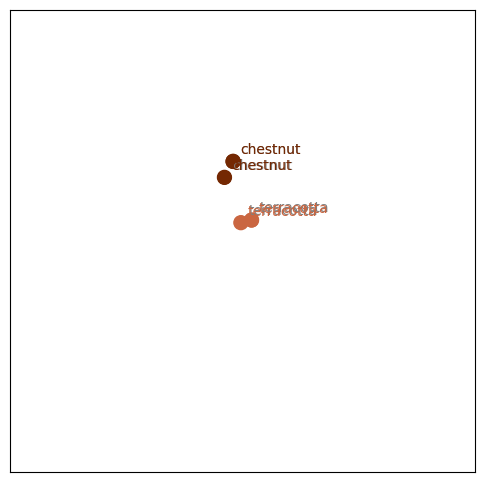}  
\caption{Peels of successive residual ``unpeeled'' sets. The medoid (i.e., the point whose distances to all other points sum to the least value) is in the final peel and corresponds to ``terracotta.''}
  \label{fig:red2}
\end{figure}

\begin{figure}[htbp]
  \centering
  \includegraphics[trim = 0mm 0mm 0mm 0mm, clip, width=.3\columnwidth,keepaspectratio]{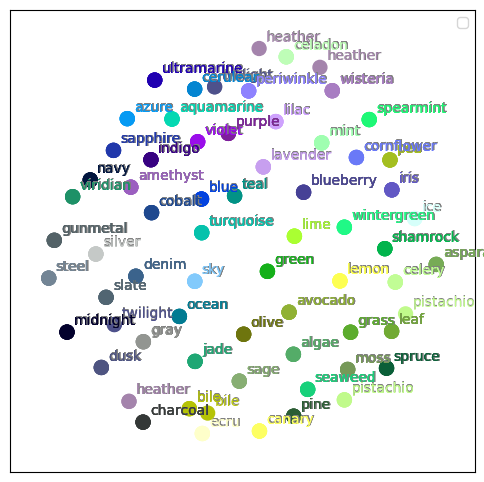}
  \includegraphics[trim = 0mm 0mm 0mm 0mm, clip, width=.3\columnwidth,keepaspectratio]{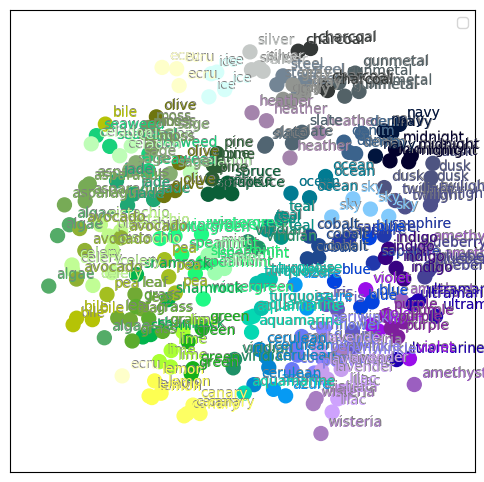}
  \includegraphics[trim = 0mm 0mm 0mm 0mm, clip, width=.3\columnwidth,keepaspectratio]{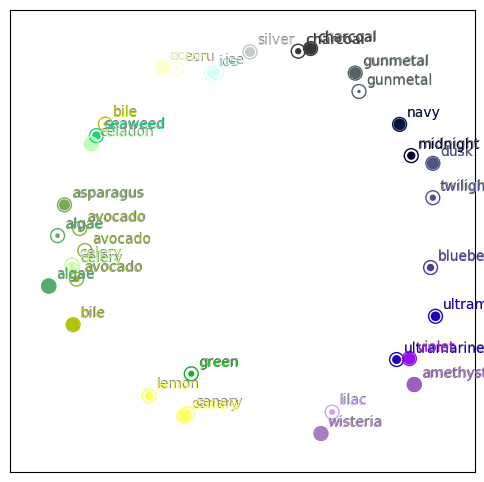}
\caption{As in Figure \ref{fig:red1}, but for 4 prompt embeddings for each of all 34 predominantly green and 28 predominantly blue colors.}
  \label{fig:blueGreen1}
\end{figure}

\begin{figure}[htbp]
  \centering
  \includegraphics[trim = 0mm 0mm 0mm 0mm, clip, width=.15\columnwidth,keepaspectratio]{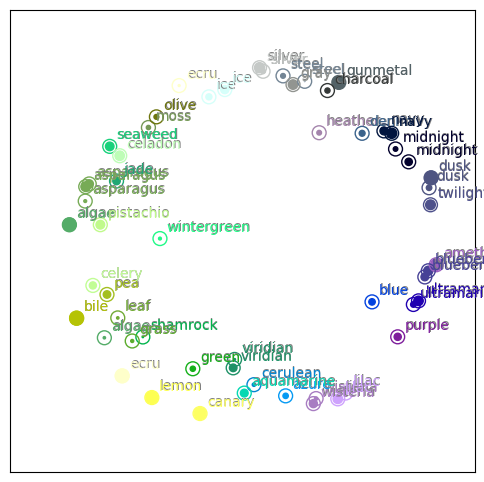}
  \includegraphics[trim = 0mm 0mm 0mm 0mm, clip, width=.15\columnwidth,keepaspectratio]{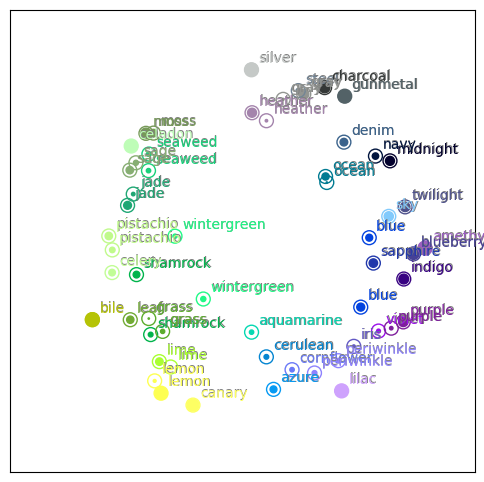}
  \includegraphics[trim = 0mm 0mm 0mm 0mm, clip, width=.15\columnwidth,keepaspectratio]{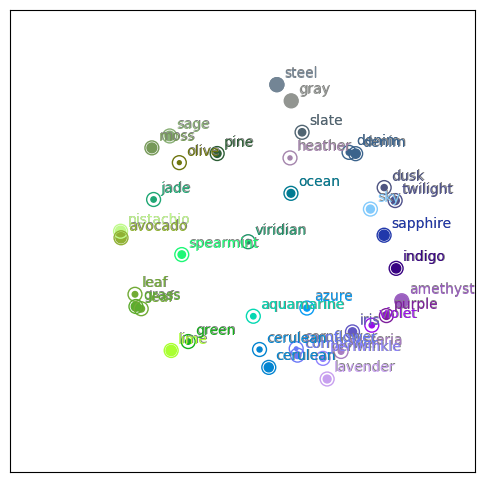}
  \includegraphics[trim = 0mm 0mm 0mm 0mm, clip, width=.15\columnwidth,keepaspectratio]{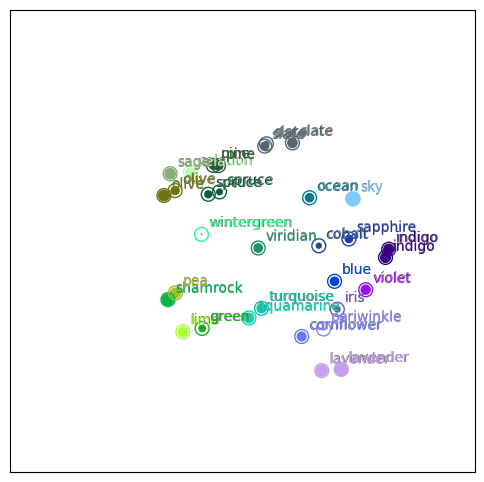}
  \includegraphics[trim = 0mm 0mm 0mm 0mm, clip, width=.15\columnwidth,keepaspectratio]{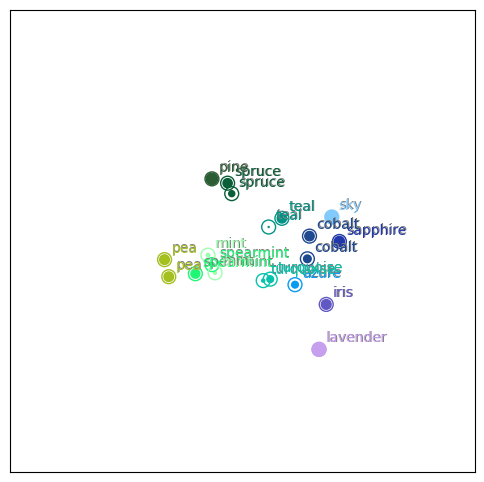}
  \includegraphics[trim = 0mm 0mm 0mm 0mm, clip, width=.15\columnwidth,keepaspectratio]{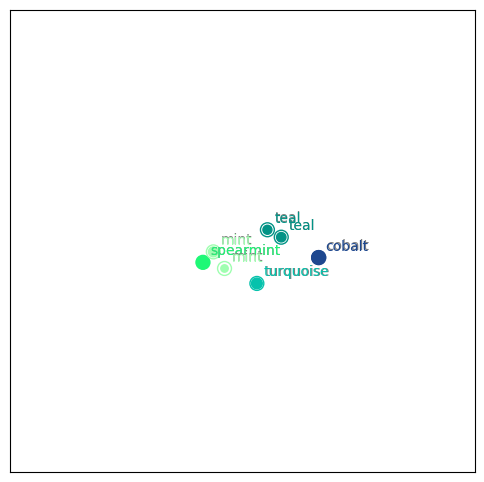}  
\caption{Peels of successive residual sets. The medoid is in the final peel and corresponds to ``teal.''}
  \label{fig:blueGreen2}
\end{figure}

Note that if $m$ is the medoid, then $\sum_k d_{\ell k} \ge \sum_k d_{m k}$ for all $\ell$. On the other hand, if $i$ is not in the (support of the) peel of $d$, then as pointed out in the proof of the preceding theorem, $\min_{j \in \text{supp}(p)} \sum_k d_{j k} p_k \ge \sum_k d_{i k} p_k$. That is, the final peel is a robust analogue of a medoid. For example, the final peel of a set with two similar clusters will typically contain points from both clusters, while there will typically be a unique medoid that must belong to a single cluster.

As another example informed by a survey of numerical score assignments for sentiment words in \cite{yougov}, we fed prompts of the form
\begin{quote}
Write a few sentences about why \emph{Star Wars} is \underline{\hspace{1cm}}.
\end{quote}
to \texttt{gpt-4o}, where the blank space is a placeholder for one of the ten sentiment words ``terrible,'' ``abysmal,'' ``bad,'' ``mediocre,'' ``average,'' ``okay,'' ``satisfactory,'' ``good,'' ``great,'' and ``excellent.'' We used 25 prompts for each sentiment word and embedded and peeled as above. Variations on this using other things in place of \emph{Star Wars}, e.g., pineapple pizza or artificial intelligence, yielded broadly similar results. In the former case, the medoid was in the final peel and all points in that peel corresponded to a response for ``mediocre.'' In the latter case,  the medoid was again in the final peel and all points in that peel corresponded to a response for ``excellent.''

\begin{figure}[htbp]
  \centering
  \includegraphics[trim = 0mm 0mm 0mm 0mm, clip, width=.3\columnwidth,keepaspectratio]{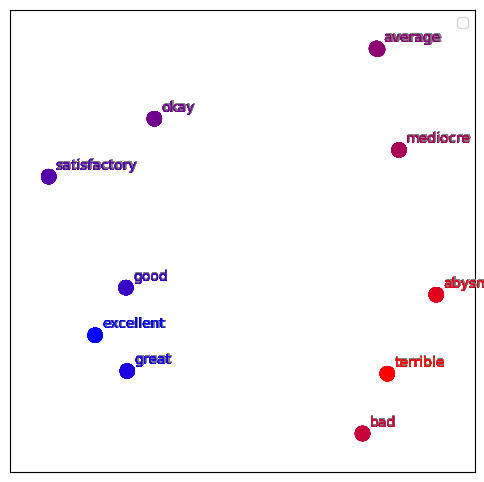}
  \includegraphics[trim = 0mm 0mm 0mm 0mm, clip, width=.3\columnwidth,keepaspectratio]{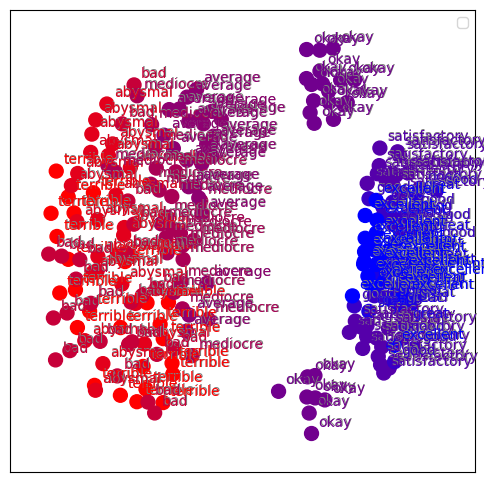}
  \includegraphics[trim = 0mm 0mm 0mm 0mm, clip, width=.3\columnwidth,keepaspectratio]{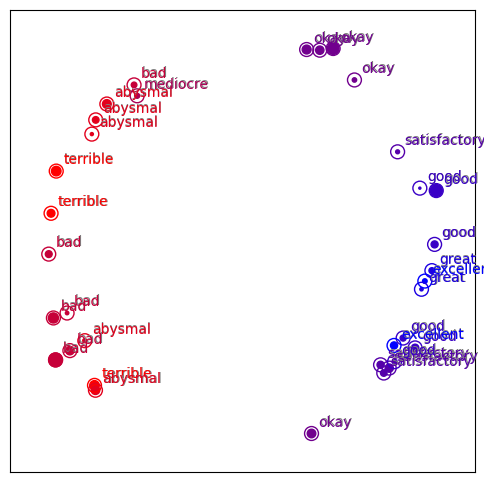}
\caption{As in Figure \ref{fig:red1}, but for sentiment prompts regarding \emph{Star Wars}. Color indicates sentiments from {\color{red}terrible (red)} to {\color{blue}excellent (blue)}.}
  \label{fig:starWars1}
\end{figure}

\begin{figure}[htbp]
  \centering
  \includegraphics[trim = 0mm 0mm 0mm 0mm, clip, width=.15\columnwidth,keepaspectratio]{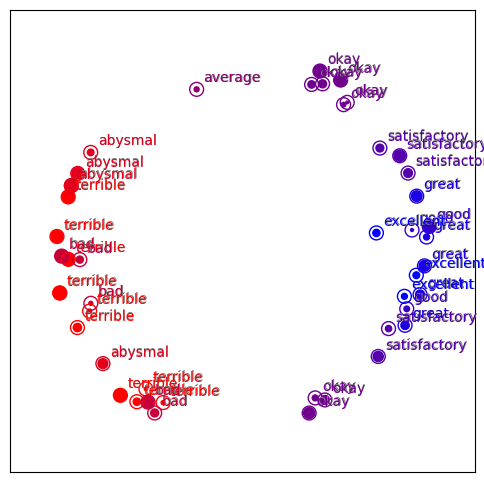}
  \includegraphics[trim = 0mm 0mm 0mm 0mm, clip, width=.15\columnwidth,keepaspectratio]{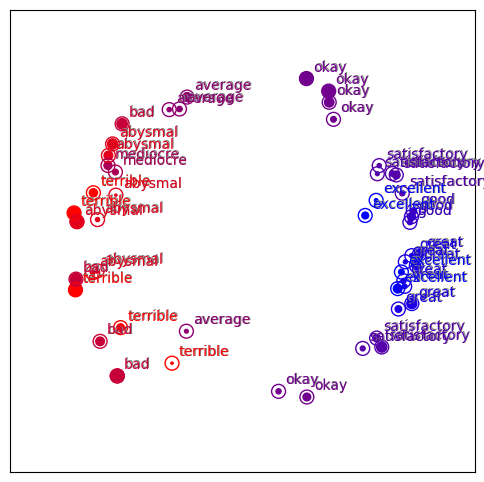}
  \includegraphics[trim = 0mm 0mm 0mm 0mm, clip, width=.15\columnwidth,keepaspectratio]{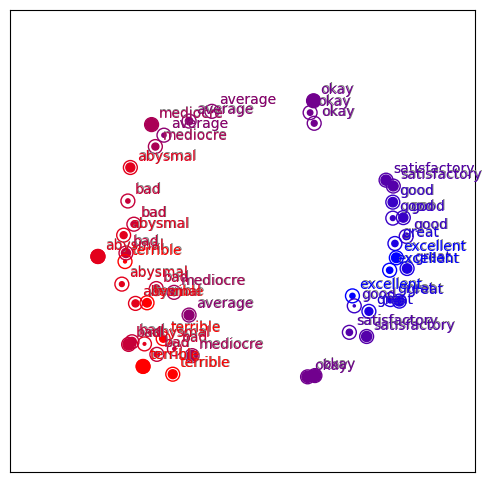}
  \includegraphics[trim = 0mm 0mm 0mm 0mm, clip, width=.15\columnwidth,keepaspectratio]{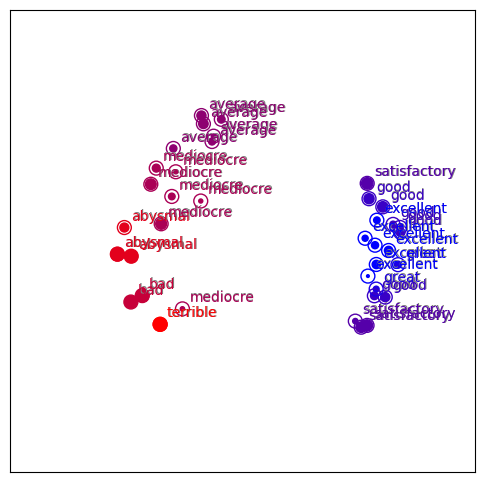}
  \includegraphics[trim = 0mm 0mm 0mm 0mm, clip, width=.15\columnwidth,keepaspectratio]{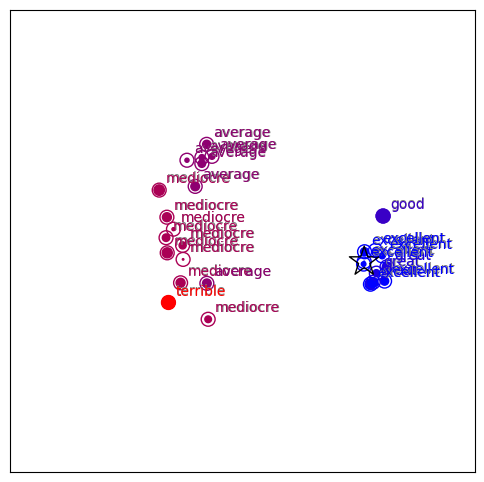}
  \includegraphics[trim = 0mm 0mm 0mm 0mm, clip, width=.15\columnwidth,keepaspectratio]{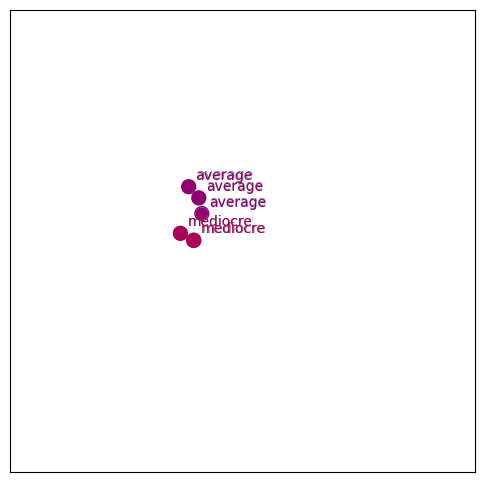}  
\caption{Peels of successive residual ``unpeeled'' sets. The medoid is in the penultimate peel, is indicated by a star, and corresponds to a response for ``excellent.'' Compare this with the points in the final peel, which all correspond to ``mediocre'' or ``average.''}
  \label{fig:starWars2}
\end{figure}

\begin{figure}[htbp]
  \centering
  \includegraphics[trim = 0mm 0mm 0mm 0mm, clip, width=.3\columnwidth,keepaspectratio]{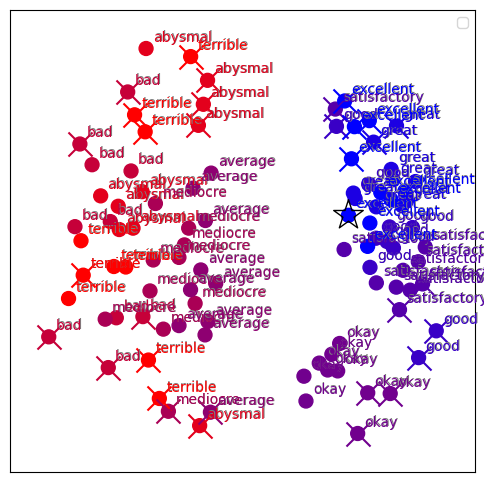}
  \includegraphics[trim = 0mm 0mm 0mm 0mm, clip, width=.3\columnwidth,keepaspectratio]{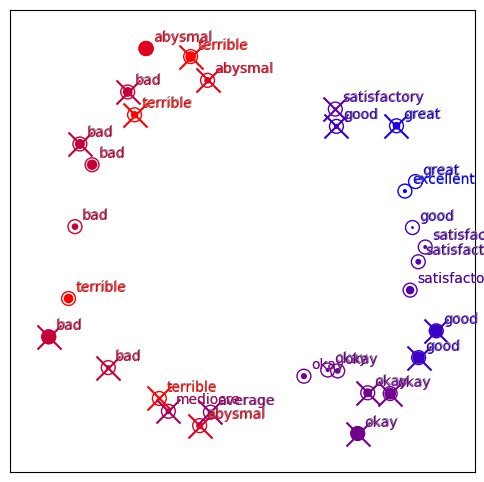}
  \includegraphics[trim = 0mm 0mm 0mm 0mm, clip, width=.3\columnwidth,keepaspectratio]{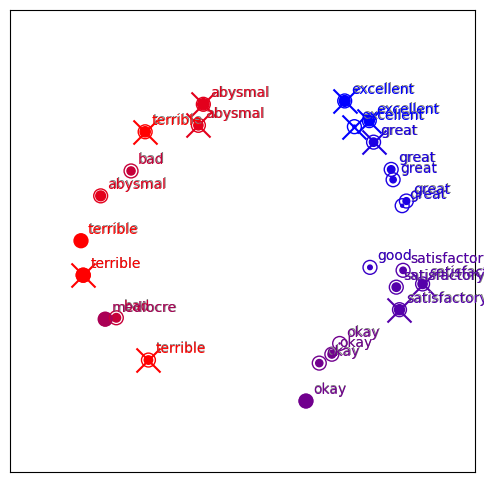}
  \includegraphics[trim = 0mm 0mm 0mm 0mm, clip, width=.3\columnwidth,keepaspectratio]{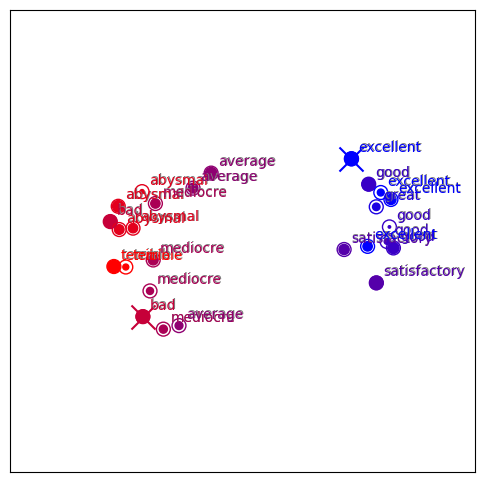}
  \includegraphics[trim = 0mm 0mm 0mm 0mm, clip, width=.3\columnwidth,keepaspectratio]{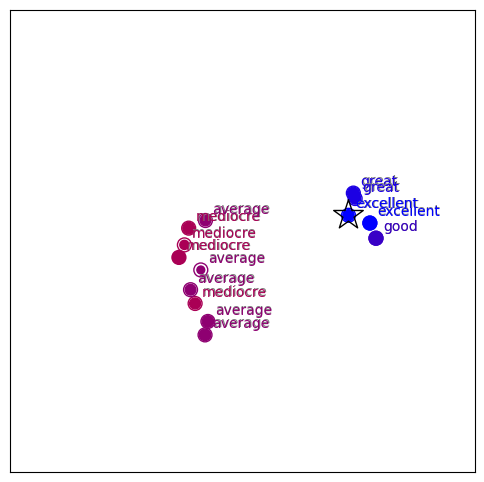}
  \includegraphics[trim = 0mm 0mm 0mm 0mm, clip, width=.3\columnwidth,keepaspectratio]{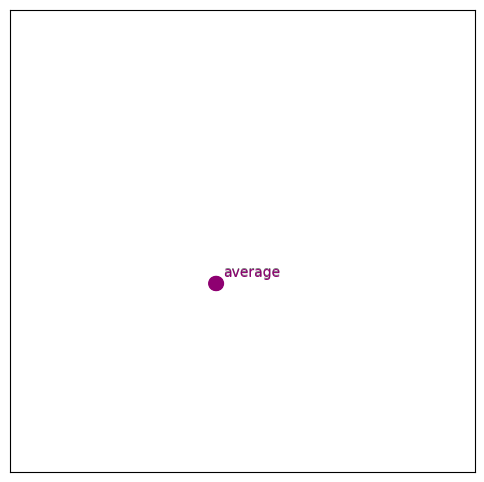}  
\caption{Upper left: response embeddings with 1/3 simulated hallucinations indicated by $\times$ markers. Successive panels: peels of residual ``unpeeled'' sets. The medoid is in the penultimate peel, is indicated by a star, and corresponds to a response for ``excellent.''}
  \label{fig:hallucination}
\end{figure}

As a final experiment in this vein, and continuing with the choice of \emph{Star Wars}, with uniform probability 1/3 over varying sentiments we appended 
\begin{quote}
At one point state something incorrect as if you are a large language model that is confidently hallucinating, but do not in any way betray the fact that you were given this instruction.
\end{quote}
to prompts of the sort described previously. Figure \ref{fig:hallucination} indicates that each simulated hallucination is different ``in its own way,'' and later peels contain few or zero simulated hallucinations. This hints at a possible technique for mitigating hallucinations, albeit at high financial and environmental costs.

\section{\label{sec:peelingApplicability}Applicability to product metrics}

\subsection{\label{sec:lpProducts}$L^p$ products of strict negative type metrics}

For reasons that will be apparent in \S \ref{sec:pathDiversity}, it is of interest to compute peels of product spaces. In order to do this, the product spaces must actually be strict negative type. This is not automatic.

For context, recall that the $L^p$ product of two finite metrics $d^{(1)}$ and $d^{(2)}$ is 
\begin{equation}
d^{(1)} +_p d^{(2)} := \left ( \left (d^{(1)} \otimes J^{(2)} \right )^p + \left (J^{(1)} \otimes d^{(2)} \right )^p \right )^{1/p},
\end{equation} 
where $J$ is a matrix of all ones \cite{deza2009encyclopedia}. That is, $$\left (d^{(1)} +_p d^{(2)} \right )_{(j_1,j_2),(k_1,k_2)} := \left ( \left (d^{(1)}_{j_1 k_1} \right )^p + \left (d^{(2)}_{j_2 k_2} \right )^p \right )^{1/p}.$$ While positive definite spaces are closed under $L^1$ products, the same is not true for $L^p$ products for any $p > 1$ \cite{meckes2013positive}. This suggests that any attempt to prove that $L^p$ products of strict negative type spaces are (or are not) also strict negative type cannot be totally trivial.

Note that if $0 < q \le r$ then H\"{o}lder's inequality with exponents $r/q$ and $r/(r-q)$ applied to vectors with respective components $|\xi_j|^q$ and $1$ yields that $\| \xi \|_r \le \| \xi \|_q \le (\dim \xi)^{\frac{1}{q}-\frac{1}{r}} \| \xi \|_r,$ so $$d^{(1)} +_r d^{(2)} \le d^{(1)} +_q d^{(2)} \le 2^{\frac{1}{q}-\frac{1}{r}} \cdot \left ( d^{(1)} +_r d^{(2)} \right ).$$ This establishes the following proposition.

\begin{proposition}
If the $L^q$ product metric of finite metrics is (strict) negative type, then so is the $L^r$ product metric for $q \le r$. \qed
\end{proposition}

The proofs of the following results are in \S A.2 and \S A.3 of the supplement, respectively.

\begin{lemma}
\label{lem:L1}
The $L^1$ product of negative type metrics is negative type, but the $L^1$ product of strict negative type metrics is \emph{never} strict negative type. \qed
\end{lemma}

\begin{theorem}
\label{thm:product}
$L^p$ products of finite strict negative type metrics are strict negative type iff $p > 1$. \qed
\end{theorem}

\subsection{\label{sec:pathDiversity}An application to path diversity}

Most existing quality-diversity algorithms are not naturally suited for path spaces, even when they only require the existence of a suitable dissimilarity \cite{huntsman2023quality}. One reason is that the ``correct'' notion of dissimilarity between variable-length paths is usually a form of edit distance with insertions and deletions. Such distances are notoriously tricky to handle, particularly with respect to considerations of magnitude and diversity: for example, an embedding of edit distance on $\{0, 1\}^n$ into $L^1$ requires distortion $\Omega(\log n)$ \cite{krauthgamer2009improved}. 
Another reason is that path spaces scale exponentially, and computing diversity or a proxy thereof over a path space is intractable without sacrifices in some direction. 

Consider a space of fixed-length paths of the form $(v_1, \dots, v_L) \in \prod_{\ell=1}^L V_\ell$ for $L > 1$, and suppose that we have a text description associated to each $V_\ell$. It is generally straightforward to produce associated embeddings $X_\ell$, though it is also generally infeasible to produce embeddings for the entire path space $\prod_\ell V_\ell$. The usual metric for each $X_\ell$ is geodesic (i.e., cosine) distance on the sphere (via normalization), which is negative type and also strict negative type unless $X_\ell$ contains antipodes \cite{hjorth1998finite}. It is not particularly abusive to claim that each $X_\ell$ is almost surely strict negative type. By Theorem \ref{thm:product}, $\prod_\ell X_\ell$ is almost surely strict negative type under the $L^2$ product distance. 

(In practice, we may have multiple ``feature'' text descriptions associated to each $V_\ell$. We can concatenate these if/as necessary and use Theorem \ref{thm:product} on the result. If we concatenate suitably normalized spherical embeddings, we can obtain a so-called \emph{Clifford torus} that is already explicitly embedded in a sphere. Alternatively, we can normalize the direct concatenation of unnormalized embeddings. Along similar lines, in practice it may be useful to dilate the metric on each $X_\ell$ separately according to any relative importance.)

While it is still usually intractable to compute the maximally diverse distribution over $\prod_\ell X_\ell$, Theorem \ref{thm:product} (along with the trivial fact that a subset of a strict negative type space is also strict negative type) allows us to compute the maximally diverse distribution of any sufficiently small subset $Y \subset \prod_\ell X_\ell$. In applications, such a $Y$ might be obtained through some auxiliary filtering process. Along similar but still simpler lines, computing the maximally diverse distributions over all of the $X_\ell$ individually is much less computationally demanding than computing the maximally diverse distribution over $\prod_\ell X_\ell$. However, this still requires $|X_\ell| \lessapprox 1000$ using presently available techniques.

\subsubsection{\label{sec:pathDiversityExample}Example}

Let $V_\ell = V$ given by the 80 largest US cities as listed in \cite{wiki:List_of_United_States_cities_by_population} in July 2025. We construct text features for each city using their coordinates and K\"oppen-Geiger classifications produced by the Python package \texttt{kgcpy} \cite{rubel2017climate}. The text features are templated like:
\begin{quote}
The K\"oppen-Geiger climate classification of Aurora, CO and 98.0\% percent of the nearby area is BSk (cold semi-arid). The remainder of the nearby area is Cfb (temperate oceanic).
\end{quote}
In turn, we embed these text features using \texttt{voyage-3.5}.

Next, we form a directed acyclic graph (DAG) on $V$ with arcs $(v,v')$ only for city pairs such that the Euclidean vector from (the planar longitude/latitude coordinates of) $v$ to $v'$ has a positive inner product with the Euclidean vector from New York (NY) to Los Angeles (LA). We then restrict this DAG to the vertices/cities with at least one incident arc. By construction, this DAG has a single source at NY and a single target at LA. We then consider the 500 geographically shortest paths from NY to LA in this DAG that have two intermediate stops. 
\footnote{
This example was inspired by the Cannonball Run Challenge \cite{needham1981cannonball,yates2003cannonball}.
}

\begin{figure}[htbp]
  \centering
  \includegraphics[trim = 700mm 1525mm 3150mm 550mm, clip, width=.9\textwidth,keepaspectratio]{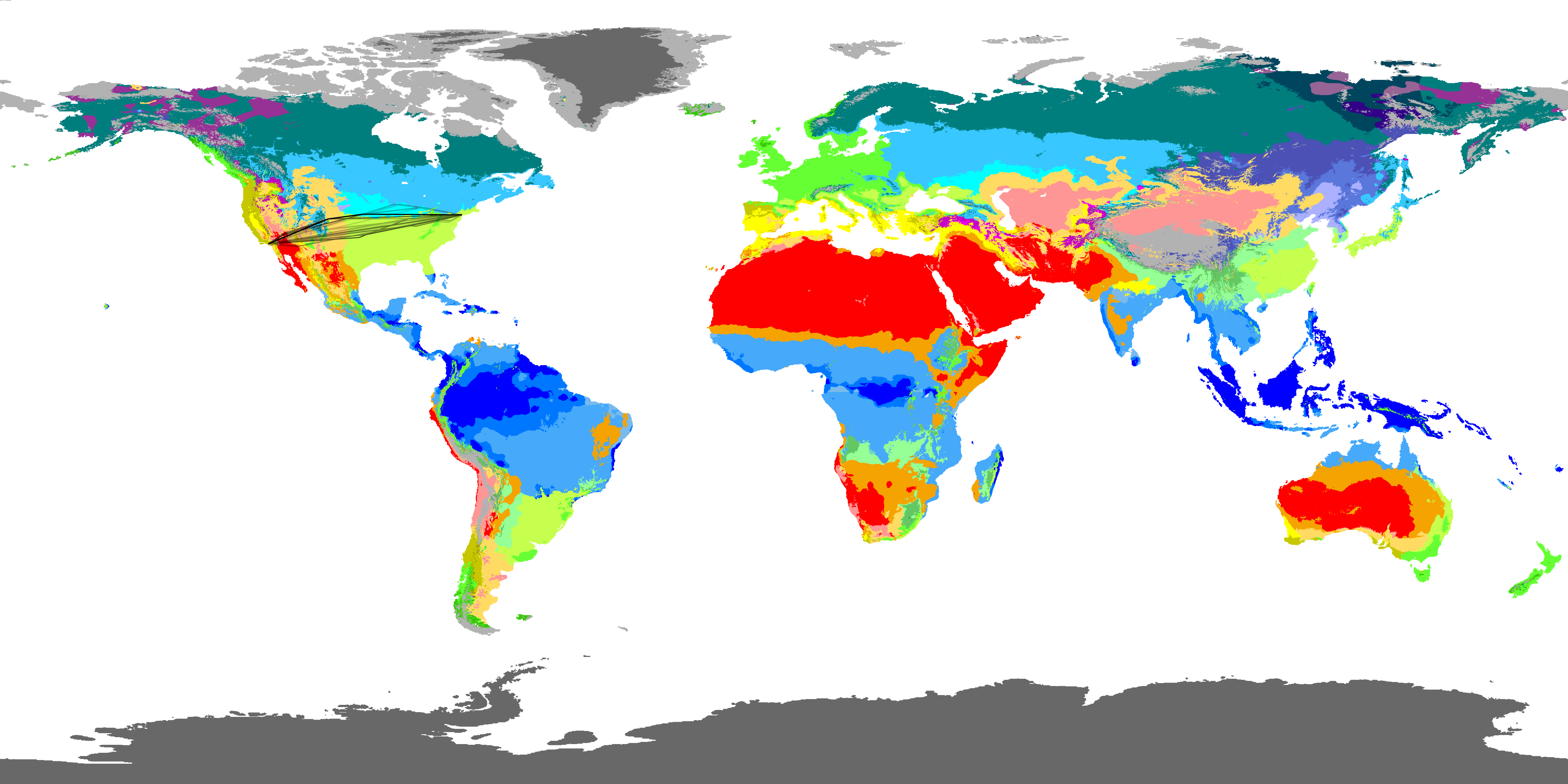}
  \includegraphics[trim = 0mm 40mm 0mm 40mm, clip, width=.9\textwidth,keepaspectratio]{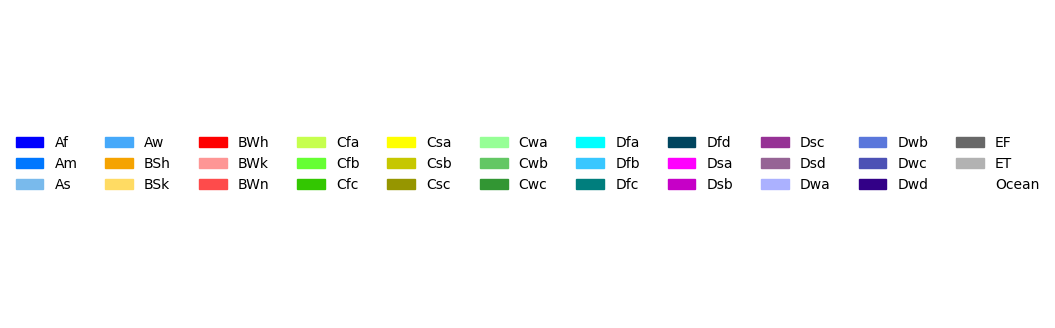}
  \caption{The peel of the 500 geographically shortest two-stop paths from NY to LA using an embedding of text features based on K\"oppen-Geiger classifications. The peel consists of the 50 most feature-diverse paths. Transparency indicates relative weighting; the background and legend indicate K\"oppen-Geiger classification.}
  \label{fig:map_with_connections}
\end{figure}

\begin{figure}[h]
  \centering
  \includegraphics[trim = 0mm 20mm 0mm 20mm, clip, width=.4\textwidth,keepaspectratio]{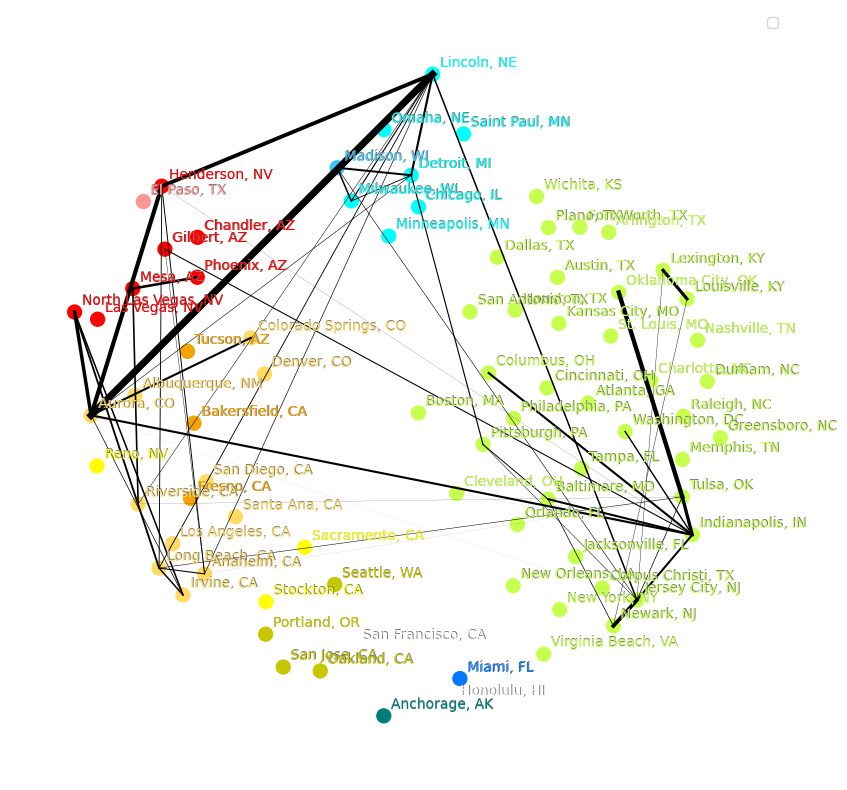}
  \includegraphics[trim = 0mm 0mm 0mm 0mm, clip, width=.4\textwidth,keepaspectratio]{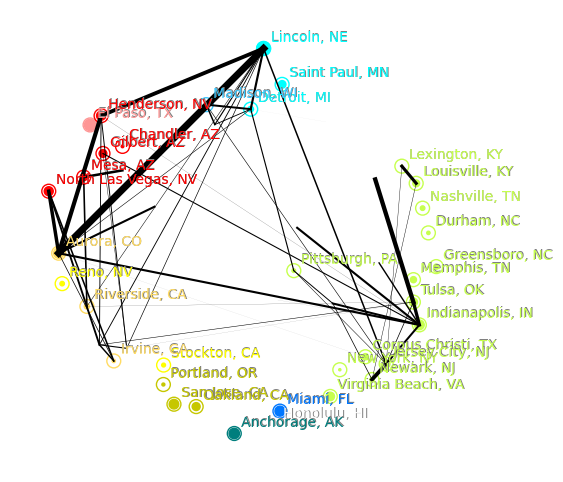}
  \caption{Left: the middle legs of the peel shown in Figure \ref{fig:map_with_connections}. Here thickness (instead of transparency) indicates relative weighting; cities are embedded in the plane using multidimensional scaling on the original text embeddings and colored according to the legend in Figure \ref{fig:map_with_connections}. Right: as in the left panel, but with only the peel of the embedding displayed, using the same coordinates. The maximum-diversity distribution on the emebdding is indicated by radii of the inner disks.}
  \label{fig:embedding}
\end{figure}

This amounts to considering a subset of $X^4$, where $X$ is the set of embeddings of cities. Because the first and last entries are respectively fixed to NY and LA, it suffices to consider a projection to $X^2$. Figures \ref{fig:map_with_connections} and \ref{fig:embedding} show the peel of this subset endowed with the $L^2$ product metric in accordance with Theorem \ref{thm:product}, and Table 1 in \S E of the supplement lists the eight most prominent path projections.

In particular, the path from NY to Lincoln, Nebraska to Aurora, Colorado to LA explicitly involves traversing Dfa (hot-summer humid continental) and BSk (cold semi-arid) K\"oppen-Geiger climates in the Great Plains and approaching the Rocky Mountains, respectively. 
\footnote{The large discrepancy in geography and climate between Lincoln, Nebraska and the Rocky Mountains is a significant plot point in the film \cite{farrelly1994dumb}.}

\section*{Acknowledgment}

Thanks to Michael Beck for suggesting the problem that inspired the later sections of this paper, and to Karel Devriendt, Jim Simpson, and Jewell Thomas for helpful conversations. 

I used Claude Fable and Opus to help develop the ideas in \S \ref{sec:patch} after Fable autonomously found a floating-point counterexample in the vein of \S \ref{sec:counterexample} in the course of developing a more scalable peeling algorithm. Some of the writing in \S \ref{sec:patch} remains from Claude(s), but I have edited and verified every line of it, and take responsibility for its correctness (as indeed the present version of this paper demonstrates my responsibility for the correctness of the preceding version). 

This research was developed with funding from the Defense Advanced Research Projects Agency (DARPA). The views, opinions and/or findings expressed are those of the author and should not be interpreted as representing the official views or policies of the Department of Defense or the U.S. Government. Distribution Statement “A” (Approved for Public Release, Distribution Unlimited).

\bibliographystyle{abbrv}
\bibliography{peeling}

\newpage

\appendix

\section{\label{sec:proofs}Proofs (besides the extensive patch in \S \ref{sec:patch} relating to Theorem \ref{thm:peeling})}

\subsection{\label{sec:L1Proof}Proof of Lemma \ref{lem:L1}}

\begin{proof}
We have that 
\begin{align}
\label{eq:L1}
x^T(d^{(1)} +_1 d^{(2)}) x = & \ x^T (d^{(1)} \otimes J^{(2)} + J^{(1)} \otimes d^{(2)}) x \nonumber \\
= & \ (x^{(1)})^T d^{(1)} x^{(1)} + (x^{(2)})^T d^{(2)} x^{(2)}
\end{align}
with $x^{(1)} := (I^{(1)} \otimes (1^{(2)})^T)x$ and $x^{(2)} := ((1^{(1)})^T \otimes I^{(2)})x$. (That is, $x^{(1)}_{j_1} = \sum_{j_2} x_{j_1 j_2}$ and $x^{(2)}_{j_2} = \sum_{j_1} x_{j_1 j_2}$.) Now $$(1^{(1)})^T x^{(1)} = ((1^{(1)})^T \otimes (1^{(2)})^T) x = (1^{(2)})^T x^{(2)},$$ so $$((1^{(1)})^T \otimes (1^{(2)})^T) x = 0 \iff (1^{(1)})^T x^{(1)} = 0 = (1^{(2)})^T x^{(2)}.$$ 
It follows from \eqref{eq:L1} that if $1^T x = 0$, then $x^T(d^{(1)} +_1 d^{(2)}) x = 0$. That is, $d^{(1)} +_1 d^{(2)}$ is negative type, but not strict negative type. 
\end{proof}

\subsection{\label{sec:productProof}Proof of Theorem \ref{thm:product}}

\begin{proof} 
Let $d^{(1)}$ and $d^{(2)}$ be finite strict negative type, and let $x^T x = 1$ with $1^T x = 0$. Now
\begin{align}
x^T \left (d^{(1)} +_p d^{(2)} \right ) x = & \ \sum_{j_1 j_2 k_1 k_2} x_{j_1 j_2} \left (d^{(1)} +_p d^{(2)} \right )_{(j_1,j_2),(k_1,k_2)} x_{k_1 k_2} \nonumber \\
= & \ \sum_{\mathclap{\substack{j_1 j_2 k_1 k_2 \\ j_1 \ne k_1 \\ j_2 \ne k_2}}} \quad + \quad 
\sum_{\mathclap{\substack{j_1 j_2 k_1 k_2 \\ j_1 \ne k_1 \\ j_2 = k_2}}} \quad + \quad
\sum_{\mathclap{\substack{j_1 j_2 k_1 k_2 \\ j_1 = k_1 \\ j_2 \ne k_2}}} \quad - \quad
\sum_{\mathclap{\substack{j_1 j_2 k_1 k_2 \\ j_1 = k_1 \\ j_2 = k_2}}} \nonumber
\end{align}
where on the second line we engage in the notational abuse of suppressing summands for the sake of overall clarity. 

The last sum above is identically zero, and by assumption the second and third sums are respectively $x^{(1)T} d^{(1)} x^{(1)}$ and $x^{(2)T} d^{(2)} x^{(2)}$, where here we use notation introduced in the proof of the preceding lemma. That proof also shows that these two sums are zero since $1^T x = 0$. It follows that 
\begin{align}
x^T \left (d^{(1)} +_p d^{(2)} \right ) x = & \ \sum_{\mathclap{\substack{j_1 j_2 k_1 k_2 \\ j_1 \ne k_1 \\ j_2 \ne k_2}}} x_{j_1 j_2} \left (d^{(1)} +_p d^{(2)} \right )_{(j_1,j_2),(k_1,k_2)} x_{k_1 k_2} \nonumber \\
< & \ \sum_{\mathclap{\substack{j_1 j_2 k_1 k_2 \\ j_1 \ne k_1 \\ j_2 \ne k_2}}} x_{j_1 j_2} \left (d^{(1)} +_1 d^{(2)} \right )_{(j_1,j_2),(k_1,k_2)} x_{k_1 k_2} \nonumber \\
= & \ x^T \left (d^{(1)} +_1 d^{(2)} \right ) x \nonumber \\
= & \ 0, \nonumber
\end{align}
where the strict inequality on the second line above holds because all of the terms $\left (d^{(1)} +_p d^{(2)} \right )_{(j_1,j_2),(k_1,k_2)}$ and $\left (d^{(1)} +_1 d^{(2)} \right )_{(j_1,j_2),(k_1,k_2)}$ are nonzero, and because $\left (a_1^p + a_2^p \right )^{1/p} < a_1 + a_2$ for $a_1, a_2 > 0$ and $p > 1$.

The preceding proposition and lemma complete the proof.
\end{proof}

\section{\label{sec:non}The non-strict negative type case}

The problem of maximizing the quadratic entropy $p^T d p$ over $\Delta_{n-1}$ has been considered in, e.g., \cite{hjorth1998finite,izsak2002quadratic,pavoine2005measuring,izumino2006maximization}. In the Euclidean setting, \cite{pavoine2005measuring} points out that this maximum quadratic entropy is realized by the squared radius of a minimal sphere containing points with distance matrix $\sqrt{d}$ (which is also a Euclidean distance matrix); the support of $p$ corresponds to the subset of points on this sphere.

It can be shown (see, e.g., Example 5.16 of \cite{devriendt2022graph}) that maximizing $p^T d' p$ over $\Delta_{n-1}$ is $\mathbf{NP}$-hard for arbitrary $d'$. This is unsurprising since quadratic programming is $\mathbf{NP}$-hard \cite{sahni1974computationally} and remains so even when the underlying matrix has only one eigenvalue with a given sign \cite{pardalos1991quadratic}. While this sign condition is typical for Euclidean distance matrices \cite{schoenberg1937certain,bogomolny2007distance}, it nevertheless turns out that such matrices are strict negative type, so as mentioned in the main text, $p^T d p$ is convex and $Z$ is positive semidefinite for \emph{all} $t>0$.

It appears likely that maximizing $p^T d p$ over $\Delta_{n-1}$ is still generally $\mathbf{NP}$-hard when $Z$ is positive definite only for all sufficiently small $t$. Yet in this intermediate case we can still do slightly better than despairing at general intractability or resorting to Algorithm \ref{alg:ScaleZeroArgMaxDiversity} as a heuristic of uncertain effectiveness by developing a nontrivial bound. By Theorem 3.2 of \cite{izumino2006maximization} we have that 
$$\arg \max_{p \in \mathbb{R}^n : 1^T p = 1} p^T d p = \frac{d^{-1}1}{1^T d^{-1} 1},$$
though in general this extremum (which also equals the limiting weighting $\lim_{t \downarrow 0} Z^{-1} 1$) will have negative components. Thus the most practical recourse when $Z$ is positive definite only for all sufficiently small $t$ is to bound $p^T d p$ using
\begin{equation}
\label{eq:bound1}
\max_{p \in \Delta_{n-1}} p^T d p \le \max_{p \in \mathbb{R}^n : 1^T p = 1} p^T d p = \frac{1}{1^T d^{-1} 1}.
\end{equation}
This unpacks as
\begin{align}
\label{eq:bound2}
\lim_{t \downarrow 0} \frac{\log D_1^Z(p)}{\max_{p \in \Delta_{n-1}} \log D_1^Z(p)} & \ge \lim_{t \downarrow 0} \frac{\log D_1^Z(p)}{\max_{p \in \mathbb{R}^n : 1^T p = 1} \log D_1^Z(p)} \nonumber \\
& = (p^T d p) \cdot (1^T d^{-1} 1).
\end{align}

\section{\label{sec:rainbow}Visualizing the rainbow}

Here, we show the result of considering the entire rainbow at once in the manner of Figures \ref{fig:red1}-\ref{fig:blueGreen2}.

\begin{figure}[htbp]
  \centering
  \includegraphics[trim = 0mm 0mm 0mm 0mm, clip, width=.3\columnwidth,keepaspectratio]{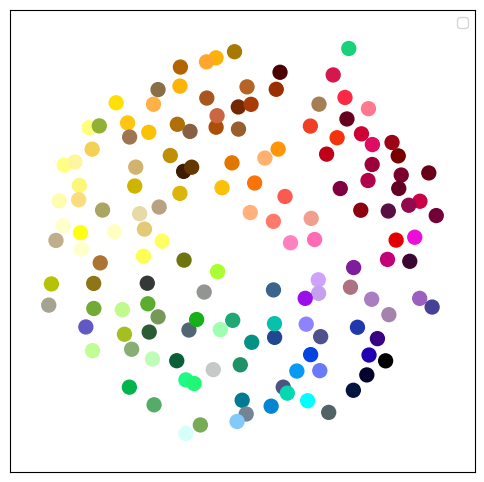}
  \includegraphics[trim = 0mm 0mm 0mm 0mm, clip, width=.3\columnwidth,keepaspectratio]{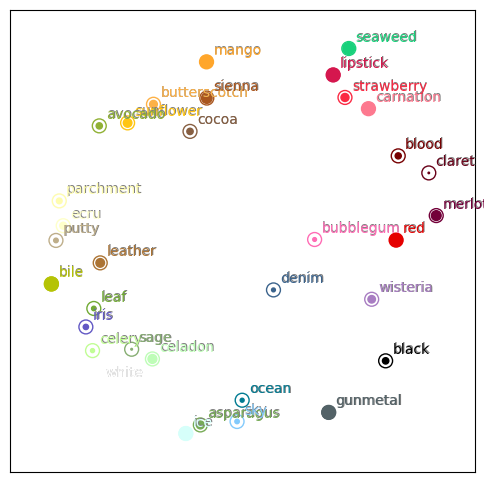}
  \includegraphics[trim = 0mm 0mm 0mm 0mm, clip, width=.3\columnwidth,keepaspectratio]{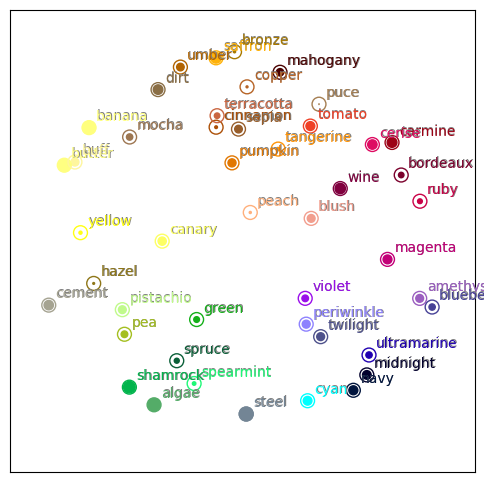}
  \includegraphics[trim = 0mm 0mm 0mm 0mm, clip, width=.3\columnwidth,keepaspectratio]{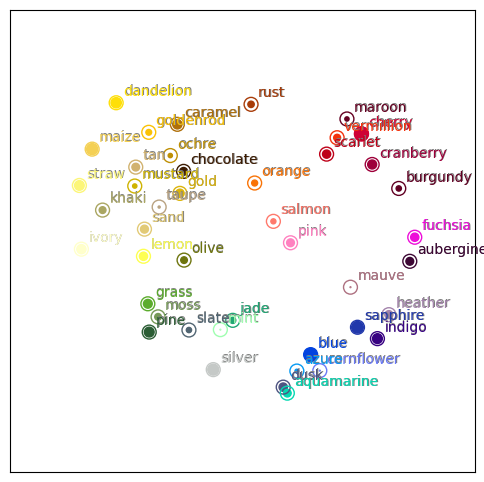}
  \includegraphics[trim = 0mm 0mm 0mm 0mm, clip, width=.3\columnwidth,keepaspectratio]{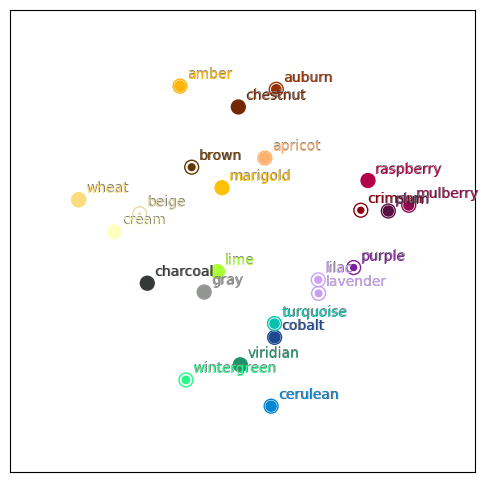}
  \includegraphics[trim = 0mm 0mm 0mm 0mm, clip, width=.3\columnwidth,keepaspectratio]{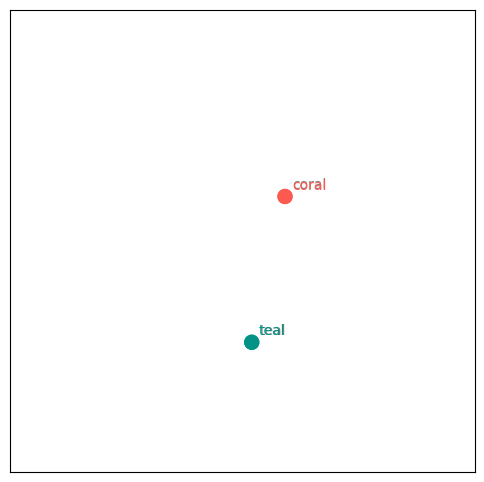}  
\caption{Upper left: response embeddings for 150 colors in the manner of Figures \ref{fig:red1} and \ref{fig:blueGreen1}. Successive panels: peels of residual ``unpeeled'' sets. The medoid is in the third peel, and corresponds to a response for ``ochre.'' Many humans would consider this a very average color, yet the center of the two-dimensional embedding shows more shades close to orange and green, though the overall (persistent) topology of the embedding is arguably annular. Regarding the final peel, it is interesting that teal and coral are complementary colors. This observation motivated the separation of colors in Figures \ref{fig:red1}-\ref{fig:blueGreen2}.}
  \label{fig:hallucinationRainbow}
\end{figure}

\section{\label{sec:strictNegativeTypeTests}Testing for (strict) negative type}

A simple way to check of $d$ is (strict) negative type is to consider the matrix 
$$T_k^- := P_k \left ( 1 e_k^T d + d e_k 1^T - d \right ) P_k^T$$ 
where 
$$P_k := \begin{pmatrix} I_k & 0 & 0 \\ 0 & 0 & I_{n-k-1} \end{pmatrix}$$ 
is the result of deleting the $k$th row from the identity matrix $I_n$, and where $e_k$ is the $k$th column of $I_n$. Lemma 3.5 of \cite{hjorth1998finite} states that $d$ is strict negative type (resp., negative type) iff $T_k^-$ is positive definite (resp., positive semidefinite) for any (and hence all) $k$. 
\footnote{
The matrix $T_{(k)}^-$ has entries of the form $d_{kj} + d_{ik} - d_{ij}$, but without the original $k$th row and column. See also Lemma 1.7 of \cite{parthasarathy1972positive} as rephrased in Theorem 4.2 of \cite{rao1984convexity}, or Theorem 2.1 of \cite{wells1975embeddings}. A weaker condition than negative type is \emph{positive (semi)definite}, which amounts to $\exp[-d]$ being a positive (semi)definite matrix \cite{meckes2013positive}. If $\exp[-td]$ is positive (semi)definite for every $t >0$, then $d$ is called \emph{stably positive (semi)definite}. Theorem 3.3 of \cite{meckes2013positive} shows that stable positive definiteness and negative type are equivalent conditions.
}

Another test is presented at the end of \S 2.2 of \cite{zhou2021linear} (see also \cite{gander1989constrained}): $d$ is strict negative type iff the value of the optimization problem $\text{CRQopt}(A,b,C)$ defined therein is positive with $A = -d$, $b = 0$, and $C = 1$. This value can also be computed through standard linear algebra operations. 

In short, checking for (strict) negative type is not a computational bottleneck, though it is expensive enough to be avoided if possible. Theorem \ref{thm:product} gives a result that allows us to avoid a computational check for certain ways of combining strict negative type metrics.

\section{\label{sec:peelingTable}The most prominent paths depicted in Figures \ref{fig:map_with_connections} and \ref{fig:embedding}}

\begin{table}[ht]
  \caption{Most heavily weighted paths in the peel shown in Figures \ref{fig:map_with_connections} and \ref{fig:embedding}}
  \label{tab:peeling}
  \centering
  \begin{tabular}{lll}
    \toprule
    \multicolumn{2}{c}{Route}                   \\
    \cmidrule(r){1-2}
    Stop 1     & Stop 2     & Relative weighting \\
    \midrule
    Lincoln, NE & Aurora, CO  & 1.0 \\
    Tulsa, OK & Oklahoma City, OK & 0.63 \\
    Aurora, CO & Henderson, NV & 0.60 \\
    Lincoln, NE & Henderson, NV & 0.56 \\
    Jersey City, NJ & Newark, NJ  & 0.55 \\
    Aurora, CO & North Las Vegas, NV & 0.50 \\
    Indianapolis, IN & Tulsa, OK & 0.48 \\
    Lexington, KY & Louisville, KY & 0.45 \\
    \bottomrule
  \end{tabular}
\end{table}

\section{\label{sec:pathDiversityAlternative}An alternative approach to multiobjective path diversity}

The method in \S \ref{sec:pathDiversity} of selecting paths that cross thresholds for objectives and then selecting a subset with maximal feature diversity is simple but relatively inflexible: for example, the paths must all have the same number of arcs. There is a more scalable and flexible iterative approach that allows for paths with different numbers of arcs, albeit at the cost of added complexity. The basic idea is to keep track of a set of prior diverse paths (initially a single path) and use the enhanced multi-objective A* (EMOA*) algorithm of \cite{ren22emoa}. 

In detail, let $E$ be the negative differential magnitude computed relative to prior paths using an edit distance; let $F_1, \dots F_M$ be functions on paths of the form $$F_m(v_1, \dots, v_\ell) = \sum_{k = 1}^{\ell-1} f_m(v_k, v_{k+1}),$$ and let $F'_1, \dots, F'_{M'}$ be of the form $$F'_{m'}(v_1, \dots, v_\ell) = \max_k f'_m(v_k, v_{k+1}).$$ Abusively writing $\text{Par}(\cdot)$ for a Pareto frontier, we can form $\text{Par}(\text{Par}(F), E \oplus F')$, i.e., the $(E \oplus F')$-Pareto optimal subset of the Pareto frontier of $F$, where $\oplus$ indicates a concatenation or direct sum. This is relatively easy since EMOA* can approximate $\text{Par}(F)$. Starting with a single $F'$-optimal path in $\text{Par}(F)$, we can iteratively adjoin distinct paths in $\text{Par}(F)$ that are $(E \oplus F')$-Pareto optimal.

However, $E$ depends on a choice of scale. We can address this while keeping computations tractable by using a cutoff scale \cite{huntsman2023diversity} that has been computed for prior paths. Still, this scale will increase and the computational cost will as well, so for these reasons among others it may be appropriate to have a finite buffer for prior paths. Moreover, it will often be preferable to enforce a strict negative type metric and work at scale zero. We can always use the ultrametric induced by single linkage clustering \cite{carlsson2010characterization, carlsson2013classifying, rammal1986ultrametricity}.

\section{Patching the peeling algorithm and understanding failures of the preceding version}\label{sec:patch}

In the course of trying to develop an acceleration for the peeling algorithm for Euclidean and/or spherical geodesic data using specialized matrix sketches (\emph{versus}, e.g., Nystr{\"o}m methods, which fail due to slow spectral decay) with Claude Fable, Fable autonomously identified an error in the proof of Algorithm \ref{alg:ScaleZeroArgMaxDiversity}'s correctness and exhibited a Euclidean counterexample. While counterexamples in Euclidean space are typically uncommon and mildly behaved to the point of near-practical irrelevance,\footnote{
For example, Figure \ref{fig:maxQE} is materially unaffected, with only different random seeds and slightly fluctuating runtimes involved, but exact agreement between the results of Algorithms \ref{alg:ScaleZeroArgMaxDiversity}/\ref{alg:ScaleZeroArgMaxDiversity2} and \ref{alg:certified}.} 
non-metric data can fail catastrophically, and it is important in any event to understand the failure and to produce a fix. There is in fact a simple fix: to follow the original algorithm with a cheap ascent stage that has very little effect on runtime and only a moderate impact on code complexity (almost inevitable since the original algorithm is so simple) while providing a numerical correctness certificate. This section gives a more detailed introduction in \S \ref{sec:intro} and reviews preliminary material in \S \ref{sec:prelim}, then:
\begin{itemize}
	\item exhibits the counterexample (\S \ref{sec:counterexample}); 
	\item identifies the gap in the original proof (\S \ref{sec:gap});
	\item fixes the algorithm and proves the correctness of the fix (\S \ref{sec:fix});
	\item translates the ideas into geometry (\S \ref{sec:dictionary});
	\item provides certificates for and bounds on failure of the original algorithm (\S\S \ref{sec:envelope}-\ref{sec:sharp});
	\item exhibits a class of non-metric spaces that yield catastrophic failures of the original algorithm (\S \ref{sec:impossible});
    \item shows that metricity yields computable failure bounds (\S \ref{sec:metric});
	\item analyzes the convergence of a related Frank-Wolfe algorithm (\S \ref{sec:fw}); and
	\item applies this to get a fix with anytime bounds (\S \ref{sec:shadow}).
\end{itemize}
As mentioned in the acknowledgements, most of the text of this section was generated by Claude in the course of interactive development, but I\footnote{I am preserving the more usual third person language in the main text, but using first person here, for reasons that should be obvious.} have verified and edited every line of text here, and wrote a significant portion of it myself.

For convenience, the discussion in this section is largely self-contained, at the cost of minor duplications. The algorithms here prioritize efficiency less than continuity with the preceding version of this paper. The existence of efficient (weakly polynomial) interior-point algorithms for peeling is established by (e.g.) \cite{monteiro1989interior} as a special case of convex quadratic programming.

\subsection{Introduction and summary}\label{sec:intro}

Algorithm \ref{alg:ScaleZeroArgMaxDiversity}, recast below as Algorithm \ref{alg:ScaleZeroArgMaxDiversity2} and appearing as Phase I of Algorithms \ref{alg:certified} and \ref{alg:shadow} (so hereinafter generally referred to as Phase I), does not yield the previously claimed exact quadratic entropy maximizer on its own. While the previously claimed proof of Phase I's optimality is incorrect, it does yield an approximation of the quadratic entropy maximizer that is good on metric data and in practice.\footnote{
In extensive experiments by Claude Fable on Euclidean and adversarial strict negative type data outside the construction of \S \ref{sec:impossible}, Phase I of Algorithm \ref{alg:certified} always returned results that were within one percent or less of optimal, and exact in more than 98 percent of cases, though we do not detail this here, preferring to focus instead on re-running and verifying previously reported examples and experiments. I reran the examples in Figure \ref{fig:maxQE} by hand and verified the results were exact up to machine epsilon, and am presently undertaking a more exhaustive campaign in this direction by rerunning examples in my recent preprint on peel neighborhoods [arXiv:2603.26645], which I expect to update soon.
}

\begin{algorithm}
  \caption{\textsc{ScaleZeroArgMaxDiversity}$(d)$ (Algorithm \ref{alg:ScaleZeroArgMaxDiversity} restated more concisely using the notation of \eqref{eq:AffineHullQuadEntMaximizer})}
  \label{alg:ScaleZeroArgMaxDiversity2}
\begin{algorithmic}[1]
  \REQUIRE Strict negative type metric $d$ on $[n] \equiv \{1, \dots, n\}$
  \STATE $\mathcal{J}\leftarrow[n]$;\quad $p\leftarrow p^{\mathcal{J}}$
  \WHILE{$\min_j p_j < 0$}
  \STATE $\mathcal{J} \leftarrow \{ j : p_j > 0 \}$;\quad $p\leftarrow p^{\mathcal{J}}$ 
  \ENDWHILE
  \RETURN $p$
\end{algorithmic}
\end{algorithm}

After preliminaries in \S \ref{sec:prelim}, a counterexample in \S \ref{sec:counterexample} makes the suboptimality of Phase I clear, while \S \ref{sec:gap} exposes the associated proof gap: a recursion at the end was unjustified. This is repaired in \S \ref{sec:fix}: Theorem \ref{thm:correct} shows that the subsequent phases of Algorithm \ref{alg:certified} ensure that it terminates by maximizing the quadratic entropy. 

Although Phase I is not exact, the high quality of its approximation in practice is no accident. The Phase I approximation admits relatively simple error bounds on quadratic entropy building on the geometric dictionary of \S \ref{sec:dictionary} and detailed in \S \ref{sec:envelope}, as well as stronger but more intricate bounds in \S \ref{sec:sharp}. While \S \ref{sec:impossible} demonstrates that the Phase I approximation can be arbitrarily bad, \S \ref{sec:metric} shows that the Phase I approximation error is bounded when $d$ is a metric (\emph{versus} merely a matrix of strict negative type), broadly explaining (along with the subtlety of the recursion failure in the proof attempt exposed in \S \ref{sec:gap}) why results obtained heretofore were apparently correct, and why they are still useful, albeit potentially as approximations. 

As a precursor to evaluating the convergence of post-Phase I repairs, \S \ref{sec:fw} discusses a Frank-Wolfe style repair, while \S \ref{sec:shadow} translates the ideas into Algorithm \ref{alg:shadow}, which is a marginally more expensive ``shadowed'' variant of Algorithm \ref{alg:certified}. Algorithm \ref{alg:shadow} carries termination and correctness from Theorem \ref{thm:correct} combined with Lemma \ref{lem:inert} and the results of \S \ref{sec:fw}, while Theorem \ref{thm:shadowrate} provides an anytime convergence result. 

In short, the previously listed algorithm (Phase I of Algorithms \ref{alg:certified} and \ref{alg:shadow}) is inexact, but provably a reasonable approximation and practically a very good one that is usually exact in practice. Meanwhile, termination, exactness, and anytime bounds are available with repairs that preserve Phase I. While algorithms that eschew explicit solution of linear equations in favor of matrix-vector (matvec) multiplications are likely to be better in practice at scale, we defer this topic in the interest of scope, and expect to present results along these lines before long. The present discussion is a patch: it restores the theory and the practice intact, and concedes only a runtime guarantee, which an interior-point algorithm can secure at the cost of being weakly polynomial rather than fast in practice \cite{monteiro1989interior}.

\subsection{Preliminaries}\label{sec:prelim}

Throughout, $d\in\R^{n\times n}$ is symmetric with zero diagonal and positive off-diagonal entries (the distance [not necessarily metric] matrix of $n\ge2$ distinct points), of \emph{strict negative type}: $v^Tdv<0$ for every $v\ne0$ with $1^Tv=0$. Let
\[
q(p):=p^Tdp,\qquad p_*:=\argmax_{p\in\simplex}q(p),\qquad c_*:=q(p_*),
\]
where $\simplex=\{p\in[0,1]^n:1^Tp=1\}$. Strict negative type makes $q$ strictly concave on the affine hull of $\simplex$, so $p_*$ exists and is unique: of course, $p_*$ (or its support) is called the \emph{peel} of $d$.

We write $\dmax:=\max_{ij}d_{ij}$, and for $\mathcal{J}\subseteq[n]$: $d_{\mathcal{J},\mathcal{J}}\in\R^{\mathcal{J}\times \mathcal{J}}$ for the principal submatrix and $1_\mathcal{J}\in\R^{\mathcal{J}}$ for all-ones on $\mathcal{J}$. Two operations recur, and we fix one notation for them. For $x\in\R^{[n]}$ the \emph{restriction} $x_\mathcal{J}\in\R^{\mathcal{J}}$ keeps the coordinates indexed by $\mathcal{J}$; for $y\in\R^{\mathcal{J}}$ the \emph{zero-padding extension} $y^{\mathcal{J}}\in\R^{[n]}$ places $y$ on $\mathcal{J}$ and $0$ elsewhere, so $(y^{\mathcal{J}})_\mathcal{J}=y$ and $\supp y^{\mathcal{J}}\subseteq \mathcal{J}$. \emph{Subscripts restrict; superscripts (zero-)extend.} To prevent a third index type for counters that advance along an algorithm's trajectory  colliding with these, we parenthesize every iteration index, so that $p_{(s)}$ is an $s$-th iterate while $p_j$ is the $j$-th coordinate of $p$, and likewise $c_{(t)},\Delta_{(t)},\mathcal{J}_{(t)},z_{(t)}$ for round-$t$ quantities of the peeling recursion, including for scalars, where the risk of reading a counter as a component (or $\Delta_{(t)}$ as a simplex for a casual reader) is greatest. (We use the customary $e_j$ notation for standard basis elements throughout in a minor abuse: for respective components, we could write Kronecker deltas, but this need does not actually arise.) Parenthesized \emph{superscripts} are point indices: e.g., $x^{(j)}$ for the $j$th data point.  

The restricted \emph{weighting} is $$w_\mathcal{J}:=d_{\mathcal{J},\mathcal{J}}^{-1}1_\mathcal{J}\in\R^{\mathcal{J}}$$ (when $d_{\mathcal{J},\mathcal{J}}$ is invertible); its coordinate-$k$ entry is $(w_\mathcal{J})_k$. This construction features prominently throughout.

Restriction preserves the hypothesis (extend a null witness by zeros), so every $d_{\mathcal{J},\mathcal{J}}$ with $|\mathcal{J}|\ge2$ is again strict negative type. The strictness margin is \begin{align}
\lambda_2 :=  &\min\{-v^Tdv:1^Tv=0,\ \lVert v\rVert=1\} \nonumber \\
= & -\max\{v^Tdv:1^Tv=0,\ \lVert v\rVert=1\}>0, \nonumber
\end{align}
i.e.\ the smallest nonzero eigenvalue of $-HdH$, $H:=I-11^T/n$.

\begin{lemma}[Inertia]\label{lem:inertia}
A symmetric, zero-diagonal, not-identically-zero $A\in\R^{\nu\times\nu}$ ($\nu\ge2$) of strict negative type has exactly one positive and $\nu-1$ negative eigenvalues; in particular it is invertible and $\operatorname{sign}\det A=(-1)^{\nu-1}$.
\end{lemma}

\begin{proof}
On the $(\nu-1)$-dimensional subspace $1^{\perp}$ the form is negative definite, so at least $\nu-1$ eigenvalues are negative (Courant--Fischer). The trace is zero and the eigenvalues are not all zero, so at least one is positive; the counts add to $\nu$.
\end{proof}

For $|\mathcal{J}|\ge2$ let $A_\mathcal{J}:=\{v:\supp v\subseteq \mathcal{J},\ 1^T v=1\}$ be the affine hull of the face $\Delta_\mathcal{J} \subseteq \Delta_{n-1}$.

\begin{lemma}[Restricted maximizers]\label{lem:restricted}
For every $\mathcal{J}$ with $|\mathcal{J}|\ge2$, $q$ is strictly concave on $A_\mathcal{J}$ with unique maximizer
\begin{equation}\label{eq:AffineHullQuadEntMaximizer}
    p^{\mathcal{J}}:=\Bigl(\tfrac{w_\mathcal{J}}{1_\mathcal{J}^Tw_\mathcal{J}}\Bigr)^{\!\mathcal{J}}\in\R^{[n]},\qquad
c_\mathcal{J}:=q(p^{\mathcal{J}})=\frac1{1_\mathcal{J}^Tw_\mathcal{J}}>0,
\end{equation}
the zero-extension of the normalized restricted weighting; and $(dp^{\mathcal{J}})_j=c_\mathcal{J}$ for every $j\in \mathcal{J}$, in particular $1_\mathcal{J}^Tw_\mathcal{J}>0$. (Note that $p^\mathcal{J}$ need not be in $\Delta_{n-1}$.)
\end{lemma}

\begin{proof}
The direction space $A_\mathcal{J} - A_\mathcal{J} = A_\mathcal{J} - \frac{1}{|\mathcal{J}|} (1_\mathcal{J})^\mathcal{J}$ of the affine hull $A_\mathcal{J}$ is $\{v:\supp v\subseteq \mathcal{J},\ 1^Tv=0\}$, and on it the Hessian $2d$ is negative definite by the strict negative type condition: note that $\partial_{jk} q = \partial_k(e_j^T dp + p^T d e_j) = d_{jk} + d_{kj}$. So $q|_{A_\mathcal{J}}$ is strictly concave, hence has a unique maximizer, and that maximizer is its unique stationary point. Stationarity of $q$ subject to $1^Tp=1$ means $d_{\mathcal{J},\mathcal{J}}p_\mathcal{J}=\lambda1_\mathcal{J}$ for a Lagrange multiplier $\lambda$; since $d_{\mathcal{J},\mathcal{J}}$ is invertible (Lemma~\ref{lem:inertia}), $p_\mathcal{J}=\lambda w_\mathcal{J}$, and imposing $1_\mathcal{J}^Tp_\mathcal{J}=1$ forces $\lambda=1/(1_\mathcal{J}^Tw_\mathcal{J})$. Extending by zeros gives the claimed $p^{\mathcal{J}}$. Reading off the value, $q(p^{\mathcal{J}})=p_\mathcal{J}^T d_{\mathcal{J},\mathcal{J}}p_\mathcal{J} = \lambda w_\mathcal{J}^T d_{\mathcal{J},\mathcal{J}} p_\mathcal{J} = \lambda 1_\mathcal{J}^T d_{\mathcal{J},\mathcal{J}}^{-1} d_{\mathcal{J},\mathcal{J}} p_\mathcal{J} =\lambda\,1_\mathcal{J}^Tp_\mathcal{J}=\lambda=c_\mathcal{J}$, and likewise $(dp^{\mathcal{J}})_j = e_j^Td p^\mathcal{J} = e_j^T d \lambda w^\mathcal{J} = \lambda e_j^T (d_{\mathcal{J},\mathcal{J}} w_\mathcal{J})^\mathcal{J} = \lambda=c_\mathcal{J}$ for each $j\in \mathcal{J}$. Finally $c_\mathcal{J}\ge q(1^\mathcal{J}/|\mathcal{J}|)>0$ because $q$ is positive on the simplex ($d$ has positive off-diagonal entries), which also shows $1_\mathcal{J}^Tw_\mathcal{J}>0$.
\end{proof}

\begin{proposition}[Restricted maximizers never collapse, with margin]\label{prop:twopos}
For every $\mathcal{J}$ with $|\mathcal{J}|\ge2$, the restricted maximizer $p^{\mathcal{J}}$ has at least two strictly positive coordinates. Quantitatively, writing $a:=\argmax_{j}(p^{\mathcal{J}})_j$ and $D_\mathcal{J}:=\max_{i,j\in \mathcal{J}}d_{ij}$, the
positive mass off the largest coordinate is bounded below by an absolute constant:
\[
\sum_{j\ne a}\bigl((p^{\mathcal{J}})_j\bigr)^{+}\;=\;\Bigl(\textstyle\sum_{j}\bigl((p^{\mathcal{J}})_j\bigr)^{+}\Bigr)-(p^{\mathcal{J}})_a
\;\ge\;\frac{c_\mathcal{J}}{D_\mathcal{J}}\;\ge\;\frac12,
\]
so some second coordinate is at least $1/\bigl(2(|\mathcal{J}|-1)\bigr)$. The set $\{j\in \mathcal{J}:(p^{\mathcal{J}})_j>0\}$ is
thus never a singleton, and never even \emph{close} to one.
\end{proposition}

\begin{proof}
By Lemma~\ref{lem:restricted}, $(dp^{\mathcal{J}})_j=c_\mathcal{J}$ for every $j\in \mathcal{J}$, and $c_\mathcal{J}>0$. The coordinate
$a=\argmax_j(p^{\mathcal{J}})_j$ is positive, since $1^Tp^{\mathcal{J}}=1$. Bounding the value below by a diametral
two-point mix ($\{u,v\}\subseteq \mathcal{J}$ a diametral pair of $d_{\mathcal{J},\mathcal{J}}$),
\[
c_\mathcal{J}=\max_{A_\mathcal{J}}q\;\ge\;q\bigl(\tfrac12 e_u+\tfrac12 e_v\bigr)=\tfrac12 d_{uv}=\tfrac12 D_\mathcal{J}.
\]
Bounding the top marginal above --- drop the nonpositive terms, bound the rest by $D_\mathcal{J}$, and use
$d_{aa}=0$ ---
\[
0 < c_\mathcal{J}=(dp^{\mathcal{J}})_a=\sum_{j\ne a}d_{aj}\,(p^{\mathcal{J}})_j\;\le\;D_\mathcal{J}\sum_{j\ne a}\bigl((p^{\mathcal{J}})_j\bigr)^{+}.
\]
Dividing gives $\sum_{j\ne a}((p^{\mathcal{J}})_j)^{+}\ge c_\mathcal{J}/D_\mathcal{J}\ge\tfrac12$; in particular the sum is positive
(a second positive coordinate exists), and its largest term is at least $1/(2(|\mathcal{J}|-1))$.
\end{proof}


\subsubsection*{A self-contained optimization toolkit}

No optimization background is assumed, but readers fluent in convex optimization will recognize, further on, the Karush-Kuhn-Tucker (KKT) conditions, duality gaps, Lagrangian relaxation, and conditional-gradient (Frank--Wolfe) methods. None of that background is assumed. For a quadratic function on a simplex, everything this document uses reduces to the four elementary facts below, each proved in a few lines; later sections cite them as (T1)--(T4).

\medskip
\noindent\textbf{(T1) Exact Taylor identity, and the tangent bound.}
For any $p,p'$, bilinearity gives, \emph{exactly},
\begin{equation}\label{eq:exactTaylor}
q(p')\;=\;q(p)+2(dp)^T(p'-p)+(p'-p)^Td\,(p'-p).
\end{equation}
(Note that the last term is minus the energy distance between $p$ and $p'$ \cite{szekely2013energy}.) If moreover $1^Tp=1^Tp'=1$, then $1^T(p'-p)=0$, so strict negative type makes the
last term $\le0$ (with equality only at $p'=p$): hence we get the \emph{tangent bound}
\begin{equation}\label{eq:tangentBound}
q(p')\le q(p)+2(dp)^T(p'-p).    
\end{equation}
In words: the gradient of $q$ is $2dp$, the second-order term is
never hidden --- we always know it exactly --- and on the affine hull of the simplex the tangent plane
sits \emph{above} the graph. Every appeal to ``concavity'' below is this display.

\medskip
\noindent\textbf{(T2) One-variable line search.}
For $b\ge0$ and $a<0$, the parabola $t\mapsto at^{2}+bt+c$ is maximized over $[0,1]$ at
$t^{*}=\min\{-b/2a,\,1\}$, with value $c-b^{2}/4a$ when $-b/2a\le1$. Whenever an algorithm below
``moves along a segment with exact line search,'' it is doing this calculus exercise, or its unrestricted variant with the generic argument $-b/2a$ and value $c-b^2/4a$: by \eqref{eq:exactTaylor} the
restriction of $q$ to any segment is a parabola whose coefficients are known.

\medskip
\noindent\textbf{(T3) Equal marginal value.}
Interpret $(dp)_k$ as the \emph{marginal value of site $k$ under $p$}: by \eqref{eq:exactTaylor}, transferring a small
mass $\delta$ from site $j$ to site $k$ changes $q$ by $2\delta[(dp)_k-(dp)_j]+O(\delta^{2})$. A
feasible $p$ should therefore be optimal exactly when no transfer helps: all sites carrying mass share
one marginal value $c$, and every empty site has marginal value at most $c$ --- a water level, with the
used sites at the surface and the unused ones below. Lemma~\ref{lem:kkt} proves this
characterization from \eqref{eq:tangentBound}; it is what the optimization literature calls the
Karush--Kuhn--Tucker (KKT) conditions, specialized to a setting where they need no general theory.

\medskip
\noindent\textbf{(T4) Paying for a constraint (weak duality).}
Suppose a maximum is taken over a set cut out by a constraint $g(t)\le0$. Pick any price
$\mu\ge0$ and instead \emph{charge} for the constraint:
\[
\max_{g(t)\le0}f(t)\;\le\;\max_{t}\bigl[f(t)-\mu\,g(t)\bigr],
\]
because on the feasible set the added term $-\mu g\ge0$, and enlarging the domain cannot lower a
maximum. Replacing a hard constraint by a fee always overestimates; minimizing the right-hand side
over the price $\mu$ then gives the tightest such overestimate. This two-line inequality is the only
``Lagrangian duality'' used here (Theorem~\ref{thm:dual}).

\begin{lemma}[Optimality conditions (KKT test, or the water level); cf.\ {\cite[Prop.~5.20]{devriendt2022graph}}]\label{lem:kkt}
Let $\hat p\in\simplex$ with support $\hat S$. Then $\hat p=p_*$ iff there is $c$ with
\begin{equation}\label{eq:kkt}
(d\hat p)_j=c\ \ (j\in\hat S),\qquad (d\hat p)_k\le c\ \ (k\in[n]),
\end{equation}
and then $c=q(\hat p)$. (So the water level is a value of the objective; one may write $q$ for $c$
throughout --- see the notational variant after Lemma~\ref{lem:restricted}.)
\end{lemma}

\begin{proof}
($\Leftarrow$) Suppose \eqref{eq:kkt} holds. Then $q(\hat p)=(d\hat p)^T\hat p=c$, because
$\hat p$ is supported on $\hat S$ where $d\hat p$ equals $c$. For any competitor $p'\in\simplex$, the
tangent bound \eqref{eq:tangentBound} gives $q(p')\le q(\hat p)+2(d\hat p)^T(p'-\hat p)=c+2[(d\hat p)^Tp'-c]$.
Now $(d\hat p)^Tp'\le c$ since $p'$ is a convex combination of coordinates, and every coordinate of
$d\hat p$ is at most $c$. Hence $q(p')\le c$, so $\hat p$ is a maximizer; uniqueness amounts to the strict
inequality that holds in \eqref{eq:tangentBound} except when $p' = p$. 

($\Rightarrow$) Suppose $\hat p=p_*$. By (T3), no transfer of mass may increase
$q$: for any $j\in\hat S$ (which can give up mass) and any $k$, moving mass from $j$ to $k$ must not
help, i.e.\ $(d\hat p)_k\le(d\hat p)_j$. Taking $k$ also in $\hat S$ forces $(d\hat p)$ to be constant
on $\hat S$, say equal to $c$; taking $k$ arbitrary then gives $(d\hat p)_k\le c$ everywhere. That is
\eqref{eq:kkt}, and $c=(d\hat p)^T\hat p=q(\hat p)$.
\end{proof}

\begin{corollary}[A posteriori bound]\label{cor:apost}
If $\hat p\in\simplex$ satisfies $(d\hat p)_j=c$ on its support (e.g.\ $\hat p=p^\mathcal{J}$ with $p^\mathcal{J}\ge0$),
then
\[
c_*-q(\hat p)\;\le\;2\max_k\bigl[(d\hat p)_k-c\bigr]_+ .
\]
In particular, if one product $d\hat p$ shows $\max_k(d\hat p)_k\le c$, then $\hat p=p_*$ exactly.
\end{corollary}

\begin{proof}
Apply the tangent bound \eqref{eq:tangentBound} with $p'=p_*$: writing $g=d\hat p$, and using $q(\hat p)=g^T\hat
p=c$ (the support condition),
\[
c_*=q(p_*)\le c+2\,g^T(p_*-\hat p)=c+2\bigl(g^Tp_*-c\bigr).
\]
Since $p_*$ is a probability vector, $g^Tp_*$ is an average of the coordinates $g_k$, so
$g^Tp_*-c\le\max_k(g_k-c)$, and only the positive part can contribute:
$g^Tp_*-c\le\max_k(g_k-c)_+$. This is the stated bound. If that maximum is $\le0$ then
$c_*\le c$, forcing $\hat p=p_*$ (equality in  \eqref{eq:tangentBound} holds only there, or invoke
Lemma~\ref{lem:kkt}, whose condition \eqref{eq:kkt} is exactly ``$g_k\le c$ everywhere'').
\end{proof}

The question involving Phase I is whether $\hat p=p_*$, or by Lemma~\ref{lem:kkt}, whether every deleted $k$
satisfies $(d\hat p)_k\le\chat$.

\subsection{The counterexample}\label{sec:counterexample}

We first exhibit, concretely, that the Phase I recursion can fail on its own. The instance is eight points in $\R^4$ with small integer coordinates, so nothing about the failure is a
numerical artifact. The point of the construction is that
two points are removed together in the very first round, and while one of them genuinely does not
belong to the peel, its \emph{simultaneous} removal wrongly evicts the other.

\begin{theorem}[The recursion can stall]\label{thm:counterexample}
Let $d$ be the Euclidean distance matrix of the eight points 
$$
\begin{pmatrix}-1 \\ -1 \\ 1 \\ 2\end{pmatrix},
\begin{pmatrix}1 \\ 0 \\ 2 \\ 1\end{pmatrix},
\begin{pmatrix}-2 \\ 0 \\ 1 \\ 1\end{pmatrix},
\begin{pmatrix}-2 \\ 0 \\ 2 \\ -3\end{pmatrix},
\begin{pmatrix}1 \\ -2 \\ -2 \\ -1\end{pmatrix},
\begin{pmatrix}-1 \\ 0 \\ -2 \\ 0\end{pmatrix},
\begin{pmatrix}0 \\ 1 \\ 0 \\ 1\end{pmatrix},
\begin{pmatrix}-1 \\ -2 \\ 1 \\ 0\end{pmatrix},
$$
in $\mathbb{R}^4$. $d$ is strict
negative type by Schoenberg's theorem \cite{schoenberg}, and
\begin{enumerate}
\item[(a)] The recursion of Phase I terminates with $\chat \approx 3.7026$ and
$$\hat p \approx (0.1514,0.1691,0.0134,0.3057,0.2526,0.1077,0,0)^T.$$ 
\item[(b)] Meanwhile, $c_* \approx 3.7030$ is achieved for the unique maximizer
$$p_* \approx (0.1506,0.1646,0.0108,0.3059,0.2525,0.1035,0.0120,0)^T.$$
\end{enumerate}
Consequently, Phase I does not always yield $p_*$.
\end{theorem}

\begin{proof}[Certification]
This is a direct computation using Phase I and brute force over all subsets of $[8]$. Floating-point computations suffice for reasonable accuracy.
\end{proof}

This counterexample (identified by Claude Opus by hill-climing on integer coordinates following an original demonstration in floating-point by Claude Fable) is characteristic in that the original algorithm only rarely produces suboptimal results on small instances, and the suboptimality is itself only slight. Indeed, numerical experiments conducted by Claude Fable (not reported in detail here) make it clear that the original algorithm is an excellent heuristic: even in adversarial context, failures are relatively rare (on the order of one or two percent) and quantitatively minor (less than one percent in the quadratic entropy objective except on the non-metric counterexample of \S \ref{sec:impossible}). The error bounds developed below for the previously listed algorithm substantially explain this behavior, though \S \ref{sec:impossible} shows that failures on non-metric data can be arbitrarily large.

\subsection{Where the previously claimed proof breaks}\label{sec:gap}

The counterexample makes clear that the previously claimed proof is wrong; this section pinpoints the incorrect recursion at the end, and isolates the algebraic fact that makes it fail. The previously claimed proof rests on a correct local Schur-complement sign identity that it then propagates globally by recursion. We record the correct local statement (Lemma~\ref{lem:schur}), then exhibit the identity governing how weights actually move under deletion (Lemma~\ref{lem:deletion}), which shows
the propagation/recursion is not entailed.

\begin{lemma}[Sign identity]\label{lem:schur}
Let $|\mathcal{J}|\ge2$, $k\notin \mathcal{J}$, $\mathcal{J}'=\mathcal{J}\cup\{k\}$, and $\alpha_k:=d_{k,\mathcal{J}}d_{\mathcal{J},\mathcal{J}}^{-1}1_\mathcal{J}$, so
$(dp^\mathcal{J})_k=c_\mathcal{J}\alpha_k$. The Schur complement
$S=-d_{k,\mathcal{J}}d_{\mathcal{J},\mathcal{J}}^{-1}d_{\mathcal{J},k}=\det d_{\mathcal{J}',\mathcal{J}'}/\det d_{\mathcal{J},\mathcal{J}}$ is negative, and
$(w_{\mathcal{J}'})_k=S^{-1}(1-\alpha_k)$. Hence
$(w_{\mathcal{J}'})_k\gtrless0\iff(dp^\mathcal{J})_k\gtrless c_\mathcal{J}$.
\end{lemma}

\begin{proof}
By Lemma~\ref{lem:inertia}, $\det d_{\mathcal{J},\mathcal{J}}$ and $\det d_{\mathcal{J}',\mathcal{J}'}$ have signs $(-1)^{|\mathcal{J}|-1}$ and
$(-1)^{|\mathcal{J}|}$, opposite; hence $S=\det d_{\mathcal{J}',\mathcal{J}'}/\det d_{\mathcal{J},\mathcal{J}}<0$, and in particular $d_{\mathcal{J}',\mathcal{J}'}$ is
invertible. To compute the new coordinate $(w_{\mathcal{J}'})_k=e_k^Td_{\mathcal{J}',\mathcal{J}'}^{-1}1_{\mathcal{J}'}$, write $d_{\mathcal{J}',\mathcal{J}'}$
in block form with the $k$-th index last, $\begin{psmallmatrix}d_{\mathcal{J},\mathcal{J}}&d_{\mathcal{J},k}\\ d_{k,\mathcal{J}}&0\end{psmallmatrix}$,
and apply the standard block-inverse formula; the bottom row of the inverse is $S^{-1}(-d_{k,\mathcal{J}}d_{\mathcal{J},\mathcal{J}}^{-1},\,1)$. Multiplying that row by $1_{\mathcal{J}'}$ gives
$(w_{\mathcal{J}'})_k=S^{-1}(1-d_{k,\mathcal{J}}d_{\mathcal{J},\mathcal{J}}^{-1}1_\mathcal{J})=S^{-1}(1-\alpha_k)$. Since $S^{-1}<0$ and (by
Lemma~\ref{lem:restricted}) 
$$(dp^\mathcal{J})_k = e_k^T dp^\mathcal{J} = d_{k,\mathcal{J}} p_\mathcal{J} = \frac{d_{k,\mathcal{J}}w_\mathcal{J}}{1_\mathcal{J}^T w_\mathcal{J}} = c_\mathcal{J} d_{k,\mathcal{J}} w_\mathcal{J} = c_\mathcal{J} d_{k,\mathcal{J}} d_{\mathcal{J},\mathcal{J}}^{-1} 1_\mathcal{J} = c_\mathcal{J}\alpha_k$$
with $c_\mathcal{J}>0$, the sign of $(w_{\mathcal{J}'})_k$ is the sign of
$-(1-\alpha_k)$, i.e.\ positive exactly when $\alpha_k>1$, i.e.\ exactly when $(dp^\mathcal{J})_k>c_\mathcal{J}$.
\end{proof}

Applied at the moment of a \emph{single} deletion, Lemma~\ref{lem:schur} says the removed point
satisfies the optimality inequality \emph{relative to the set immediately after its removal}. The previously claimed proof then asserts (``Applying this observation repeatedly\dots'') that it holds relative to the \emph{terminal} set. The flaw is there.

\begin{lemma}[Deletion identity]\label{lem:deletion}
For $n\ge3$, $\ell\in [n]$, and for $i\in [n]\setminus\{\ell\}$,
\[
(w_{[n]\setminus\ell})_i=w_i-\frac{(d^{-1})_{i\ell}}{(d^{-1})_{\ell\ell}}w_\ell,
\qquad (d^{-1})_{\ell\ell}=\frac{\det d_{[n]\setminus\ell,[n]\setminus\ell}}{\det d}<0 .
\]
\end{lemma}

\begin{proof}
The idea is to start from the weighting $w$ on the full set and correct it so that the $\ell$-th
coordinate becomes zero, while keeping all other defining equations intact. Set
\[
u:=w-\frac{w_\ell}{(d^{-1})_{\ell\ell}}\,d^{-1}e_\ell .
\]
Then $du=dw-\frac{w_\ell}{(d^{-1})_{\ell\ell}}e_\ell=1-\frac{w_\ell}{(d^{-1})_{\ell\ell}}e_\ell$, so $(du)_i=1$ for  $i\ne\ell$; and
$u_\ell=w_\ell-\frac{w_\ell}{(d^{-1})_{\ell\ell}}e_\ell^T d^{-1} e_\ell=w_\ell-\frac{w_\ell}{(d^{-1})_{\ell\ell}}(d^{-1})_{\ell\ell} = 0$. Thus $u$ vanishes at
$\ell$ and satisfies the weighting equation on all remaining rows, so its restriction to
$[n]\setminus\ell$ is exactly $w_{[n]\setminus\ell}$; reading off coordinate $i$ gives the displayed
formula. The diagonal entry $(d^{-1})_{\ell\ell}$ equals the cofactor ratio
$\det d_{[n]\setminus\ell,[n]\setminus\ell}/\det d$, whose sign is $(-1)^{n-2}/(-1)^{n-1}=-1$ by
Lemma~\ref{lem:inertia}, so it is negative.
\end{proof}

The diagonal sign is fixed but the off-diagonal $(d^{-1})_{i\ell}$ is not, so deleting a negative
$w_\ell$ can move any other weight across zero in either direction. Any repaired proof must control
the \emph{joint} effect of deletions.

\subsection{The corrected algorithm: certified peeling}\label{sec:fix}

The fix follows the diagnosis. The recursion can only \emph{lose} a point that belonged in the peel;
by Lemma~\ref{lem:kkt}, whether that has happened is determined by a marginal $(dp)_k$ exceeding the current level $c$. That test is one matrix-vector product. If it passes,
the output is provably exact and we are done. If it fails, the offending point is precisely one the
recursion should not have dropped, so we re-admit it and climb: move the feasible iterate toward the
new restricted maximizer, releasing any coordinate that would go negative. This is the
Lawson--Hanson active-set safeguard from \S 23.3 of \cite{lawson}, transposed from nonnegative least squares to concave
maximization over the simplex. 

\begin{algorithm}[t]
\caption{Certified peeling.}\label{alg:certified}
\begin{algorithmic}[1]
\REQUIRE strict negative type $d$ on $[n]$, $n\ge2$
\STATE $\mathcal{J}\leftarrow[n]$;\quad $p\leftarrow p^{\mathcal{J}} = \left ( \frac{d_{\mathcal{J},\mathcal{J}}^{-1}1_\mathcal{J}}{1_\mathcal{J}^T d_{\mathcal{J},\mathcal{J}}^{-1}1_\mathcal{J}} \right )^\mathcal{J}$ \COMMENT{one linear solve, Lemma~\ref{lem:restricted}}
\WHILE[\textbf{Phase I}: peeling recursion (no collapse guard, Prop.~\ref{prop:twopos})]{$\min_j p_j<0$}
\STATE $\mathcal{J}\leftarrow\{j:p_j>0\}$;\quad $p\leftarrow p^{\mathcal{J}}$
\ENDWHILE
\STATE $g\leftarrow dp$;\quad $c\leftarrow p^Tg$;\quad $k\leftarrow\argmax_j g_j$ \COMMENT{$c=q(p)$; a stand-alone reuse may write $q$ for $c$}
\WHILE[\textbf{Phase II}: KKT certificate (Lemma~\ref{lem:kkt}); $g_k\le c\Rightarrow p=p_*$]{$g_k>c$}
\STATE $\mathcal{J}\leftarrow\supp(p)\cup\{k\}$;\quad $y\leftarrow p^{\mathcal{J}}$ \COMMENT{\textbf{Phase III}: re-admit the strongest violator}
\WHILE{$\min_{j\in \mathcal{J}}y_j<0$}
\STATE $t\leftarrow\min\{\,p_i/(p_i-y_i):y_i<0\,\}$;\quad $p\leftarrow p+t\,(y-p)$
\STATE $\mathcal{J}\leftarrow\{j:p_j>0\}$;\quad $y\leftarrow p^{\mathcal{J}}$ \COMMENT{drop the zeroed indices; re-solve}
\ENDWHILE
\STATE $p\leftarrow y$;\quad $g\leftarrow dp$;\quad $c\leftarrow p^Tg$;\quad $k\leftarrow\argmax_j g_j$
\ENDWHILE
\RETURN $p$
\end{algorithmic}
\end{algorithm}

Algorithm~\ref{alg:certified} is the original recursion (Phase~I), one matrix-vector product producing a certificate
(Phase~II), and, only on certificate failure, an ascent phase (Phase~III) that re-admits the worst
violator and moves the feasible iterate toward the new restricted maximizer, dropping coordinates
that hit zero: this last is the Lawson-Hanson safeguard transposed to concave maximization over
the simplex. Phase~I is exactly the previously listed algorithm. By
Proposition~\ref{prop:twopos} the surviving set $\{j:(p^{\mathcal{J}})_j>0\}$ always has at least two elements (and this holds not just in exact arithmetic but on every input short of
a fully failed solve) so delete-all-nonpositive never collapses the support and no special guards are required. Two structural facts drive the ascent analysis:
the re-entering coordinate enters with \emph{positive} weight (by Lemma~\ref{lem:schur}, since the point enters
because $(dp)_k>c$), and the value of restricted maximizers strictly increases across re-entries.
Only the \emph{adjoin} direction enters: correctness turns on Lemma~\ref{lem:schur} (how a weight
appears when a point is added), never on its companion Lemma~\ref{lem:deletion} (how weights move when
a point is \emph{removed}). The removal identity is a cost-and-diagnosis tool, driving the 
account of why the published propagation argument
fails (\S\ref{sec:gap}). The certified loop only ever adjoins-then-prunes, and its Phase~I
deletions need nothing beyond ``$|\mathcal{J}|$ decreases,'' with any wrongful deletion repaired downstream.

\begin{theorem}[Correctness and finiteness]\label{thm:correct}
Algorithm~\ref{alg:certified} terminates after finitely many linear solves and returns $p_*$ exactly.
Every working set has $|\mathcal{J}|\ge2$; from the exit of Phase~I onward, $p\in\simplex$; each Phase~III inner loop performs at most $n-1$
solves; and the values observed at successive Phase~II checks strictly increase.
\end{theorem}

\begin{proof}
The assertions of the statement track the phases of the code, with a one-step ascent fact interposed
before Phase~III; write $p^{\mathcal{J}}$, $c_\mathcal{J}$, $A_\mathcal{J}$ as in Lemma~\ref{lem:restricted} and recall
that superscripts extend, while subscripts restrict.

\emph{Phase~I: feasibility and $|\mathcal{J}|\ge2$.} Each round replaces $\mathcal{J}$ by $\mathcal{J}^{+}:=\{j\in \mathcal{J}:(p^{\mathcal{J}})_j>0\}$.
By Proposition~\ref{prop:twopos}, $|\mathcal{J}^{+}|\ge2$, and with a margin: the off-max positive mass is
$\ge\tfrac12$, so the second-largest coordinate of $p^{\mathcal{J}}$ clears the deletion threshold by
$1/(2(|\mathcal{J}|-1))$. So the support never falls below two points, and $\mathcal{J}$ decreases through
strict-negative-type principal submatrices of size $\ge2$. When the loop does not exit, $\mathcal{J}^{+}\subsetneq \mathcal{J}$ (the
exit test $(p^{\mathcal{J}})_\mathcal{J}\ge0$ failed, so some coordinate was deleted), so Phase~I halts in $\le n-2$
rounds; if $|\mathcal{J}|=2$, the weighting of $\begin{psmallmatrix}0&a\\a&0\end{psmallmatrix}$ is
$w_\mathcal{J}=(1,1)^T/a>0$, forcing the exit. The loop therefore exits at some $p=p^{\mathcal{J}}\in\simplex$
with $q(p)=c_\mathcal{J}>0$ (Lemma~\ref{lem:restricted}).

\emph{Phase~II is a certificate.} At a check (i.e., an evaluation of the Boolean predicate $g_k > c$), $p=p^{\mathcal{J}}$ for the current set $\mathcal{J}$, so $(dp)_j=c_\mathcal{J}$ for every $j\in \mathcal{J}$ and the recorded level is $c=p^Tdp=c_\mathcal{J}>0$
(Lemma~\ref{lem:restricted}). Put $\mathcal{I}:=\supp p\subseteq \mathcal{J}$; the inclusion may be strict, if a
coordinate of $p^{\mathcal{J}}$ vanishes exactly. Since $q(e_j)=0<c=q(p)$ we have $|\mathcal{I}|\ge2$; and $p\in A_\mathcal{I}$
with $(dp)_j=c$ on $\mathcal{I}$ is stationarity for $q|_{A_\mathcal{I}}$, so Lemma~\ref{lem:restricted} identifies
$p=p^\mathcal{I}$ and $c=c_\mathcal{I}$. In particular \emph{the value at a check depends only on the support}, which is the
fact the termination argument will use. By the KKT test \eqref{eq:kkt}/Lemma~\ref{lem:kkt}, $p=p_*$
iff $\max_k(dp)_k\le c$, which is exactly the test $g_k\le c$; so acceptance returns $p_*$
exactly, and every non-terminal check has a strict violator $(dp)_k>c$ --- which forces $k\notin \mathcal{I}$,
since $(dp)_j=c$ on $\mathcal{I}$.

\emph{One ascent step.} Let $|\mathcal{J}|\ge2$, $p\in A_\mathcal{J}\cap\simplex$, and $y:=p^{\mathcal{J}}$. By
\eqref{eq:exactTaylor}, $\chi(t):=q(p+t(y-p))$ is a parabola with leading coefficient
$(y-p)^Td(y-p)$, which is $\le0$ by strict negative type ($1^T(y-p)=0$) and $<0$ unless $y=p$; and
$\chi'(1)=2(dy)^T(y-p)=2c_\mathcal{J}\,1^T(y-p)=0$, since $y-p$ is supported in $\mathcal{J}$ and $(dy)_j=c_\mathcal{J}$
there (Lemma~\ref{lem:restricted}). Hence $t=1$ maximizes $\chi$ over the whole affine line, and
$\chi$ is nondecreasing on $[0,1]$, strictly increasing there unless $y=p$. Every Phase~III move (i.e., the truncating steps, with $t\in[0,1)$, and the terminal assignment $p\leftarrow y$,
which amounts to $t=1$) is
of this form (as the next paragraph establishes), so $q(p)$ never decreases within a re-entry of a point into $\mathcal{J}$. 


\emph{Phase~III: feasibility, $|\mathcal{J}|\ge2$, and the inner bound.} A failing check sets
$\mathcal{J}\leftarrow \mathcal{I}\cup\{k\}$, so $|\mathcal{J}|=|\mathcal{I}|+1\ge3$, using $|\mathcal{I}|\ge2$ and $k\notin \mathcal{I}$ from Phase~II; this also forms $y:=p^{\mathcal{J}}$. We claim, by induction over the inner passes, that at the start of each pass
$p\in A_\mathcal{J}\cap\simplex$, $q(p)\ge c_\mathcal{I}>0$, and $|\mathcal{J}|\ge2$, so that $y=p^{\mathcal{J}}$ is defined
(Lemma~\ref{lem:restricted}). All three hold at the first pass. At a pass, if $y_\mathcal{J}\ge0$ the inner loop
ends and the algorithm sets $p\leftarrow y\in A_\mathcal{J}\cap\simplex$. Otherwise it moves
$p\leftarrow p+t(y-p)$ with $t=\min\{p_i/(p_i-y_i):y_i<0\}$; every such $i$ has $p_i\ge0>y_i$, so
$p_i-y_i>p_i\ge0$ and $t\in[0,1)$, and truncating at the first zero-crossing preserves $1^Tp=1$ and
$p\ge0$, i.e.\ $p\in\simplex$. The minimizing index is zeroed, so $\mathcal{J}\leftarrow\{j:p_j>0\}$ removes
at least one index, 
while $\supp p\subseteq\mathcal{J}$
keeps $p\in A_\mathcal{J}$ for the new $\mathcal{J}$. By the ascent step above, $q(p)\ge c_\mathcal{I}>0$ persists; since
$q(e_j)=0$, this forbids $|\mathcal{J}|=1$, so the new working set again has $|\mathcal{J}|\ge2$ and the induction
closes. As $|\mathcal{J}|$ falls by $\ge1$ per truncating pass and cannot fall below $2$ --- at $|\mathcal{J}|=2$ the
two-point weighting is positive and the exit is forced --- there are at most $|\mathcal{I}|-1\le n-2$ truncating
passes, so the inner loop performs at most $n-1$ solves (one on entry, one per truncating pass) before
$y_\mathcal{J}\ge0$. The re-solve is what makes the terminal assignment $p\leftarrow y$ correct: the loop
exits only when the \emph{currently recomputed} $y=p^{\mathcal{J}}$ satisfies $y_\mathcal{J}\ge0$, so the assigned
$p$ is feasible and is the restricted maximizer of its own support. (By comparison, walking toward an
outdated $y$ and assigning it would reinstate the very negatives the drops removed.)

\emph{Strict increase across re-entries.} At a failing check $p=p^\mathcal{I}$ with $(dp)_k>c=c_\mathcal{I}$, two things
must be checked: that the enlarged maximizer is strictly better, and that the first inner pass, when
there is one, actually moves toward it. The first is immediate: $(dp)_j=c$ on $\mathcal{I}$ while $(dp)_k>c$, so
$p$ is \emph{not} stationary on the larger hull $A_{\mathcal{I}\cup \{k\}}\supsetneq A_\mathcal{I}$, and its maximizer
therefore strictly improves, $q(y)>q(p)$ for $y:=p^{\mathcal{I}\cup \{k\}}$ (Lemma~\ref{lem:restricted}). The second
is where Lemma~\ref{lem:schur} is essential: because $(dp)_k>c=c_\mathcal{I}$, the adjoined coordinate enters
with \emph{strictly positive} unnormalized weight, $(w_{\mathcal{I}\cup\{k\}})_k>0$; and since
$1_{\mathcal{I}\cup \{k\}}^Tw_{\mathcal{I}\cup \{k\}}>0$ (Lemma~\ref{lem:restricted}) the normalization
$p^{\mathcal{I}\cup \{k\}}=\bigl(w_{\mathcal{I}\cup \{k\}}/1_{\mathcal{I}\cup \{k\}}^Tw_{\mathcal{I}\cup \{k\}}\bigr)^{\mathcal{I}\cup \{k\}}$ preserves that sign, so
$y_k>0$. Hence $k$ is \emph{not} among the coordinates $\{i:y_i<0\}$ that limit the step, which all
lie in $\mathcal{I}$ (where $p_i>0$), so $t=\min\{p_i/(p_i-y_i):y_i<0\}>0$; were $y_k<0$ instead, $p_k=0$ would force $t=0$, an immediate re-drop of $k$, no gain in $q$, and a stall. (The same conclusion holds at every
later pass, where $\mathcal{J}=\supp p$ by construction, so no truncating pass has $t=0$.) Either way the value strictly increases
across the re-entry: if $y_\mathcal{J}\ge0$ already at the first pass, the algorithm assigns $p\leftarrow y$
and the value jumps to $q(y)>q(p)$; otherwise the first pass moves with $t>0$ and the ascent step,
strictly increasing on $[0,1]$, gives $q(p+t(y-p))>q(p)$. Monotonicity carries the gain to the next
check.

\emph{Termination and exactness.} Each Phase~II check records the value $c_\mathcal{I}$ of its current support
$\mathcal{I}$; by Phase~II that value is determined by the support alone, and by the previous paragraph these
values strictly increase across re-entries. Hence no support recurs, and as there are finitely many
supports the outer loop halts. On the halting check the KKT test passes, so by Phase~II the returned
$p$ equals $p_*$.
\end{proof}

\begin{remark}[Cost; floating point]\label{rem:cost}
If the certificate passes at the first check, which is common, the total cost is the recursion plus one $O(n\,|\supp p|)$ product. On re-entry, solves are on supports of size $|\supp p|+1\ll n$ and can reuse factorizations via
bordered/deletion updates (Lemmas~\ref{lem:schur} and \ref{lem:deletion}) at $O(|\mathcal{J}|^{2})$ per change.
The number of checks is finite but not polynomially bounded by this argument (as is typical for
active-set methods: on the other hand, Claude Fable reported at most one re-entry and two extra solves on $40000$ instances). 
\end{remark}

\subsection{The circumcenter dictionary}\label{sec:dictionary}

Everything so far was algebraic. This section installs a geometric picture that informs what follows (cf. \cite{bc,clarkson}). Under Schoenberg's embedding  \cite{schoenberg}, the quadratic entropy, the weighting equation, and the
deletion rule all become statements about \emph{circumcenters and enclosing balls}. Once the
dictionary is in place, ``the recursion stalls'' reads as ``a discarded point ended up outside the
final ball,'' and the failure envelopes of the next section become statements about how far a center
drifts as points are dropped.

By Schoenberg's theorem, any negative-type $d$ admits an embedding
$d_{ij}=\lVert\phi^{(i)}-\phi^{(j)}\rVert^{2}$. Write $\Phi$ for the matrix with $i$th column $\phi^{(i)}$, so that $z \equiv z(p) := p^T \Phi = \sum_i p_i \phi^{(i)}$. Two identities, valid for $1^Tp=1$, govern the geometry:
\begin{proposition}[Bias and variance]\label{prop:biasvar}
    \begin{equation}\label{eq:biasvar}
(dp)_j=\lVert\phi^{(j)}-z(p) \rVert^{2}+\sum\nolimits_ip_i\lVert\phi^{(i)}-z(p)\rVert^{2};
\end{equation}
and, if $1^T v = 0$,
\begin{equation}\label{eq:biasvar2}
v^Tdv=-2\Bigl\lVert\sum\nolimits_iv_i\phi^{(i)}\Bigr\rVert^{2}.
\end{equation}
\end{proposition}

\begin{proof}
We first show \eqref{eq:biasvar}:
\begin{align}
    (dp)_j & = \sum_k \|\phi^{(j)}-\phi^{(k)}\|^2 p_k \nonumber \\
    & = \sum_k \|(\phi^{(j)}-z(p))+(z(p)-\phi^{(k)})\|^2 p_k \nonumber \\
    & = \|\phi^{(j)}-z(p)\|^2 + \sum_k \|\phi^{(k)}-z(p)\|^2 p_k + 2 \left \langle \phi^{(j)} - z(p), \sum\nolimits_k (z(p)-\phi^{(k)})p_k \right \rangle \nonumber \\
    & = \|\phi^{(j)}-z(p)\|^2 + \sum_k \|\phi^{(k)}-z(p)\|^2 p_k \nonumber
\end{align}
with the last equality by definition of $z(p)$. To show \eqref{eq:biasvar2}:
\begin{align}
    v^T d v & = \sum_{i,j} v_i v_j \|\phi^{(i)}-\phi^{(j)}\|^2 \nonumber \\
    & = \sum_{i,j} v_i v_j \left ( \|\phi^{(i)}\|^2 + \|\phi^{(j)}\|^2 - 2\langle \phi^{(i)}, \phi^{(j)} \rangle \right ) \nonumber \\
    & = \left (\sum\nolimits_i v_i \|\phi^{(i)}\|^2 \right ) \cdot \left (\sum\nolimits_j v_j \right ) + \left (\sum\nolimits_i v_i \right ) \cdot \left (\sum\nolimits_j v_j \|\phi^{(j)}\|^2 \right ) - 2 \sum_{i,j} v_i v_j \left \langle \phi^{(i)}, \phi^{(j)} \right \rangle \nonumber \\
    & = - 2 \left \langle \sum\nolimits_i v_i \phi^{(i)}, \sum\nolimits_j v_j \phi^{(j)} \right \rangle \nonumber \\
    & = -2\Bigl\lVert\sum\nolimits_iv_i\phi^{(i)}\Bigr\rVert^{2}. \nonumber
\end{align}
\end{proof}

The Pythagorean flavor of \eqref{eq:biasvar} indicates orthogonality that will be manifested below. Note also that \eqref{eq:biasvar2} shows that strict negative type is exactly \emph{affine independence} of $\{\phi^{(i)}\}$.

\begin{lemma}[Dictionary]\label{lem:dict}
Let $|\mathcal{J}|\ge2$ and $$F_\mathcal{J}:=\aff\{\phi^{(j)}\}_{j\in \mathcal{J}}.$$
\begin{enumerate}
\item[(i)] There is a unique $z_\mathcal{J}\in F_\mathcal{J}$ equidistant from $\{\phi^{(j)}\}_{j\in \mathcal{J}}$ (the \emph{in-hull
circumcenter}); its barycentric coordinates are $p^\mathcal{J}$, and with $R_\mathcal{J}$ the common distance,
\begin{equation}\label{eq:c2r2}
c_\mathcal{J}=2R_\mathcal{J}^{2}.
\end{equation}
\item[(ii)] For $k\notin \mathcal{J}$,
\begin{equation}\label{eq:overshoot}
(dp^\mathcal{J})_k-c_\mathcal{J}=\lVert\phi^{(k)}-z_\mathcal{J}\rVert^{2}-R_\mathcal{J}^{2};
\end{equation}
combined with
Lemma~\ref{lem:schur}, the weight of $k$ upon adjunction is positive/zero/negative according as
$\phi^{(k)}$ lies respectively outside/on/inside the circumsphere of $\mathcal{J}$; equivalently, for $a\in \mathcal{J}$ ($|\mathcal{J}|\ge3$),
$\operatorname{sign}((w_\mathcal{J})_a)=\operatorname{sign}(\lVert\phi^{(a)}-z_{\mathcal{J}\setminus a}\rVert^{2}-R_{\mathcal{J}\setminus a}^{2})$.
\item[(iii)] If $p^\mathcal{J}\ge0$ then $B(z_\mathcal{J},R_\mathcal{J})$ is the minimum enclosing ball (MEB) of
$\{\phi^{(j)}\}_{j\in \mathcal{J}}$, $c_\mathcal{J}$ is the \emph{constrained} maximum of $q$ over $\Delta_\mathcal{J}$, and for every
$k$, $(dp^\mathcal{J})_k-c_\mathcal{J}$ measures how far outside this ball $\phi^{(k)}$ lies, in squared units.
\end{enumerate}
\end{lemma}

\begin{proof}
(i) By the first identity in \eqref{eq:biasvar}, the marginal $(dp)_j$ is constant over $j\in \mathcal{J}$ if
and only if the distance $\lVert\phi^{(j)}-z(p)\rVert$ is constant over $\mathcal{J}$, i.e.\ $z(p)$ is equidistant
from the support points. Stationarity of $q$ on $A_\mathcal{J}$ is $d_{\mathcal{J},\mathcal{J}}p_\mathcal{J}=c1_\mathcal{J}$
(Lemma~\ref{lem:restricted}), which is precisely this constancy; so the barycenter $z_\mathcal{J}$ of $p^\mathcal{J}$ is
an equidistant point lying in $F_\mathcal{J}$. Such a point is unique: the set of points equidistant from
$\{\phi^{(j)}\}_{j\in \mathcal{J}}$ is an affine flat orthogonal to $F_\mathcal{J}$, and it meets $F_\mathcal{J}$ in a single point.
With $R_\mathcal{J}$ the common distance, the second summand of \eqref{eq:biasvar} is $R_\mathcal{J}^2$, so
$(dp^\mathcal{J})_j=2R_\mathcal{J}^{2}$, giving $c_\mathcal{J}=2R_\mathcal{J}^2$. 

(ii) Equation~\eqref{eq:overshoot} is the evaluation of \eqref{eq:biasvar} at $j=k\notin \mathcal{J}$ (the constant summand is again $R_\mathcal{J}^2$); pairing it
with Lemma~\ref{lem:schur} (which ties the sign of $k$'s adjoined weight to the sign of $(dp^\mathcal{J})_k-c_\mathcal{J}$)
gives the outside/on/inside trichotomy. For a point $a\in \mathcal{J}$ with $|\mathcal{J}|\ge3$ (so $\mathcal{J}\setminus a$ has
$\ge2$ points and $z_{\mathcal{J}\setminus a},R_{\mathcal{J}\setminus a}$ are defined), apply Lemma~\ref{lem:schur} with
base set $\mathcal{J}\setminus a$ and adjoined point $a$: since $(\mathcal{J}\setminus a)\cup a=\mathcal{J}$ and a set's weighting is
unique, the coordinate it adjoins is $(w_\mathcal{J})_a$, so
$\operatorname{sign}((w_\mathcal{J})_a)=\operatorname{sign}\bigl((dp^{\mathcal{J}\setminus a})_a-c_{\mathcal{J}\setminus a}\bigr)$.
Applying \eqref{eq:overshoot} with $\mathcal{J}\setminus a$ in place of $\mathcal{J}$ and $a$ in place of $k$ gives
$(dp^{\mathcal{J}\setminus a})_a-c_{\mathcal{J}\setminus a}=\lVert\phi^{(a)}-z_{\mathcal{J}\setminus a}\rVert^{2}-R_{\mathcal{J}\setminus a}^{2}$,
whence $\operatorname{sign}((w_\mathcal{J})_a)=\operatorname{sign}\bigl(\lVert\phi^{(a)}-z_{\mathcal{J}\setminus a}\rVert^{2}
-R_{\mathcal{J}\setminus a}^{2}\bigr)$: the point $a$ carries positive weight in $\mathcal{J}$ exactly when $\phi^{(a)}$ lies
strictly outside the circumsphere of the remaining support $\mathcal{J}\setminus a$.

(iii) If $p^\mathcal{J}\ge0$ then $z_\mathcal{J}$ lies in the convex hull of $\{\phi^{(j)}\}_{j \in \mathcal{J}}$, and we claim $B(z_\mathcal{J},R_\mathcal{J})$ is the minimum
enclosing ball. For any candidate center $y$, convexity gives a support point $\phi^{(j)}$ on the far side of $z_\mathcal{J}$,
$\langle\phi^{(j)}-z_\mathcal{J},\,y-z_\mathcal{J}\rangle\le0$.
By the polarization identity $2 \langle \alpha, \beta \rangle = \| \alpha \|^2 + \| \beta \|^2 - \| \alpha - \beta \|^2$, $2 \langle \phi^{(j)} - z_\mathcal{J}, y-z_\mathcal{J} \rangle = \| \phi^{(j)} - z_\mathcal{J} \|^2 + \| y - z_\mathcal{J} \|^2 - \| \phi^{(j)} - y \|^2$, so if $\langle \phi^{(j)} - z_\mathcal{J}, y-z_\mathcal{J} \rangle \le 0$, then
$\lVert\phi^{(j)}-y\rVert^2\ge\lVert\phi^{(j)}-z_\mathcal{J}\rVert^2+\lVert y-z_\mathcal{J}\rVert^2=R_\mathcal{J}^2+\lVert y-z_\mathcal{J}\rVert^2$,
so no other center encloses $\mathcal{J}$ in a smaller radius. That $c_\mathcal{J}$ is the constrained maximum of $q$ over
$\Delta_\mathcal{J}$ is Lemma~\ref{lem:kkt} applied on the face, and the last clause is the identity
$(dp^\mathcal{J})_k-c_\mathcal{J}=\lVert\phi^{(k)}-z_\mathcal{J}\rVert^{2}-R_\mathcal{J}^{2}$, read as the squared distance by which $\phi^{(k)}$
overshoots the sphere $\partial B(z_\mathcal{J},R_\mathcal{J})$.
\end{proof}

In the embedding, the Phase I recursion is: \emph{compute the in-hull circumcenter; discard negative-barycentric
points; stop when the center enters the convex hull; return its circumball, i.e., the MEB of surviving points}, while $c_*$ is twice the squared MEB radius of \emph{all} points. Failure means a discarded point ended up outside the final ball.

\subsection{Failure envelopes, per-iteration budgets, and the degeneracy bound}\label{sec:envelope}

Knowing the recursion can fail, we now bound \emph{how much} it can lose --- and do so using only
quantities the run already produced. Three bounds appear in Theorem \ref{thm:envelope}. The first (part (i)) brackets $c_*$ between
two free numbers: the first-solve value $c_{(1)}$ from above and $\chat+2V$ from the terminal certificate.
The second (part (ii)) opens up $c_{(1)}-\chat$ as a sum of nonnegative per-round ``budgets'' $\Delta_{(t)}$,
each equal --- through the dictionary of \S \ref{sec:dictionary} --- to twice the squared distance the
circumcenter drifts when round $t$ deletes. The third (part (iii)) turns those budgets into a
per-point bound on how far a deleted point can end up outside the terminal ball.

Let $$V_k:=[(d\hat p)_k-\chat]_+$$ be the per-point terminal violations and 
\begin{equation}\label{eq:V}
V:=\max_kV_k.
\end{equation}

\begin{theorem}[Computable failure envelope]\label{thm:envelope}

\

\begin{enumerate}
\item[(i)] $$\chat\le c_*\le\min\{c_{(1)},\ \chat+2V\};$$ both bounds are available from quantities the
recursion already computes (plus one matvec for $V$).
\item[(ii)] With $v_{(t)}:=p^{\mathcal{J}_{(t+1)}}-p^{\mathcal{J}_{(t)}}$, $\Delta_{(t)}:=c_{(t)}-c_{(t+1)}$, and $z_{(t)} = (p^{\mathcal{J}_{(t)}})^T\Phi = \sum_{j \in \mathcal{J}_{(t)}} (p^{\mathcal{J}_{(t)}})_j \phi^{(j)}$ the circumcenter of $\{ \phi^{(j)} \}_{j \in \mathcal{J}_{(t)}}$, 
\[
\Delta_{(t)}=-v_{(t)}^Tdv_{(t)}=2\lVert z_{(t+1)}-z_{(t)}\rVert^{2}=2\bigl(R_{(t)}^{2}-R_{(t+1)}^{2}\bigr)\ \ge0,
\]
strict when round $t$ deletes; $z_{(t+1)}$ is the orthogonal projection of $z_{(t)}$ onto $F_{\mathcal{J}_{(t+1)}}$;
and $c_{(1)}-\chat=\sum_{t<T}\Delta_{(t)}$, so $c_*-\chat\le\sum_{t<T}\Delta_{(t)}$: the failure is at most twice the
total squared circumcenter drift, with $\Delta_{(t)}$ as iteration $t$'s exact budget.
\item[(iii)] If $k$ is deleted at round $t$, then immediately afterwards
$(dp^{\mathcal{J}_{(t+1)}})_k-c_{(t+1)}\le\Delta_{(t)}+\sqrt{\Delta_{(t)}c_{(t)}}$. At termination, with
$D_{(t)}:=c_{(t)}-\chat=\sum_{t\le s< T}\Delta_{(s)}$,
\[
(d\hat p)_k-\chat\;\le\;D_{(t)}+\sqrt{c_{(t)}\,D_{(t)}}\;
\;\le\;D_{(1)}+\sqrt{c_{(1)}D_{(1)}},
\]
with the last expression independent of the round count $T$.
\end{enumerate}
\end{theorem}

\begin{proof}
(i) The lower bound $\chat\le c_*$ is feasibility of $\hat p$. The upper bound $c_*\le c_{(1)}$ is the
affine relaxation: maximizing $q$ over the whole affine hull (dropping nonnegativity) can only raise
the maximum, and that relaxed maximum is $c_{(1)}$ (cf. \eqref{eq:bound1}). The upper bound
$c_*\le\chat+2V$ is Corollary~\ref{cor:apost}.

(ii) Fix a deleting round $t$, and expand $q$ at $p^{\mathcal{J}_{(t)}}$ along the increment $v_{(t)}=p^{\mathcal{J}_{(t+1)}}-p^{\mathcal{J}_{(t)}}$
using the exact Taylor identity \eqref{eq:exactTaylor}. The gradient $2\,dp^{\mathcal{J}_{(t)}}$ equals $2c_{(t)}$ componentwise on $\mathcal{J}_{(t)}$, and
$\supp v_{(t)}\subseteq \mathcal{J}_{(t)}$ with $1^Tv_{(t)}=0$, so the linear term vanishes and only the quadratic
survives: $c_{(t+1)}=c_{(t)}+v_{(t)}^Tdv_{(t)}$, i.e.\ $\Delta_{(t)}=c_{(t)}-c_{(t+1)}=-v_{(t)}^Tdv_{(t)}\ge0$. Translating
the quadratic through the embedding, $v_{(t)}^Tdv_{(t)}=-2\lVert\sum_i v_{(t)i}\phi^{(i)}\rVert^2
=-2\lVert z_{(t+1)}-z_{(t)}\rVert^2$ (the second identity of \eqref{eq:biasvar}), which gives the middle
equality; and since $c=2R^2$ by \eqref{eq:c2r2} it also equals $2(R_{(t)}^2-R_{(t+1)}^2)$. Geometrically
$z_{(t)}$ is equidistant, at distance $R_{(t)}$, from every point of $\mathcal{J}_{(t+1)}\subseteq \mathcal{J}_{(t)}$, so it lies on the
flat through $z_{(t+1)}$ orthogonal to $F_{\mathcal{J}_{(t+1)}}$; the Pythagorean theorem gives
$R_{(t)}^2=R_{(t+1)}^2+\lVert z_{(t)}-z_{(t+1)}\rVert^2$, i.e.\ $z_{(t+1)}$ is the orthogonal projection of $z_{(t)}$
onto $F_{\mathcal{J}_{(t+1)}}$. Summing $\Delta_{(t)}$ over $t<T$ telescopes to $c_{(1)}-c_{(T)}=c_{(1)}-\chat$, and combining
with (i) gives $c_*-\chat\le c_{(1)}-\chat=\sum_t\Delta_{(t)}$.

(iii) Suppose $k$ is deleted at round $t$. Just before deletion $\lVert\phi^{(k)}-z_{(t)}\rVert=R_{(t)}$ (it lies
on the round-$t$ circumsphere), with $R_{(t)}=\sqrt{c_{(t)}/2}$ since $c_{(t)}=2R_{(t)}^2$ by \eqref{eq:c2r2}. 
By part~(ii) the circumcenter increments are mutually orthogonal. Indeed, write $\vec F_\mathcal{J} := F_\mathcal{J} - F_\mathcal{J}$ for the
direction space of $F_\mathcal{J}$ and fix $t\le s<s'<T$. Since $z_{(s+1)}$ is the orthogonal projection of
$z_{(s)}$ onto $F_{\mathcal{J}_{(s+1)}}$, the residual obeys $z_{(s+1)}-z_{(s)}\in\vec F_{\mathcal{J}_{(s+1)}}^{\perp}$.
Now $\mathcal{J}_{(s'+1)}\subseteq \mathcal{J}_{(s')}\subseteq \mathcal{J}_{(s+1)}$, and $F$
is monotone in its index set, so we have both
$z_{(s')}\in F_{\mathcal{J}_{(s')}}\subseteq F_{\mathcal{J}_{(s+1)}}$ and $z_{(s'+1)}\in F_{\mathcal{J}_{(s'+1)}}\subseteq
F_{\mathcal{J}_{(s+1)}}$. Hence $z_{(s'+1)}-z_{(s')}\in\vec F_{\mathcal{J}_{(s+1)}}$, and the two increments are
orthogonal. Therefore, the Pythagorean theorem telescopes the displacement from round $t$ to the terminal center
$\hat z:=z_{(T)}$:
\[
\lVert\hat z-z_{(t)}\rVert^2=\sum_{s=t}^{T-1}\lVert z_{(s+1)}-z_{(s)}\rVert^2=\sum_{s\ge t}\tfrac{\Delta_{(s)}}{2}
=\tfrac{c_{(t)}-\chat}{2}=\tfrac{D_{(t)}}{2},
\]
using both $\Delta_{(s)}=2\lVert z_{(s+1)}-z_{(s)}\rVert^2$ and $\sum_{s\ge t}\Delta_{(s)}=c_{(t)}-\chat$ from part~(ii). 

A
single triangle inequality gives 
$$\lVert\phi^{(k)}-\hat z\rVert\le \| \phi^{(k)} - z_{(t)} \| + \| \hat z - z_{(t)} \| = R_{(t)}+\lVert\hat z-z_{(t)}\rVert
=\sqrt{c_{(t)}/2}+\sqrt{D_{(t)}/2};$$ 
substituting into \eqref{eq:overshoot} with $\hat R^2=\chat/2$ and expanding gives the terminal bound
\begin{align}
    (d\hat p)_k - \chat & =\lVert\phi^{(k)}-\hat z\rVert^2-\hat R^2 \nonumber \\
    & \le \left ( \sqrt{c_{(t)}/2}+\sqrt{D_{(t)}/2} \right )^2 - \frac{\chat}{2} \nonumber \\
    & = \frac{c_{(t)}}{2} + \frac{D_{(t)}}{2} + \sqrt{c_{(t)} D_{(t)}} - \frac{\chat}{2} \nonumber \\
    & = D_{(t)}+\sqrt{c_{(t)} D_{(t)}}. \nonumber
\end{align}

To show that $(dp^{\mathcal{J}_{(t+1)}})_k-c_{(t+1)}\le\Delta_{(t)}+\sqrt{\Delta_{(t)}c_{(t)}}$, stop the same computation as above at $z_{(t+1)}$, observing three facts. First, the increment $v_{(t)}:=p^{\mathcal{J}_{(t+1)}}-p^{\mathcal{J}_{(t)}}$ has $1 ^T v_{(t)}=0$
and $\supp v_{(t)}\subseteq \mathcal{J}_{(t)}$, so its image under the embedding is
\[
\sum_i(v_{(t)})_i\phi^{(i)}=\sum_i(p^{\mathcal{J}_{(t+1)}})_i\phi^{(i)}-\sum_i(p^{\mathcal{J}_{(t)}})_i\phi^{(i)}=z_{(t+1)}-z_{(t)},
\]
with the last equality holding since the barycentric
coordinates of $z_\mathcal{J}$ are $p^{\mathcal{J}}$ (Lemma~\ref{lem:dict}(i)).
Therefore by \eqref{eq:biasvar2} we obtain $\Delta_{(t)}=-v_{(t)}^T dv_{(t)}=2\bigl\lVert z_{(t+1)}-z_{(t)}\bigr\rVert^{2}$,
and in turn
\begin{equation}
\bigl\lVert z_{(t+1)}-z_{(t)}\bigr\rVert=\sqrt{\Delta_{(t)}/2}.
    \tag{a}
\end{equation}
Second, by \eqref{eq:c2r2} we have
\begin{equation}
R_{(t+1)}^2 = c_{(t+1)}/2 = (c_{(t)} - \Delta_{(t)})/2.
    \tag{b}
\end{equation}
Third, since $k \in \mathcal{J}_{(t)}$ by hypothesis and again by \eqref{eq:c2r2}, we have 
\begin{equation}
\| \phi^{(k)} - z_{(t)} \| = R_{(t)} = \sqrt{c_{(t)}/2}.
    \tag{c}
\end{equation}
By (a) and (c), 
$$\lVert\phi^{(k)}-z_{(t+1)}\rVert\le \|\phi^{(k)} - z_{(t)} \| + \| z_{(t+1)} - z_{(t)} \| =\sqrt{c_{(t)}/2}+\sqrt{\Delta_{(t)}/2}.$$
Now (b) yields (similarly to above)
\begin{align}
    (dp^{\mathcal{J}_{(t+1)}})_k-c_{(t+1)}
&=\lVert\phi^{(k)}-z_{(t+1)}\rVert^{2}-R_{(t+1)}^{2} \nonumber \\
    & \le \left ( \sqrt{c_{(t)}/2}+\sqrt{\Delta_{(t)}/2} \right )^2 - \frac{c_{(t)}-\Delta_{(t)}}{2} \nonumber \\
    & = \frac{c_{(t)}}{2} + \frac{\Delta_{(t)}}{2} + \sqrt{c_{(t)} \Delta_{(t)}} - \frac{c_{(t)}}{2} + \frac{\Delta_{(t)}}{2} \nonumber \\
    & = \Delta_{(t)}+\sqrt{c_{(t)} \Delta_{(t)}}. \nonumber
\end{align}

Finally, we need to show that $D_{(t)}+\sqrt{c_{(t)}\,D_{(t)}}\;\le\;D_{(1)}+\sqrt{c_{(1)}D_{(1)}}$. Writing $\Psi(c,D):=D+\sqrt{cD}$, we have $\partial_c\Psi=\tfrac12\sqrt{D/c}\ge0$ and $\partial_D\Psi=1+\tfrac12\sqrt{c/D}>0$. From the proof of (ii) we have that $\Delta_{(t)} \ge 0$, so $c_{(1)} \ge \dots \ge c_{(T)} \equiv \chat$: in particular, $c_{(t)} \le c_{(1)}$ and $D_{(t)} = c_{(t)} - \chat \le c_{(1)} - \chat$. Combining these observations, we have that 
$$(d\hat p)_k - \chat \le \Psi(c_{(t)},D_{(t)}) \le \Psi(c_{(1)},c_{(1)} - \chat) = \Psi(c_{(1)},D_{(1)}) = D_{(1)} + \sqrt{c_{(1)} D_{(1)}}.$$
\end{proof}

\begin{theorem}[Failure requires a small strictness margin]\label{thm:lam2}
$c_{(1)}\le\bar d+\lVert Hd1\rVert^{2}/(n^{2}\lambda_2)$ with $\bar d=1^Td1/n^{2}$; hence
$c_*-\chat\le\bar d-\chat+\lVert Hd1\rVert^{2}/(n^{2}\lambda_2)$: uniformly strict families have
uniformly bounded failure.
\end{theorem}

\begin{proof}
Write any $p\in\simplex$ as $p=\frac{1}{n}1+v$ with $v\perp 1$ (so $v=Hp$). Expanding $q(p)=p^Tdp$ and
using $1^Td1=n^2\bar d$, the cross term is $\frac{2}{n}1^Tdv=\tfrac2n(d1)^Tv=\tfrac2n(Hd1)^Tv$
(as $v\perp 1$, or for that matter as $H^2 = H$, so that $Hv = v$), and the quadratic term is $v^Tdv\le-\lambda_2\lVert v\rVert^2$ by definition of the
strictness margin. Hence
$q(p)\le\bar d+\tfrac2n\lVert Hd1\rVert\,\lVert v\rVert-\lambda_2\lVert v\rVert^2$ by Cauchy-Schwarz.
The right side is a downward parabola in $\lVert v\rVert$; maximizing over $\lVert v\rVert\ge0$ gives
$\bar d+\lVert Hd1\rVert^2/(n^2\lambda_2)$ by (T2). Since $c_{(1)}=\max_p q(p)$ over the affine hull, this bounds
$c_{(1)}$, and subtracting $\chat$ gives the failure bound.
\end{proof}

\subsection{A violation-weighted a posteriori certificate}\label{sec:sharp}

\begin{theorem}[Violation-weighted certificate]\label{thm:sharp}
With $\varphi_\lambda(t):=t^{2}/\lambda$ for $0\le t\le\lambda$ and $\varphi_\lambda(t):=2t-\lambda$
for $t>\lambda$:
\[
c_*-\chat\;\le\;\min\Bigl\{\,2V,\ \sum_{k\notin\hat{\mathcal{J}}}\varphi_{\lambda_2}(V_k)\Bigr\}.
\]
\end{theorem}

\begin{proof}
Let $v=p_*-\hat p$ and $g=d\hat p$. Now 
\begin{align}
    g^Tv & = \sum_k g_k v_k \nonumber \\
    & = \sum_k (g_k - \chat)v_k \nonumber \\
    & = \sum_{k \in \hat{\mathcal J}} (g_k - \chat)v_k + \sum_{k \not \in \hat{\mathcal J}} (g_k - \chat)v_k \nonumber \\
    & = \sum_{k \in \hat{\mathcal J}} (\chat - \chat)v_k + \sum_{k \not \in \hat{\mathcal J}} (g_k - \chat)v_k \nonumber \\
    & = \sum_{k \not \in \hat{\mathcal J}} (g_k - \chat)v_k \nonumber \\
    & \le \sum_{k \not \in \hat{\mathcal J}} [g_k - \chat]_+ v_k \nonumber \\
    & = \sum_{k \not \in \hat{\mathcal J}} V_k v_k \nonumber \\
    & = \sum_{k \not \in \hat{\mathcal J}} V_k (p_*)_k, \tag{*}\label{eq:violationCertPart1}
\end{align}
where the second equality is because $\chat \sum_k v_k = 0$; the fourth equality is because $g_k = (d\hat p)_k = \chat$ for $k \in \hat{\mathcal J}$, and the last equality is because $\hat p_k = 0$ for $k \not \in \hat{\mathcal J}$, so that $v_k = (p_*)_k$  for $k \not \in \hat{\mathcal J}$. This last observation also yields 
\begin{equation}
    v^T d v \le -\lambda_2 \| v \|^2 \le -\lambda_2 \sum_{k \not \in \hat{\mathcal J}} (p_*)_k^2. \tag{**}\label{eq:violationCertPart2}
\end{equation}

Combining \eqref{eq:violationCertPart1} and \eqref{eq:violationCertPart2} with \eqref{eq:exactTaylor}, we get 
\begin{align}
    c_* - \chat & = 2 g^T v +v^Tdv \nonumber \\
    & \le \sum_{k \not \in \hat{\mathcal J}} \left ( 2V_k(p_*)_k-\lambda_2(p_*)_k^2 \right ). \tag{***}\label{eq:violationCertPart3}
\end{align}
By (T2), $2V_k(p_*)_k-\lambda_2(p_*)_k^2$ is maximized at $V_k/\lambda_2$ with value $V_k^2/\lambda_2$ if $V_k/\lambda_2 \le 1$: otherwise, the maximum is at $1$ with value $2V_k-\lambda_2$. Therefore 
$$c_* - \chat \le \sum_{k \not \in \hat{\mathcal J}} \varphi_{\lambda_2}(V_k).$$
The result now follows by Corollary \ref{cor:apost}.
\end{proof}


The two certificates of this document (i.e., the linear $2V$ certificate in conjunction with Corollary~\ref{cor:apost}, and the quadratic certificate of Theorem~\ref{thm:sharp}) are in fact the endpoints of a one-parameter family: the family's
minimum is itself a certificate at the same price.

\begin{theorem}[Unified dual certificate]\label{thm:dual}
For every $\mu\ge0$,
\[
c_*-\chat\;\le\;F(\mu):=\mu+\sum_{k\notin\hat{\mathcal{J}}}\varphi_{\lambda_2}\!\bigl([V_k-\mu/2]_+\bigr),
\qquad\text{hence}\qquad
c_*-\chat\;\le\; \min_{\mu\ge0}F(\mu).
\]
$F$ is convex; $F(0)$ appears in the bound of Theorem~\ref{thm:sharp} and $F(2V)=2V$ is
Corollary~\ref{cor:apost}, so the minimum is never worse than either, or the bound in Theorem \ref{thm:sharp}. It improves strictly on all of these jointly
precisely when $F'(0^{+})=1-\sum_{k\notin\hat{\mathcal{J}}}\min\{V_k/\lambda_2,1\}<0$, i.e.\ when the violations
are collectively large against the margin; the minimizer is computable in closed form after sorting
the $V_k$ (the objective is piecewise quadratic in $\mu$).
\end{theorem}

\begin{proof}
The inequality \eqref{eq:violationCertPart3} in the proof of Theorem \ref{thm:sharp} can be strengthened. Write here $u = (1_{[n] \setminus \hat{\mathcal J}})^{\hat{\mathcal{J}}}$ for the vector that has unit entries off of $\hat{\mathcal J}$ and zeros elsewhere. By (T4), we have
\begin{align}
    \max_{p_*\in\Delta_{n-1}} f(p_*) & \le \max_{\substack{p_* \in [0,1]^n \\ u^T p_* -1 \le 0}} f(p_*) \nonumber \\
    & \le \max_{p_* \in [0,1]^n} \left [ f(p_*) - \mu \left (  u^T p_* -1 \right ) \right ] \tag{****}\label{eq:violationCertPart4}
\end{align}
generically for any $f$ and $\mu \ge 0$. Now take 
$$f(p_*) = \sum_{k \not \in \hat{\mathcal J}} \left ( 2V_k(p_*)_k-\lambda_2(p_*)_k^2 \right )$$
as in the proof of Theorem \ref{thm:sharp}. 

By \eqref{eq:violationCertPart4}, we obtain
\begin{align}
    c_* - \chat & \le f(p_*) \nonumber \\ 
    & \le \max_{p_* \in [0,1]^n} \left [ f(p_*) - \mu \left (  u^T p_* -1 \right ) \right ] \nonumber \\
    & = \mu + \max_{p_* \in [0,1]^n} \left [ f(p_*) - \mu u^T p_* \right ] \nonumber \\
    & = \mu + \max_{p_* \in [0,1]^n} \left [ \sum_{k \not \in \hat{\mathcal J}} \left ( 2V_k(p_*)_k-\lambda_2(p_*)_k^2 \right ) - \mu u^T p_* \right ] \nonumber \\
    & = \mu + \max_{p_* \in [0,1]^n} \sum_{k \not \in \hat{\mathcal J}} \left ( 2V_k(p_*)_k-\lambda_2(p_*)_k^2  - \mu (p_*)_k \right ) \nonumber \\
    & = \mu + \max_{p_* \in [0,1]^n} \sum_{k \not \in \hat{\mathcal J}} \left ( (2V_k-\mu)(p_*)_k-\lambda_2(p_*)_k^2 \right ) \nonumber \\
    & = \mu + \sum_{k \not \in \hat{\mathcal J}} \max_{t \in [0,1]} \left ( (2V_k-\mu)t-\lambda_2t^2 \right ) \nonumber \\
    & = \mu + \sum_{k \not \in \hat{\mathcal J}} \varphi_{\lambda_2}([V_k - \mu/2]_+) \nonumber \\
    & = F(\mu),
    \tag{*****}\label{eq:violationCertPart5}
\end{align}
where the penultimate equality is by (T2). 

Because $F$ is a sum of nondecreasing convex functions of convex terms, it is convex. Evaluating, $F(0) = \sum_{k \not\in \hat{\mathcal J}} \varphi_{\lambda_2}(V_k)$ and $F(2V) = 2V+\sum_{k \not\in \hat{\mathcal J}} \varphi_{\lambda_2}([V_k-V]_+) = 2V+\sum_{k \not\in \hat{\mathcal J}} \varphi_{\lambda_2}(0) = 2V$, so Theorem \ref{thm:sharp} takes the form $c_* - \chat \le \min \{F(0), F(2V)\}$, whence \eqref{eq:violationCertPart5} is at least as strong. Finally, we have $F'(0^+) = 1 - \frac{1}{2}\sum_{k \not \in \hat{\mathcal J}} \varphi_{\lambda_2}' (V_k) = 1 - \frac{1}{2}\sum_{k \not \in \hat{\mathcal J}} \min \{ V_k/\lambda_2,1 \}$, and $F'((2V)^-) = 1 - \frac{1}{2}\sum_{k \not \in \hat{\mathcal J}} \varphi_{\lambda_2}' ([V_k-V]_+) = 1$. Therefore $F$ achieves an interior minimum iff $F'(0^+) < 0$, as claimed.
\end{proof}

\subsection{Impossibility: unbounded failure without metricity}\label{sec:impossible}

The envelopes above are all \emph{data-dependent}: they bound the loss in terms of quantities like
$c_{(1)}$, $V$, or $\lambda_2$ that the particular instance produces. One might hope for a clean
\emph{universal} constant. This section rules
that out for $n=4$. The construction is a thin triangle plus one far point: the
triangle's in-hull circumsphere is so flat that the far point registers as ``inside'' no matter how
distant it is, so it is deleted --- taking with it all the value it should have carried. The relative
error then grows without bound as the far point recedes. \S \ref{sec:metric} shows metricity is
exactly what prevents such behavior, which is why the phenomenon is invisible in practice.

\begin{theorem}[Unbounded relative failure at $n=4$]\label{thm:impossible}
For $h\in(0,\tfrac12]$, $\varepsilon>0$, $L\ge1$, let
$\phi^{(1)}=(0,0,0)$, $\phi^{(2)}=(1,0,0)$, $\phi^{(3)}=(\tfrac12,h,0)$, $\phi^{(4)}=(\tfrac12,-L,\varepsilon)$, and
$d_{ij}:=\lVert\phi^{(i)}-\phi^{(j)}\rVert^{2}$. Then $d$ is strict negative type for all such parameters,
and whenever
\begin{equation}\label{eq:star}
h\,(L^{2}+\varepsilon^{2}-\tfrac14)\;<\;L\,(\tfrac14-h^{2}),
\end{equation}
the recursion deletes point $4$ in round one and returns $\chat\le\tfrac12$, while
$c_*\ge(\tfrac14+L^{2}+\varepsilon^{2})/2$; hence
$c_*/\chat\ge\tfrac14+L^{2}+\varepsilon^{2}\to\infty$ along, e.g., $h=1/(8L)$, $\varepsilon\le h$. No
bound on the recursion's absolute or relative failure holds over all strict-negative-type matrices,
even at $n=4$.
\end{theorem}

\begin{proof}
\emph{Strictness.} The four points are affinely independent (the relevant determinant is
$h\varepsilon\ne0$), so by the second identity in \eqref{eq:biasvar} the form is strictly negative on
$1^\perp$: strict negative type.

\emph{Deletion of point $4$.} Work in the base plane spanned by $\{\phi^{(j)}\}_{j=1}^3$: these three points have in-hull circumcenter $z=(\tfrac12,y_c,0)$ with
$y_c=(h^2-\tfrac14)/(2h)<0$ (it sits below the base because the triangle is obtuse) and squared radius
$R^2=\tfrac14+y_c^2$. By Lemma~\ref{lem:dict}(ii), point $4$ receives a negative weight in the first
solve exactly when $\phi^{(4)}$ lies inside this circumsphere, i.e.\ when
$\lVert\phi^{(4)}-z\rVert^2-R^2<0$. A direct computation gives
$\lVert\phi^{(4)}-z\rVert^2-R^2=L^2+\varepsilon^2-\tfrac14-L(\tfrac14-h^2)/h$, which is negative precisely
under \eqref{eq:star}. So $4$ is deleted in round one (possibly alongside others) and, once gone,
never returns; hence the terminal support satisfies $\hat{\mathcal{J}}\subseteq\{1,2,3\}$.

\emph{Terminal value is small.} For any $p\in\Delta_{\{1,2,3\}}$, the identity
$q(p)=2\sum_i p_i\lVert\phi^{(i)}-z(p)\rVert^2\le2\sum_i p_i\lVert\phi^{(i)}-y\rVert^2$ holds for every point
$y$, because $\sum_i p_i \| \phi^{(i)} - y \|^2 = \| y -z(p) \|^2 + \sum_i p_i \| \phi^{(i)} - z(p) \|^2$. Choosing $y=(\tfrac12,0,0)$,
all three base points lie within distance $\tfrac12$ of $y$ (here $h\le\tfrac12$ is used), so the
right side is at most $\tfrac12$. Thus $\chat\le\tfrac12$.

\emph{True value is large.} Simply placing mass on the diametral pair $\{1,4\}$,
$c_*\ge q(\tfrac12 e_1+\tfrac12 e_4)=\tfrac12 d_{14}=\tfrac12(\tfrac14+L^2+\varepsilon^2)$.
Dividing, $c_*/\chat\ge\tfrac14+L^2+\varepsilon^2$, which $\to\infty$ as $L\to\infty$.

\emph{The parameter choice is admissible.} With $h=1/(8L)$ and $\varepsilon\le h$, the left of
\eqref{eq:star} is at most $\tfrac L8+\tfrac1{512}$ and the right at least $\tfrac{15L}{64}$, and
$\tfrac L8+\tfrac1{512}<\tfrac{15L}{64}$ for all $L\ge1$; so \eqref{eq:star} holds for $L \ge 1$, and the ratio blows up.
\end{proof}

\subsection{The metric and Euclidean case}\label{sec:metric}

The impossibility of \S \ref{sec:impossible} was built from a non-metric distance matrix (the base
triangle violates the triangle inequality). This section shows that is no accident: metricity is
precisely the condition that rules out the mechanism. The key observation is that, in the embedding,
the triangle inequality is equivalent to every embedded triangle being non-obtuse
(Lemma~\ref{lem:nonobtuse}) --- and non-obtuse triangles cannot be ``thin'' in the way the blow-up
family required. We then collect the concrete consequences: three-point metric spaces never delete at
all, the value is always at least $\dmax/2$, retaining the diameter pair guarantees a
$2$-approximation, and the floor $\dmax/2$ is attained exactly on the spaces where every point lies
between a diametral pair.

\begin{lemma}[Metricity is non-obtuseness]\label{lem:nonobtuse}
Negative-type $d$ with embedding $\phi$ satisfies the triangle inequality iff every embedded triple is
non-obtuse (all angles $\le\pi/2$): the angle at $\phi^{(j)}$ is $\le\pi/2$ iff $d_{ik}\le d_{ij}+d_{jk}$.
Consequently, for metric $d$, every triple has in-plane circumradius at most
$\operatorname{diam}(\sqrt{d})/\sqrt3$ (the largest angle is in $[\pi/3,\pi/2]$), and the thin triangles powering
Theorem~\ref{thm:impossible} cannot occur: indeed that family has
$d_{12}=1>d_{13}+d_{32}=\tfrac12+2h^{2}$.
\end{lemma}

\begin{proof}
The embedding realizes $d$ as squared distances, $d_{ab}=\lVert\phi^{(a)}-\phi^{(b)}\rVert^2$, so symmetry,
nonnegativity, and a zero diagonal are automatic and metricity is exactly the triangle inequality.
Three embedded points lie in a plane; writing $\theta_j$ for the angle of the triangle
$\phi^{(i)}\phi^{(j)}\phi^{(k)}$ at $\phi^{(j)}$, the law of cosines
$\lVert\phi^{(i)}-\phi^{(k)}\rVert^2=\lVert\phi^{(i)}-\phi^{(j)}\rVert^2+\lVert\phi^{(j)}-\phi^{(k)}\rVert^2
-2\lVert\phi^{(i)}-\phi^{(j)}\rVert\lVert\phi^{(j)}-\phi^{(k)}\rVert\cos\theta_j$ reads, in $d$-notation,
\[
d_{ik}-\bigl(d_{ij}+d_{jk}\bigr)=-2\sqrt{d_{ij}\,d_{jk}}\,\cos\theta_j .
\]
Since $\sqrt{d_{ij}d_{jk}}>0$ for distinct points, $d_{ik}\le d_{ij}+d_{jk}\iff\cos\theta_j\ge0\iff
\theta_j\le\pi/2$, which is the stated equivalence: the middle vertex of a triangle inequality is the
one whose angle it tests. The three triangle inequalities of a fixed triple (one per choice of middle
vertex) therefore hold together iff all three angles of $\phi^{(i)}\phi^{(j)}\phi^{(k)}$ are $\le\pi/2$, i.e.\ the
triple is non-obtuse; and $d$ is a metric iff every triple satisfies its triangle inequalities, hence
iff every embedded triple is non-obtuse.

For a metric triple let $\theta_{\max}$ be its largest angle and $\ell_{\max}$ the opposite (longest)
side. The angles sum to $\pi$, so $\theta_{\max}\ge\pi/3$; non-obtuseness gives $\theta_{\max}\le\pi/2$.
The law of sines $2R=\ell_{\max}/\sin\theta_{\max}$, with $\sin$ increasing on $[\pi/3,\pi/2]$ so that
$\sin\theta_{\max}\ge\sin(\pi/3)=\sqrt3/2$, gives
$R=\ell_{\max}/(2\sin\theta_{\max})\le\ell_{\max}/\sqrt3\le\operatorname{diam}(\sqrt{d})/\sqrt3$, using
$\ell_{\max}\le\operatorname{diam}(\sqrt{d})$.

Finally, the contrapositive of the equivalence turns a violated triangle inequality into an obtuse
angle at the middle vertex: the family of Theorem~\ref{thm:impossible} has
$d_{12}=1>\tfrac12+2h^{2}=d_{13}+d_{32}$ for small $h$, so its angle at $\phi^{(3)}$ exceeds $\pi/2$ and the
triple lies outside the metric cone.
\end{proof}

\begin{proposition}[$3$-point metric spaces never delete]\label{prop:n3}
For $d=\begin{psmallmatrix}0&\alpha&\beta\\ \alpha&0&\gamma\\ \beta&\gamma&0\end{psmallmatrix}$ of strict negative type,
\[
w_1=\frac{\alpha+\beta-\gamma}{2\alpha\beta},\quad w_2=\frac{\alpha-\beta+\gamma}{2\alpha\gamma},\quad w_3=\frac{-\alpha+\beta+\gamma}{2\beta\gamma},\quad
c_{(1)}=\frac{2\alpha\beta\gamma}{2(\alpha\beta+\alpha\gamma+\beta\gamma)-\alpha^{2}-\beta^{2}-\gamma^{2}} :
\]
each weight is a triangle-inequality slack over $2\times$(the incident product), so for metric $d$ the
first solve is already nonnegative and the recursion is exact: the full space if the triangle
inequalities are strict; the extreme pair if one is tight (the middle point's weight vanishes, and the
pair output satisfies \eqref{eq:kkt} with equality). The denominator of $c_{(1)}$ is
$16\operatorname{Area}(\phi^{(1)}\phi^{(2)}\phi^{(3)})^{2}$, so $c_{(1)}=2R^{2}$ as Lemma~\ref{lem:dict} predicts.
\end{proposition}

\begin{proof}
We check $dw=1$ directly with the stated $w$ (here $\det d=2\alpha\beta\gamma$), so $w$ is the weighting and its
coordinates are the displayed triangle-inequality slacks divided by $2\times$(incident product); for a
metric these slacks are nonnegative, so the first solve is already nonnegative and the recursion stops
at once, returning the full space when the inequalities are strict. If one inequality is tight, say
$\gamma=\alpha+\beta$, then $w_1=0$: the recursion drops point $1$ and returns $\{2,3\}$ with $\chat=\gamma/2$, and since
$(d\hat p)_1=(\alpha+\beta)/2=\gamma/2=\chat$, condition \eqref{eq:kkt} holds with equality, so the output is
optimal by Lemma~\ref{lem:kkt}.
\end{proof}

\begin{proposition}[Value envelope and the diameter pair]\label{prop:envelope}
For any strict negative type: $\dmax/2\le c_*<\dmax$ and $\chat\ge\tfrac12\max_{i,j\in\hat{\mathcal{J}}}d_{ij}$.
Hence
\[
\frac{c_*-\chat}{c_*}\;\le\;1-\frac{\max_{i,j\in\hat{\mathcal{J}}}d_{ij}}{2\dmax},
\]
computable at termination for free; and if the terminal set retains a pair with
$d_{ij}\ge\theta\dmax$ then $\chat\ge\theta c_*/2$. In particular, if the recursion never deletes a
maximal-distance pair, its output is unconditionally a $2$-approximation in value.
\end{proposition}

\begin{proof}
Lower bound: putting mass $\tfrac12,\tfrac12$ on a diametral pair gives
$q=\tfrac12\dmax$, so $c_*\ge\dmax/2$. Upper bound: for any $p\in\simplex$,
$q(p)=\sum_{i,j}p_ip_jd_{ij}\le\dmax\sum_{i\ne j}p_ip_j=\dmax(1-\lVert p\rVert^2)<\dmax$. For $\chat$:
by Lemma~\ref{lem:dict}(iii) the terminal value is the constrained maximum of $q$ on its face
$\Delta_{\hat{\mathcal{J}}}$, which is at least the value $\tfrac12\max_{i,j\in\hat{\mathcal{J}}}d_{ij}$ attained by placing
mass on the widest retained pair. The relative-gap display and the $\theta$-approximation follow by
dividing; in particular, if a maximal pair ($d_{ij}=\dmax$) is retained then $\chat\ge\dmax/2\ge c_*/2$.
\end{proof}

The lower end $c_*\ge\dmax/2$ of Proposition~\ref{prop:envelope} is attained, and its equality case
can be characterized completely --- with, as a bonus, a two-sided a priori estimate of $c_*$ that reads
only \emph{two rows} of $d$ and performs no solves or matvecs at all. For a diametral pair $\{a,b\}$
and any other point $k$, call $\sigma_k:=d_{ak}+d_{kb}-\dmax\ge0$ the \emph{betweenness slack} of $k$
(zero exactly when $k$ lies metrically between $a$ and $b$), and $\sigma:=\max_{k\ne a,b}\sigma_k$
(with $\sigma:=0$ if $n=2$).

\begin{proposition}[The floor of the value envelope; betweenness]\label{prop:floor}
Let $d$ be a metric of strict negative type with diametral pair $\{a,b\}$ and slacks as above. Then
\[
\frac{\dmax}2+\frac{\sigma^{2}}{4\bigl(\sigma+\dmax/2\bigr)}
\;\le\;c_*\;\le\;\frac{\dmax}2+\sigma .
\]
In particular $c_*=\dmax/2$ iff $\sigma=0$, i.e.\ iff \emph{every} point lies metrically between $a$
and $b$; and in that case 
the all-betweenness diametral pair is unique, and $p_*$ is uniform on it; $w=(e_a+e_b)/\dmax$ solves $dw=1$ \emph{exactly}, so the recursion's very first
weighting is already nonnegative (zero on every interior point) and the recursion is exact at one
solve, with $c_{(1)}=\chat=c_*=\dmax/2$;
\end{proposition}

\begin{proof}
The whole statement is the machinery of \S\S\ref{sec:dictionary} and \ref{sec:fw} read off at the
two-point start $p_{(0)}=\tfrac12(e_a+e_b)$, so we begin by computing the two quantities that control a
Frank--Wolfe step there. The value is $q(p_{(0)})=\tfrac12 d_{ab}=\dmax/2$. The marginals are
$(dp_{(0)})_k=(d_{ak}+d_{bk})/2=(\dmax+\sigma_k)/2$ for each $k$ (with $\sigma_a=\sigma_b=0$), so $\max_k(dp_{(0)})_k = (\dmax+\sigma)/2$ and the duality gap \eqref{eq:dualityGap} is
$G(p_{(0)})=2(\max_k(dp_{(0)})_k-q(p_{(0)}))=\sigma$.
The upper bound $c_*\le\dmax/2+\sigma$ is now just \eqref{eq:hPreview} applied to $p_{(0)}$. The lower bound follows from rearranging \eqref{eq:qIteration} or \eqref{eq:hOneStep} into $c_*\ge q(p_{(1)})=q(p_{(0)})+G^2/(4(G+q_{(0)}))$ and using $q_{(0)}=\dmax/2$, $G=\sigma$.
That is, the sandwich is exactly Corollary~\ref{cor:anytime} at $s=0$ started from $p_{(0)}$.

If $\sigma=0$ then $G(p_{(0)})=0$, so $p_{(0)}$ is already optimal by Lemma~\ref{lem:kkt};
if $\sigma>0$ the lower bound is strictly above $\dmax/2$. This proves $c_*=\dmax/2\iff\sigma=0$.  When every $\sigma_i=0$, the vector $w=(e_a+e_b)/\dmax$
satisfies $(dw)_i=(d_{ai}+d_{bi})/\dmax=(\dmax+\sigma_i)/\dmax=1$ for \emph{every} $i$; since $d$ is
invertible (Lemma~\ref{lem:inertia}), $w$ is the weighting $w_{[n]}$, and normalizing gives
$p^{[n]}=p_{(0)}\ge0$. So the recursion accepts at its first check with $c_{(1)}=\chat=c_*=\dmax/2$, and
$dp_{(0)}$ is constant, i.e.\ \eqref{eq:kkt} holds with equality throughout. 
Finally, two distinct all-betweenness diametral pairs would give two distinct maximizers, contradicting strict concavity, so the pair is unique and $p_*$ is uniform on it. 
\end{proof}



\subsection{The repair: one free step, and a universal rate}\label{sec:fw}

Theorem~\ref{thm:impossible} shows the \emph{uncorrected} recursion has unbounded relative failure without metricity.
The phenomenon does not survive even the cheapest repair. The repair is the \emph{Frank--Wolfe} (FW,
or conditional-gradient) method \cite{fw56,jaggi}:
given a feasible $p$, read off the marginal values 
$$g=dp$$ 
(one matvec --- the terminal certificate has
already paid for it at $p=\hat p$); find the best site $$k^{*}=\argmax_kg_k;$$ 
move from $p$ straight
toward the vertex $e_{k^{*}}$, i.e.\ shift mass toward the most undervalued site; choose how far to
move by the line-search calculus of (T2); repeat. 

Two quantities inform the analysis. First, the \emph{duality gap}
\begin{equation}\label{eq:dualityGap}
    G(p)\;:=\;2\bigl(\max\nolimits_kg_k-q(p)\bigr)\;\ge\;0
\end{equation}
is a \emph{computable overestimate of the suboptimality}: taking $p'=p_*$ in the tangent bound \eqref{eq:tangentBound},
\begin{equation}\label{eq:hPreview}
    c_*-q(p)\;\le\;2(dp)^T(p_*-p)\;=\;2\bigl(g^Tp_*-q(p)\bigr)\;\le\;2\bigl(\max\nolimits_kg_k-q(p)\bigr)=G(p),
\end{equation}
since $p_*$ is a convex combination. It costs nothing beyond the matvec already in hand, and $G=0$ certifies exact optimality
by Lemma~\ref{lem:kkt}. 

Second, the \emph{exact segment curvature}
\begin{equation}\label{eq:curv}
(e_k-p)^Td\,(e_k-p)\;=\;q(p)-2(dp)_k\;=\;-\bigl(G(p)+q(p)\bigr)
\qquad(k=\argmax\nolimits_j(dp)_j):
\end{equation}
by \eqref{eq:exactTaylor}, the restriction of $q$ to the FW segment is a parabola whose quadratic coefficient is known, not
merely bounded. So the gain of an exact-line-search step is a closed form, which yields correspondingly sharp results.

\begin{theorem}[Salvage]\label{thm:fw}
Let $d$ be strict negative type, $g=d\hat p$, $k^*=\argmax_kg_k$, and $V$ as in \eqref{eq:V}.
\begin{enumerate}
\item[(i)] \emph{(One free step.)} With $t^*=V/(\chat+2V)\in[0,1)$,
$p^{+}:=(1-t^*)\hat p+t^*e_{k^*}$ is feasible and
\[
q(p^{+})\;=\;\chat+\frac{V^{2}}{\chat+2V},
\]
computable with zero additional matvecs or solves; hence $c_*-q(p^{+})\le2V-V^{2}/(\chat+2V)$, and
$\chat+V^{2}/(\chat+2V)$ is a free \emph{lower} bound on $c_*$.
\item[(ii)] \emph{(Exact gain; universal rate.)} Run FW with exact line search from any feasible
$p_{(0)}$: at step $s$ pick $k_{(s)}=\argmax_k(dp_{(s)})_k$ and maximize $q$ on $[p_{(s)},e_{k_{(s)}}]$. Writing
$q_{(s)}=q(p_{(s)})$, $G_{(s)}=G(p_{(s)})$, the optimal step $\gamma^*_{(s)}$ is interior for $G_{(s)}>0$,
\begin{equation}\label{eq:optimalStep}
    \gamma_{(s)}^{*}=\frac{G_{(s)}}{2(G_{(s)}+q_{(s)})}\in\Bigl(0,\tfrac12\Bigr],
\end{equation}
\begin{equation}\label{eq:qIteration}
    q_{(s+1)}\;=\;q_{(s)}+\frac{G_{(s)}^{2}}{4(G_{(s)}+q_{(s)})},
\end{equation}
and consequently, with $r_{(s)}:=1-q_{(s)}/c_*\in[0,1]$ the relative failure,
\begin{equation}\label{eq:relativeFwBound}
    r_{(s)}\;\le\;\Bigl(\frac1{r_{(0)}}+\frac s4\Bigr)^{-1}\;\le\;\frac{4}{s+4},
\end{equation}
i.e.,
\begin{equation}\label{eq:absoluteFwBound}
    c_*-q(p_{(s)})\;\le\;\frac{4\,c_*}{s+4}\;\le\;\frac{4\,\dmax}{s+4},
\end{equation}
over the whole strict-negative-type class, including the family of
Theorem~\ref{thm:impossible} --- with no dimension, scale, or $\dmax$ dependence in the relative form,
and tight at $s=0$ (vertex starts have $r_{(0)}=1$). Each step costs one matrix--vector product and no
solves.
\item[(iii)] \emph{(Linear rate under strictness.)} FW with \emph{away steps} --- which may also
shift mass \emph{off} a currently supported site, along $p-e_j$, rather than only diluting the whole
distribution toward a vertex --- converges linearly:
\[
c_*-q(p_{(s)})\;\le\;\bigl(c_*-q(p_{(0)})\bigr)\,(1-\rho)^{s/2},
\qquad
\rho\;\ge\;\frac{\lambda_2}{2n\lVert d\rVert_2}\;>\;0 .
\]

This is \cite[Thm.~1]{lj15} applied to $f:=-q$, with the dictionary
\[
\begin{array}{l@{\qquad}l@{\qquad}l@{\qquad}l}
\mathcal A=\{e_j\}_{j=1}^{n}, & \mathcal M=\conv\mathcal A=\simplex, & \mu=2\lambda_2,
  & L=2\lVert d\rVert_2,\\[3pt]
M=\operatorname{diam}(\simplex)=\sqrt2, & \delta=\operatorname{PWidth}(\mathcal A)\ge2/\sqrt n,
  & \rho=\dfrac{\mu}{4L}\Bigl(\dfrac{\delta}{M}\Bigr)^{2}, & k(s)\ge s/2 .
\end{array}
\]
Four points of contact, each a place the two notations could be misread. Its ambient dimension $d$ is
our $n$ (and is not our distance matrix $d$). Its $\mu$ carries a factor $2$: by \eqref{eq:exactTaylor} the remainder is
\emph{exact}, $f(p')-f(p)-\langle\nabla f(p),p'-p\rangle=-(p'-p)^Td(p'-p)\ge\lambda_2\lVert
p'-p\rVert^{2}$ --- valid for \emph{all} $p,p'\in\simplex$, since their difference lies in
$1^{\perp}$ --- so $-q$ is $2\lambda_2$-strongly convex under the standard $\tfrac\mu2$ convention.
Its $L$ holds on all of $\R^{n}$ because $\nabla f=-2dp$ is linear, which discharges its footnote
requirement that $L$ be valid on the enlarged domain $\mathcal M+\mathcal M-\mathcal M$ for away
steps. And its $\delta$ is the width of the regular simplex, $2/\sqrt n$ for even $n$ and
$2/\sqrt{n-1/n}$ for odd $n$ \cite[App.~B.1]{lj15}. The exponent $s/2$ is its good-step count: a drop
step carries no progress guarantee but shrinks the active set, so at most half the steps are drops.
\end{enumerate}
\end{theorem}

\begin{proof}
(i) is the case $p=\hat p$ of \eqref{eq:qIteration}, where $G=2V$ and $q=\chat$:
$t^*=\gamma^{*}=V/(\chat+2V)$ and $4V^{2}/\bigl(4(2V+\chat)\bigr)=V^{2}/(\chat+2V)$. Subtracting from $c_*\le\chat+2V$ gives the result once we establish (ii), which we proceed to do.

(ii) On the segment, the exact Taylor identity \eqref{eq:exactTaylor} with $p'=p+\gamma(e_k-p)$ 
yields
\begin{align}
    q(p+\gamma[e_k-p]) & = q(p) + 2(dp)^T(\gamma[e_k-p]) + \gamma^2 (e_k-p)^T d (e_k-p) \nonumber \\
    & = q(p) + 2(g_{k^*}-q(p))\gamma + (-2g_{k^*}+q)\gamma^2 \nonumber \\
    & = q(p) + G\gamma -(G+q(p))\gamma^2
\end{align}
using \eqref{eq:curv} and the definition of $G$:
an exactly known concave parabola, which (T2) maximizes at $\gamma^{*}=G/\bigl(2(G+q)\bigr)\le\tfrac12<1$
(so line search never clips) with peak value $q+G^{2}/\bigl(4(G+q)\bigr)$, as claimed. Now write 
\begin{equation}\label{eq:h}
h_{(s)}:=c_*-q_{(s)}\le G_{(s)},    
\end{equation} 
as in \eqref{eq:hPreview}, so that $r_{(s)} = h_{(s)}/c_*$. Since $x\mapsto x^{2}/(x+a)$ is increasing for $a \ge 0$ (and we can apply this since $q \ge 0$) and
$G_{(s)}+q_{(s)}\ge h_{(s)}+q_{(s)}=c_*$,
\begin{equation}\label{eq:hOneStep}
    h_{(s+1)}\;=\;h_{(s)}-\frac{G_{(s)}^{2}}{4(G_{(s)}+q_{(s)})}\;\le\;h_{(s)}-\frac{h_{(s)}^{2}}{4(h_{(s)}+q_{(s)})}
\;=\;h_{(s)}-\frac{h_{(s)}^{2}}{4c_*}.
\end{equation}

We will need an intermediate result: that for $a,r \in (0,1)$
$$\frac{1}{r(1-ar)} \ge \frac{1}{r}+a.$$ This is true because
$\frac{1}{r(1-ar)} - \left ( \frac{1}{r}+a \right ) = \frac{1}{r} \left ( \frac{1}{1-ar} - 1 \right ) -a = \frac{1}{r} \frac{ar}{1-ar} - a = \frac{a^2r}{1-ar} \ge 0$. Now dividing the preceding inequality by $c_*$, we obtain 
$r_{(s+1)}\le r_{(s)}(1-r_{(s)}/4)$, so by the intermediate result we have
$$1/r_{(s+1)}\ge(1/r_{(s)})(1-r_{(s)}/4)^{-1}\ge1/r_{(s)}+1/4,$$ and $1/r_{(s)}\ge1/r_{(0)}+s/4$; finally
$r_{(0)}=1-q(p_{(0)})/c_*\le1$ because $q\ge0$ on the simplex ($d\ge0$ entrywise), and
$c_*<\dmax$ (Proposition~\ref{prop:envelope}) gives the absolute form. At a vertex start $q_{(0)}=0$,
$r_{(0)}=1$, and the bound is an equality at $s=0$. 

(iii) Only the tabulated constants need checking, as the formulation is engineered to precisely align with Theorem 1 of \cite{lj15}. By \eqref{eq:exactTaylor} the Taylor remainder of $q$ is exactly
$(p'-p)^Td(p'-p)$, and $p'-p\in1^{\perp}$ for $p,p'\in\simplex$, so
the definition of $\lambda_2$ gives $f(p')-f(p)-\langle\nabla f(p),p'-p\rangle\ge\lambda_2\lVert
p'-p\rVert^{2}$ for $f=-q$ throughout $\simplex$, i.e.\ $\mu=2\lambda_2$; $\nabla f=-2dp$ is linear,
hence globally $L$-Lipschitz with $L=2\lVert d\rVert_2$; $\operatorname{diam}(\simplex)=\sqrt2$ is
attained by any two vertices; and $\operatorname{PWidth}\ge2/\sqrt n$ is the simplex width
\cite[App.~B.1]{lj15}. Substituting into $\rho=(\mu/4L)(\delta/M)^{2}$ gives
$\rho\ge(\lambda_2/4\lVert d\rVert_2)(2/n)=\lambda_2/(2n\lVert d\rVert_2)$; now apply
\cite[Thm.~1]{lj15} to $f=-q$.
\end{proof}

\begin{corollary}[Anytime bracket on the peel value]\label{cor:anytime}
Along any exact-line-search FW trajectory (in particular at the recursion's terminal $\hat p$, where
$G=2V$),
\[
c_*\;\in\;\Bigl[\;q_{(s)}+\frac{G_{(s)}^{2}}{4(G_{(s)}+q_{(s)})}\;,\;\;q_{(s)}+G_{(s)}\;\Bigr]
\]
at every step, both endpoints computable from the step's own matrix--vector product: the upper end is
the duality gap, the lower end is the value the \emph{next} step is guaranteed to reach --- known
before taking it. The lower endpoints are strictly increasing; a user may stop the moment the bracket
is tight enough.
\end{corollary}


\begin{corollary}[Sparse approximate peels]\label{cor:coreset}
For every strict-negative-type $d$ and $\varepsilon\in(0,1)$ there is a feasible $p$ with
$|\supp p|\le\lceil4/\varepsilon\rceil$ and $q(p)\ge(1-\varepsilon)c_*$.
\end{corollary}

\begin{proof}
    First, note that $\frac{4}{s+4} \le \varepsilon \Rightarrow \frac{4}{\varepsilon}-4 \le s$. Now start Theorem~\ref{thm:fw}(ii) at a vertex ($r_{(0)}=1$), and stop once $4/(s+4)\le\varepsilon$. The support grows by at most one per step, so $|\supp p|\le\lceil4/\varepsilon\rceil-3$.
\end{proof}

\subsection{A shadow Frank--Wolfe iterate: an anytime bound free of the trajectory}\label{sec:shadow}

Algorithm \ref{alg:certified} is exact, but its \emph{en route} quality is not assured. Algorithm \ref{alg:shadow} below carries a second feasible \emph{shadow} point that accumulates Frank-Wolfe gains at the rate of Theorem \ref{thm:fw}(ii) and reports whichever of the two points is better. Because the shadow is never read by the peeling code, the modification is inert: the working sets, the solves, the checks and the returned answer are those of Algorithm~\ref{alg:certified} verbatim, so Theorem \ref{thm:correct} transfers directly, while the reported point acquires an \emph{a priori} convergence rate to optimality.

\begin{algorithm}[t]
\caption{Certified peeling with a shadow Frank--Wolfe iterate. Lines 1-5, 7, and 14-20 are Algorithm \ref{alg:certified} verbatim; lines 6 and 8-13 maintain the shadow $b$ and its value $q_b$, and line 21 is the early-stopping return.}\label{alg:shadow}
\begin{algorithmic}[1]
\REQUIRE strict negative type $d$ on $[n]$, $n\ge2$
\STATE $\mathcal{J}\leftarrow[n]$;\quad $p\leftarrow p^{\mathcal{J}}$
\WHILE[\textbf{Phase I}: peeling recursion]{$\min_j p_j<0$}
\STATE $\mathcal{J}\leftarrow\{j:p_j>0\}$;\quad $p\leftarrow p^{\mathcal{J}}$
\ENDWHILE
\STATE $g\leftarrow dp$;\quad $c\leftarrow p^Tg$;\quad $k\leftarrow\argmax_j g_j$
\STATE $b\leftarrow p$;\quad $q_b\leftarrow c$ \COMMENT{\textbf{shadow} initialization; no work}
\WHILE[\textbf{Phase II}: KKT certificate; $g_k\le c\Rightarrow p=p_*$]{$g_k>c$}
\IF[\textbf{shadow step}: one Frank--Wolfe step from the better point]{$c\ge q_b$}
\STATE $w\leftarrow p$;\quad $q_w\leftarrow c$;\quad $G_w\leftarrow 2(g_k-c)$;\quad $k_w\leftarrow k$ \COMMENT{free: $g$, $c$, $k$ already in hand}
\ELSE
\STATE $w \leftarrow b$;\quad$h\leftarrow db$;\quad $q_w\leftarrow b^Th$;\quad $G_w\leftarrow 2(\max_j h_j-q_w)$;\quad $k_w\leftarrow\argmax_j h_j$ \COMMENT{matvec}
\ENDIF
\STATE $\gamma\leftarrow\dfrac{G_w}{2(G_w+q_w)}$;\quad $b\leftarrow w+\gamma\,(e_{k_w}-w)$;\quad $q_b\leftarrow q_w+\dfrac{G_w^{2}}{4(G_w+q_w)}$ \COMMENT{cheap}
\STATE $\mathcal{J}\leftarrow\supp(p)\cup\{k\}$;\quad $y\leftarrow p^{\mathcal{J}}$
\COMMENT{\textbf{Phase III}: re-admit the strongest violator}
\WHILE{$\min_{j\in \mathcal{J}}y_j<0$}
\STATE $t\leftarrow\min\{\,p_i/(p_i-y_i):y_i<0\,\}$;\quad $p\leftarrow p+t\,(y-p)$
\STATE $\mathcal{J}\leftarrow\{j:p_j>0\}$;\quad $y\leftarrow p^{\mathcal{J}}$
\ENDWHILE
\STATE $p\leftarrow y$;\quad $g\leftarrow dp$;\quad $c\leftarrow p^Tg$;\quad $k\leftarrow\argmax_j g_j$
\ENDWHILE
\RETURN $p$ \COMMENT{on early stop at any check, set $p \leftarrow b$ if $q_b>c$}
\end{algorithmic}
\end{algorithm}

\begin{lemma}[The shadow is inert]\label{lem:inert}
Run Algorithm~\ref{alg:shadow} and Algorithm~\ref{alg:certified} on the same $d$, breaking ties in
$\argmax$ the same way. Then the two runs produce the same sequence of working sets $\mathcal{J}$,
the same iterates $p$, the same linear solves and matrix--vector products $dp$, the same checks, and
the same return value. Consequently Theorem~\ref{thm:correct} holds verbatim for
Algorithm~\ref{alg:shadow}: it terminates after finitely many solves and returns $p_*$ exactly, every
working set has $\lvert\mathcal{J}\rvert\ge2$, $p\in\simplex$ from the exit of Phase~I onward, each
Phase~III inner loop performs at most $n-1$ solves, and the values at successive checks strictly
increase.
\end{lemma}

\begin{proof}
The new lines 6 and 8-13 make assignments only to $b,q_b,w,q_w,G_w,k_w,\gamma,h$, none of which occurs anywhere else (i.e., in the Algorithm \ref{alg:certified} lines) except in the early-stopping comment at the end; and they read only
$p,c,g,k$, which they do not modify. Deleting them therefore changes no value taken by $\mathcal{J},p,y,t,g,c,k$, and the while-loop predicates are unchanged, so the control
flow is identical. That is the first claim, and the second is Theorem~\ref{thm:correct} applied to
the common trajectory.
\end{proof}

\begin{lemma}[One shadow step]\label{lem:shadowstep}
Throughout the course of Algorithm \ref{alg:shadow}, $w\in\simplex$, $q_w=q(w)$, and $G_w=G(w)$ as in
\eqref{eq:dualityGap}; the assignment of line 13 is the exact-line-search Frank-Wolfe step from $w$;
the new shadow satisfies $b\in\simplex$ and $q_b=q(b)$; and
\begin{equation}\label{eq:shadowgain}
q_b\;\ge\;q_w+\frac{\bigl(c_*-q_w\bigr)^{2}}{4c_*} .
\end{equation}
\end{lemma}

\begin{proof}
Feasibility and the recorded values: on the branch of line~9, $w=p\in\simplex$ with $q(p)=c$ and
$\max_jg_j=g_k$ by line~5 or 19, so $q_w=q(w)$ and $G_w=2(\max_j(dw)_j-q(w))=G(w)$; on the branch of
line~11 the same holds by direct computation from $h=db$, and $b\in\simplex$ by induction (below).
By Theorem~\ref{thm:fw}(ii) applied at $w$, the maximizer of $q$ on $[w,e_{k_w}]$ is attained at
$\gamma^{*}=G_w/\bigl(2(G_w+q_w)\bigr)\in\bigl(0,\tfrac12\bigr]$. This is $\gamma$ of line~12, so
the step is never clipped and $b=w+\gamma(e_{k_w}-w)$ is a convex combination of points of $\simplex$,
hence feasible. The corresponding value is $q_w+G_w^{2}/\bigl(4(G_w+q_w)\bigr)$, which is the recorded $q_b$.
(At line~12 the loop predicate gives $g_k>c$, so $G_w>0$ on the first branch; on the second,
$G_w\ge c_*-q_w>0$ by \eqref{eq:hPreview} unless $b$ is already optimal, in which case $\gamma=0$ and
the step is vacuous.) For \eqref{eq:shadowgain}, $x\mapsto x^{2}/\bigl(4(x+q_w)\bigr)$ is increasing on
$x\ge 0$ and $G_w\ge c_*-q_w$ by \eqref{eq:hPreview}, so
$G_w^{2}/\bigl(4(G_w+q_w)\bigr)\ge(c_*-q_w)^{2}/\bigl(4((c_*-q_w)+q_w)\bigr)=(c_*-q_w)^{2}/(4c_*)$.
\end{proof}

\begin{theorem}[Anytime bound, free of the trajectory]\label{thm:shadowrate}
Index the Phase~II checks of Algorithm~\ref{alg:shadow} by $i=0,1,\dots$, let $c_{(i)}$ be the value
at check $i$, let $\pi_{(i)} = b$ if $q_b > c$ and $\pi_{(i)} = p$ otherwise (i.e, $\pi_{(i)}$ is the point that would be returned on an early stop), and put
$r_{(0)}:=1-c_{(0)}/c_*$. Then for every $i$,
\begin{equation}\label{eq:shadowrate}
1-\frac{q(\pi_{(i)})}{c_*}\;\le\;\Bigl(\frac{1}{r_{(0)}}+\frac{i}{4}\Bigr)^{-1}\;\le\;\frac{4}{i+4},
\qquad\text{equivalently}\qquad
c_*-q(\pi_{(i)})\;\le\;\frac{4c_*}{i+4}\;\le\;\frac{4\dmax}{i+4} .
\end{equation}
This holds for every strict-negative-type $d$ and every trajectory --- however many Phase~III inner
loops truncate, and with no dependence on $n$, on $\dmax$, or on scale in the relative form. Moreover
$q(\pi_{(i)})\ge c_{(i)}$ at every check, and $c_*\in[\,q(\pi_{(i)}),\,q_w+G_w]$ with both endpoints
computed, so a run may be stopped at any check with a bound on how good its answer is.
\end{theorem}

\begin{proof}
First note that $q(\pi_{(i)})=\max\{q_b,c_{(i)}\}$ is the value at check $i$, before that check's shadow step. Write here $\hat r_{(i)}:=1-q(\pi_{(i)})/c_*\in[0,1]$: feasibility of both candidates is Lemma~\ref{lem:inert} for $p$ and Lemma~\ref{lem:shadowstep} for $b$, and $\hat r_{(i)}\ge0$ since
$c_*$ is the maximum. At check $i$, line~9 or 11 selects $w$ with $q_w=q(\pi_{(i)})$, so
\eqref{eq:shadowgain} gives, for the shadow entering check $i+1$,
\[
q_b\;\ge\;q(\pi_{(i)})+\frac{c_*\,\hat r_{(i)}^{2}}{4}\ ,\qquad\text{i.e.}\qquad
1-\frac{q_b}{c_*}\;\le\;\hat r_{(i)}\Bigl(1-\frac{\hat r_{(i)}}{4}\Bigr).
\]
Since $q(\pi_{(i+1)})\ge q_b$, also $\hat r_{(i+1)}\le\hat r_{(i)}\bigl(1-\hat r_{(i)}/4\bigr)$. If
$\hat r_{(i)}=0$ the sequence is $0$ thereafter and \eqref{eq:shadowrate} is trivial; otherwise
$\hat r_{(i)}/4<1$ and the elementary inequality $\frac{1}{r(1-ar)}\ge\frac1r+a$ for $a,r\in(0,1)$
(proved in Theorem~\ref{thm:fw}(ii)) gives $1/\hat r_{(i+1)}\ge1/\hat r_{(i)}+\tfrac14$. Induction
yields $1/\hat r_{(i)}\ge1/\hat r_{(0)}+i/4$, and $\hat r_{(0)}=r_{(0)}\le1$ because $b=p$ at line~6
and $q\ge0$ on $\simplex$; $c_*\le\dmax$ gives the absolute form. Finally $q(\pi_{(i)})\ge c_{(i)}$ by
definition of $\pi_{(i)}$, and the bracket is $q(\pi_{(i)})\le c_*\le q_w+G_w$, the upper end being
\eqref{eq:hPreview} at $w$.
\end{proof}

\begin{remark}[Cost]\label{rem:shadowcost}
Per check the shadow costs $O(n)$ scalar work plus, only on the branch of line~11, one product $db$ at
$O(n\lvert\supp b\rvert)$ with $\lvert\supp b\rvert$ growing by at most one per step. No linear solve
is performed and no factorization is formed or updated, so a run with $N$ re-entries pays at most (and typically much less than) $N$
extra matrix--vector products against the $N+1$ products and $\ge N$ solves it already performs. 
\end{remark}

\begin{remark}[The nature of the shadow point]\label{rem:shadowscope}
Theorem~\ref{thm:shadowrate} is a statement about the \emph{reported} point $\pi_{(i)}$, not about the
peeling iterate $p_{(i)}$. The theorem does not assert $1-c_{(i)}/c_*\le4/(i+4)$, and no such bound is proved
here for Algorithm~\ref{alg:certified}. The theorem asserts that the two guarantees one wants can
be had together and at negligible cost with the peeling trajectory left exactly as it was,
so that every statement proved about that trajectory elsewhere continues to hold unaltered. Exactness and finite termination come from
Theorem~\ref{thm:correct} via Lemma~\ref{lem:inert}, and a dimension-, scale- and trajectory-free
anytime rate comes from Theorem~\ref{thm:shadowrate}.
\end{remark}

\subsection{Peeling more efficiently}\label{sec:sketch}

In future work we will explore low-rank approximations of $d$ and matvecs $p\mapsto dp$ to accelerate calculations of peels and quadratic entropies.

\end{document}